\documentclass[11pt]{amsart}
\usepackage{amsmath}
\usepackage{amssymb}
\usepackage{latexsym}
\usepackage{amscd}
\usepackage{enumerate}
\usepackage[all]{xy}
\input xy
\xyoption{all}
\usepackage[bookmarks=false]{hyperref}
\numberwithin{equation}{section}

\newcommand{\A}{{\cal A}}

\def\j{\mathbf{j}}

\def\N{\mathcal{N}}
\def\X{\mathcal{X}}

\def\proof{\noindent\emph{Proof:}\hs}

\newcommand{\om}{\omega}

\newcommand{\tens}{\otimes}
\newcommand{\Ham}{\text{Ham}}

\def\ov{\overline}
\def\Z{\mathbb{Z}}
\def\Q{\mathbb{Q}}
\def\R{\mathbb{R}}
\def\C{\mathbb{C}}
\def\G{PSL_2(\C)}
\def\S{\mathcal{S}}
\def\J{\mathcal{J}}

\def\H{\mathcal{H}}
\def\F{\mathcal{F}}
\def\B{\mathcal{B}}
\def\O{\mathcal{O}}

\def\S{\mathcal{S}}
\def\F{\mathcal{F}}
\def\H{\mathcal{H}}
\def\L{\mathcal{L}}
\def\C{\mathbb{C}}
\def\E{\mathcal{E}}
\def\A{\mathcal{A}}
\def\M{\mathcal{M}}
\def\io{\iota}
\def\hra{\hookrightarrow}

\def\la{\langle}
\def\wt{\widetilde}
\def\ra{\rangle}
\def\LA{\left\langle}
\def\RA{\right\rangle}
\def\({\left(}
\def\){\right)}
\def\/{{/\!\!/}}
\def\sra{\rightarrow}

\def\CP{{\mathbb{CP}}}
\def\CC{{\mathcal{C}}}

\def\MMC{{\mathcal{M}}}
\def\coker{{\text{coker }}}

\def\ind{{\text{ind }}}

\def\into{{\hookrightarrow}}

\def\lra{\longrightarrow}

\def\delbar{\ov{\partial}}
\def\hs{\hspace{0.1cm}}
\def\hsp{\hspace{0.2cm}}
\def\hspa{\hspace{0.3cm}}
\def\hspac{\hspace{0.5cm}}

\def\and{\hsp\text{and}\hsp}
\def\ov{\overline}

\def\cal{\mathcal}
\def\hra{\hookrightarrow}

\def\aut{\text{Aut}}
\def\t{\textbf{\text{t}}}
\def\m{\textbf{\text{m}}}

\def\D{\mathcal{D}}
\def\and{\hspa\text{and}\hspa}
\def\prgl{\text{pgl}}

\def\proof{\noindent\emph{Proof:}\,\,\,\,\,}
\def\bs{\backslash}
\def\j{\mathbf{j}}

\newtheorem{defn}{Definition}[section]
\newtheorem{theorem}{Theorem}[section]
\newtheorem*{implicitfunctiontheorem}{Implicit Function Theorem}
\newtheorem*{transversality}{Theorem B}
\newtheorem*{productformula}{Theorem A}
\newtheorem*{structurethm}{Theorem B'}
\newtheorem*{csplitting}{Corollary}
\newtheorem*{fibrationmodule}{Theorem C}
\newtheorem*{theorem*}{Theorem}

\newtheorem{assumption}[theorem]{Assumption}

\newtheorem{ex}[theorem]{Example}

\newtheorem{lem}[theorem]{Lemma}
\newtheorem{prop}[theorem]{Proposition}
\theoremstyle{definition}
\newtheorem{rem}[theorem]{Remark}  

\def\begineq{\begin{equation}}
\def\endeq{\end{equation}}
\def\beginit{\begin{itemize}}
\def\endit{\end{itemize}}
\def\begineqn{\begin{eqnarray}}
\def\endeqn{\end{eqnarray}}

\def\dti{\left.\frac{d}{dt}\right|_{t=0}}

\begin{document}

\title{A product formula for  Gromov-Witten invariants.}

\author{Cl\'ement Hyvrier}
\address{
Universit\'e de Montpellier 2\\
Institut de Math\'ematiques et de Mod\'elisation de Montpellier, I3M\\
Place Eug\`ene Bataillon\\
34095 Montpellier\\
France
}
\email{
chyvrier@math.univ-montp2.fr
}

%\begin{itemize}
%\item For a given symplectic manifold $(B,\om_B)$, and an almost complex structure $J_B$ that tames $\om_B$, the $GW$-invariants are given by counting the (algebraic) number of unparametrized  $J_B$-holomorphic genus 0 maps with $l$ marked points, representing a given spherical homology class  $\sigma\in H_2(B,\Z)$, and intersecting transversally $l$ given cycles of $B$ at the marked points. More precisely, these invariants arise as an intersection number between
%\begin{itemize}
%\item[1)] the moduli space $\M_{0,l}(B,\sigma,J_B)$ of unparametrized $J_B$-holomorphic maps representing $\sigma$ with $l$ marked points, i.e the set
%$$\{(u,x_1,...,x_l)\in C^{\infty}(S^2,B)\times (S^2)^l | \delbar_{J_B}u=0, x_i\neq x_j, [u(S^2)]=\sigma\},$$
%\item[2)]and a product of pseudocycles  $f_i:V_i\sra B$, $i=1,...,l$, of $B$, that represent homology classes $c_i^B\in H_*(B)$ (Here $H_*(B)$ stands for the homology of $B$ modulo torsion).
%\end{itemize}
%%\item Two coordinate charts coming from different chart datas are  in fact $C^{\infty}$ compatible, since we consider only smooth elements  \cite{CL} ?
%\end{itemize}

\begin{abstract} We  establish a product formula for Gromov-Witten invariants for closed relatively semi-positive Hamiltonian fibrations, with connected fiber, and  over any connected symplectic base. Furthermore, we show that the fibration projection induces a locally trivial (orbi-)fibration map from the moduli space of pseudo-holomorphic maps with marked points in the total space of the Hamiltonian fibration to the corresponding moduli space of  pseudo-holomorphic maps with marked points in the base.  We use this induced map to recover the product formula by means of integration. Finally, we give  applications to $c$-splitting and symplectic uniruledness.
\end{abstract}

\maketitle

%\tableofcontents

\section{Introduction} 
%Let  $(F,\om)$ be a  symplectic manifolds, and let $J$ be a tame almost complex structure on $F$. For $\sigma\in H_2(F,\Z)$, let 
%$$\M_{0,l}(F,\sigma,J):=\{(u,x_1,...,x_l)\in C^{\infty}(S^2,F)\times (S^2)^l | \delbar_Ju=0, x_i\neq x_j, [u(S^2)]=\sigma\}$$
%where 
% the moduli space of genus zero  
 We consider rational Gromov-Witten invariants ($GW$-invariants) of closed  Hamiltonian fibrations with connected fiber over any  connected symplectic base.  
In  a symplectic manifold $(X,\om)$ with $\om$-tame almost complex structure $J$,   $GW$-invariants are given by counting the (algebraic) number of 
%$$(u:S^2\sra X,\mathbf{x}:=(x_1,...,x_l))\in C^{\infty}(S),$$
of unparametrized  $J$-(pseudo)-holomorphic genus 0 maps with $l$ distinct marked points,     representing a fixed spherical homology class  $A\in H_2(X,\Z)$,  and intersecting transversally $l$ given cycles of $X$ at the marked points. Roughly speaking, if $\M_{0,l}(X,A,J)$ denotes  the moduli space of unparametrized (genus 0) $J$-holomorphic maps with $l$ markings, $(u,x_1,...,x_l)$, representing the class $A$, and if $M_1,...,M_l$ are cycles in $X$ representing given classes $c_1,...,c_l\in H_*(X,\Q)$, then 
these invariants  can be seen as  the values of the multilinear homomorphism:
\begin{eqnarray*}
\la\phantom\cdots \ra^X_{0,l,A}:(H_*(X,\Q))^{\otimes l}&\sra &\Q \\
c_1\otimes...\otimes c_l&\mapsto& ev^X_l\cdot(M_1\times...\times M_l)
\end{eqnarray*}
the intersection pairing $\cdot$ being taken with respect to the evaluation at the marked points:
$$ev^X_{l}:\M_{0,l}(X,A,J)\sra X^l,\hsp (u,x_1,...,x_l)\mapsto (u(x_1),...,u(x_n)).$$
It is well known  \cite{MS}, \cite{M}, \cite{RT}, that
when the symplectic manifold  is semi-positive the  $GW$-invariants are generically  well-defined  $\Z$-valued invariants of the symplectic manifold defined on $H_*(X)$, the singular homology of $X$ modulo torsion.

In the case where the Hamiltonian fibration is a product of two symplectic manifolds, Ruan and Tian ~\cite{RT},  and shortly after Kontsevich and Manin  \cite{KM2}, showed that when both the base and the fiber are semi-positive, the $GW$-invariants of the total space are given  by products of  related $GW$-invariants in the base and in the fiber, giving rise to splitting of the quantum product. A priori, we cannot expect for such a splitting to hold in the non-trivial case, as even the cup product may not split already. Nevertheless, one may still ask about  the algebraic relations that can be established out of the invariants of the base, the fiber, and the total space.   We give such a relation under some assumptions, in particular  when the reference fiber of the Hamiltonian fibration is semi-positive relative to the total space.  We now explain in  details the setting, our main results, and some consequences.
 \\

\noindent\textbf{Hamiltonian fibrations.} By definition, a symplectic  fibration is a  smooth locally trivial fibration $\pi: P\sra B$  
with symplectic reference fiber $(F,\omega)$,  and which structure  group lies in the group of symplectic  diffeomorphisms of the fiber, denoted $\mathrm{Symp}(F,\omega)$. It follows that each  fiber  $F_b:=\pi^{-1}(b)$  is naturally  equipped with a symplectic form $\omega_b$. A symplectic fibration is  Hamiltonian if the structure group can be reduced to the group $\mathrm{Ham}(F,\omega)$ of Hamiltonian diffeomorphisms. By a result of Guillemin-Lerman-Sternberg \cite{GLS}, this is equivalent to satisfying the  following two conditions:
\begin{itemize}
\item[($H_1$):]  \emph{$P$ is symplectically trivial over the $1$-skeleton, $B_1$, of $B$;}
\item[($H_2$):] \emph{there exists a connection $Hor_{\tau}\subset TP$  with holonomy in $\mathrm{Ham}(F,\omega)$,  induced by a (canonical) closed $2$-form $\tau\in \Omega^2(P)$ extending the family $\{\omega_b\}_{b\in B}$.}
\end{itemize}

 The closed $2$-form is usually refered to as the \emph{coupling form}. Since $B$ is assumed to be closed and symplectic, the existence of such a form  is sufficient to give $P$ a symplectic structure compatible with the family $\{\omega_b\}_{b\in B}$  by considering the form $$\omega_{P,\kappa}:=\tau+\kappa \pi^*\omega_B,$$ where $\kappa>0$ is a real number chosen large enough so that $\omega_{P,\kappa}$ is non-degenerate.  A simple example of such fibration is the trivial product of the base with the fiber. A less trivial class of examples is that of Hamiltonian fibrations over $S^2$. These latter fibrations correspond  (up to isomorphism) to homotopy classes  of loops in  $\mathrm{Ham}(F,\omega)$. This direct connection with the fundamental group of  $\mathrm{Ham}(F,\omega)$ makes these objects particularly interesting  from the point of view of symplectic topology as pointed out by Seidel  \cite{Se}.  \\

\noindent\textbf{Product formula.} When the minimal Chern number of the fiber, $N_F$, satisfies the following \emph{semi-positivity relative to $P$},
 \begin{equation}\tag{$\star$}\label{ssp} N_F\geq\frac{1}{2}\dim P-2, 
 \end{equation}
we give  a \emph{product formula} relating  the rational $GW$-invariants of the base  with the $GW$-invariants  of the total space, as suggested in  \cite{LM}. In order to do so we equip the fibration $P$ with an \emph{(almost) complex structure} $J_P$, chosen compatibly with the fibration structure and a Hamiltonian connection.
Concretely, $J_P$ is uniquely given by the choice of coupling form $\tau$, an $\om_B$-tame complex structure $J_B$ on $B$ and a family of  $\om_b$-tame almost complex structures $J_b$ in $F_b$.
%More precisely we ask that $J_P$ projects to an  $\omega_B$-tame almost complex structure $J_B$ via $d\pi$, that it preserves the horizontal distribution  $\textrm{Hor}_{\tau}$ (for a given $\tau$), and such that its restriction to any fiber $F_b$ yields an  $\omega_b$-tame almost complex structure.  
Such structures are said to be \emph{compatible with $\pi$ and $\tau$} or just \emph{fibered}. 
%Let $\M_{0,l}(P,\sigma,J_P)$ be the moduli space of  unparametrized pairs $(u,\mathbf{x})$,  where $u:S^2\sra P$ is a (smooth) $J_P$-holomorphic map  representing $\sigma\in H_2(P,\Z)$, and $\mathbf{x}$ is the data of $l$ pairwise distinct marked points on $S^2$, and let $\MMC_{0,l}(B,\pi_*\sigma,J_B)$ be defined similarly.   The projection $\pi$ naturally induces a map 
% $$\pi:\MMC_{0,l}(P,\sigma,J_P)\lra\MMC_{0,l}(B,\pi_*\sigma,J_B).$$
%In general,  the genus 0  $GW$-invariants are given by the values of  a
% multilinear homomorphism:
% $$\la\phantom\cdots \ra^P_{0,l,\sigma}:(H_*(P,\Q))^{\otimes l}\lra \Q.$$
%Roughly speaking, for cycles $c_1^P,...,c_l^P$, the corresponding value is  given by 
%%Furthermore if $ev:\M_{0,l}(P,\sigma,J_P)\sra P^l$ stands for the natural evaluation at the marked points $(u,\mathbf{x})\mapsto u(\mathbf{x})$, then 
%%$$\la c_1^P,...,c_l^P\ra^P_{0,l,\sigma}:=\int_{[\M_{0,l}(P,\sigma,J_P)]} ev^*(PD(c_1^P)\wedge...\wedge PD(c_l^P))$$
%  the algebraic number of rational $J_P$-holomorphic spheres  representing  $\sigma$ and intersecting  transversally the  given cycles  at the marked points. It is well known \cite{MS}, \cite{M}, \cite{RT}, that
%when the manifold is semi-positive the  $GW$-invariants are  well-defined  $\Z$-valued invariants of the symplectic manifold, defined on $H_*(P)$ the singular homology of $P$ modulo torsion.
In the present formula  we consider classes $c_i^P\in H_*(P)$, $i=1,...,l$, that are given by product classes $c_i^B\tens c_i^F$, with  $c^B_i\in H_*(B)$ and  $c^F_i\in H_*(F)$, such that the following condition is verified for some integer $0\leq m\leq l$:
\begin{equation}\tag{$\star\star$}\label{cond11}
\begin{cases} c_i^B=pt,   & \text{for $i=1,...,m$} \\
 c_i^F=[F] & \text{for $i=m+1,...,l$.}
\end{cases}
\end{equation}

\noindent Concretely, the classes in $P$ are such that (up to multiplication by well chosen integers) their Poincar\'e duals can be  represented either by a submanifold of the fiber or by the preimage under $\pi$ of a submanifold in $B$. Let  $P_C$ denote the restriction of $P$ along the image $C$ of a smooth map from $S^2$ to $B$. This defines a Hamiltonian fibration over $S^2$ with   coupling form given by the  pull-back of  $\tau$ under the natural inclusion $\iota^P_{P_C}:P_C\hra P$. Set, $$B_{\sigma}:=\{\sigma'\in H_2(P_C,\Z)|\iota^P_{P_C}(\sigma')=\sigma\},$$
where $\iota^P_{P_C}$ is understood as  the induced map in homology, and let  $\iota_F^{P_C}:F \hra P_C$ denote  both the natural inclusion and its induced map in homology.  The product formula is as follows:

\begin{productformula}\label{productformula0} 
Let $\pi: P\sra B$ be a Hamiltonian fibration  with semi-positive fiber $(F,\om)$ relatively to $P$. Let $\sigma\in H_2(P,\Z)$ and 
suppose  $\sigma_B:=\pi_*(\sigma)\neq 0$ only admits irreducible effective decompositions for some $J_B$. 
Define $c_i^P,c_i^B,c_i^F$ as in \eqref{cond11}. 
Then for a generic fibered complex structure the following equation holds:
$$ \la c^P_1,...,c^P_l\ra_{0,l,\sigma}^{P}=\la c_1^B,...,c_l^B\ra^B_{0,l,\sigma_B}\cdot  \sum_{\sigma'\in B_{\sigma}}\LA\iota^{P_C}_{F}(c^F_1),...,\iota^{P_C}_{F}(c^F_l)\RA^{P_C}_{0,l,\sigma'}$$
where $C$ is a curve counted in $\la c_1^B,...,c_l^B\ra^B_{0,l,\sigma_B}$.
 \end{productformula}

This formula states that the number $\la\phantom\cdots\ra^{P_C}_{0,l}$ of curves above $C$  is independent of  the chosen $C$, which turns out to be essentially a consequence of ($H_1$). 
%From the divisor axiom applied to the fibres $F$ in $\PC$, we can actually rewrite the product formula as:
%$$ \la c^P_1,...,c^P_l\ra_{0,l,\sigma}^{P}=\la c_1^B,...,c_l^B\ra^B_{0,l,\sigma_B}\cdot  \sum_j \la\iota_C(c^F_1),...,\iota_C(c^F_m)\ra^{\left.P\right|_C}_{0,m,\sigma_j}.$$ 
Note that by Gromov's compactness the above sum is finite.  It is even possible to simplify the expression of the formula by considering \emph{equivalence classes} on the  preimage $\pi_*^{-1}(\sigma_B)\subset H_2(P,\Z)$. More precisely, we say that  $\sigma_1,\sigma_2\in \pi_*^{-1}(\sigma_B)$ are equivalent  if and only if
$$\tau(\sigma_1-\sigma_2)=0=c^v(\sigma_1-\sigma_2),$$
where $c^v$ denotes the first Chern class of the vertical subbundle $\ker d\pi\subset TP$.
 Let  $[\sigma]_{\sigma_B}$ denote the equivalence class of $\sigma$ in the product formula. Note that, under the pull-back by $\iota_{P_C}^P$,  any element in $\pi_*^{-1}(\sigma_B)$ defines a \emph{section class} in $P_C$ (i.e a class projecting on  $[S^2]$ under $\pi_*$), and  it is not hard to see that the preimage of $[\sigma]_{\sigma_B}$ under $\iota^{P}_{P_C}$ gives rise to an equivalence class of section classes in $P_C$. If $\sigma_C$ denotes this equivalence class, then the sum in the product formula disappears:
\begin{equation}\label{simpversionPF}\la c^P_1,...,c_l^P\ra_{0,l,[\sigma]_{\sigma_B}}^P=\la c^B_1,...,c_l^B\ra_{0,l,\sigma_B}^B\cdot \la \iota_F^{P_C}(c^F_1),...,\iota_F^{P_C}(c^F_l)\ra_{0,l,\sigma_C}^{P_C}.
\end{equation}
Finally,  we remark that the case $\sigma_B=0$ leads to  the \emph{Parametric Gromov-Witten invariants} which are well known \cite{B}, \cite{LeO}. 

An important issue in the proof of this theorem is  establishing that the  $GW$-invariants  involved are generically and simultanuously well-defined.  The irreducibilty hypothesis on $\sigma_B$, which is realized by primitive classes, arises for transversality purposes. We make this more precise. In the present context,  the projection $\pi$ naturally induces a map 
 $$ \ov{\pi}:\MMC_{0,l}(P,\sigma,J_P)\sra\MMC_{0,l}(B,\sigma_B,J_B),\hspace{0.3cm} \sigma_B:=\pi_*\sigma.$$
The issue is to realize generically both moduli spaces as smooth oriented manifolds together with preserving $\ov{\pi}$. In the general context of  a symplectic manifold $(X,\om)$,  it is well-known  that the \emph{irreducible component}, $\M^*_{0,l}(X,A,J)$ of  $\M_{0,l}(X,A,J)$, consisting  of \emph{simple} maps, i.e maps $u$ that are \emph{somewhere injective} (there exists a point $z_0\in S^2$ such that $du(z_0)$ is injective and $u^{-1}(u(z_0))=z_0$),  is actually  an oriented open manifold of finite dimension for a generic choice of $\om$-tame structure $J$  (see \cite{MS}, \cite{R}). Here we cannot  directly  apply this  result since  tame  complex structures of the  total space  do not coincide with fibered complex structures. 
Nevertheless,  generalizing a  result of McDuff and Salamon in the case of Hamiltonian fibrations over Riemann surfaces \cite{MS} we prove that:

\begin{transversality} 
Suppose $\sigma_B\neq0$. There exists a second category subset, $\J_{P,reg}$, of fibered almost complex structures such that for every $J_P\in\J_{P,reg}$
\begin{enumerate}[1)]
\item the subset $\MMC^{**}_{0,l}(P,\sigma,J_P)$ of $\MMC_{0,l}(P,\sigma,J_P)$ consisting of simple maps that project to simple maps under $\ov{\pi}$, and the moduli space $\MMC^*_{0,l}(B,\sigma_B,J_B)$,
are open oriented manifolds. 
\item
for any countable set $Z$ of elements  in $\MMC^*_{0,l}(B,\sigma_B,J_B)$, for every $u\in Z$,  the preimage $\ov{\pi}^{-1}(u)$ is an open oriented manifold.
\end{enumerate}
\end{transversality}

The dimensions of the manifolds in 1) are respectively given by the indices of  the linearizations,  $D^P$ and $D^B$, of the Cauchy Riemann operators  $\delbar_{J_P}$ and $\delbar_{J_B}$, while the dimension of the moduli space in 2) is given by the index of $D^v$, the restriction of $D^P$ to vector fields along the curves that are vertically valued. The proof of  Theorem B is based on the relation
$$\pi_*\circ D^P=D^B\circ \pi_*.$$
Consequently,  $\pi$ induces a \emph{submersion of Fredholm systems} (see Section 2) between the \emph{Fredholm systems} relative to the operators $\delbar_{J_P}$ and $\delbar_{J_B}$, as defined in \cite{CL}. From this we derive  an   exact sequence
%The  moduli spaces correspond to the sets of solutions of the associated Cauchy-Riemann operators with respect to $\delbar_{J_P}$ and $\delbar_{J_B}$.
%These first order differential operators can,  be seen as sections of some Banach fibrations over 
%$$\E_P(\sigma,J_P)\lra \B_P(\sigma)\and \E_B(\sigma_B,J_B)\lra \B_B(\sigma_B),$$
% appropriate completions  of $C^{\infty}(S^2,P)$ and $C^{\infty}(S^2,B)$.
%These fiber bundles, together with the above sections, define 
%$$ (\B_P(\sigma),\E_P(\sigma,J_P),\delbar_{J_P}^H)\and  (\B_B(\sigma_B),\E_B(\sigma_B,J_B),\delbar_{J_B})$$
% \emph{Fredholm systems} as defined in \cite{CT} (if we forget about the compacity condition on the sets of solutions). In the fibration context the projection $\pi$ induces a \emph{submersion of Fredholm systems}, as described in section 2. The  main observation made here is that the linearizations $D$ of $\delbar_{J_P}$ and $D^B$ of $\delbar_{J_B}$ satisfy the relation:
%$$\pi_*\circ D=D^B\circ \pi_*.$$
\begin{equation*}
  0\rightarrow \ker D^{v} \rightarrow \ker D^P \rightarrow \ker D^B\rightarrow \coker D^{v}\rightarrow \coker D^P\rightarrow\coker D^B\rightarrow 0.
\end{equation*}
%\begin{transversality}
%Suppose $\sigma_B\neq0$. There exists a second category Baire set, $\P_{reg}\subset\P$, consisting of fibered almost complex structures such that:
%\begin{enumerate}[{\rm 1)}]
%\item for every $J_P=(J_B,J,H)\in\P_{reg}$, the moduli spaces:
%$$\MMC^{**}_{0,l}(P,\sigma,J_P)\and \MMC^*_{0,l}(B,\sigma_B,J_B),$$
%are open oriented manifolds of respective dimensions $\text{Ind}(D)$ and $\text{Ind}(D^B)$,  and where $\MMC^{**}_{0,l}(P,\sigma,J_P)$ stands for the restriction to \emph{simple maps} that project to simple maps. 
%\item
%For any countable set $A$ of elements  in $\MMC^*_{0,l}(B,\sigma_B,J_B)$, for every $u\in A$  the preimage $\pi^{-1}(u)$ is an open oriented manifold of dimension the index of $D^v$.
%\end{enumerate}
%\end{transversality}
The result then follows by ensuring that the cokernels vanish, at least at the level of the universal moduli spaces, for example the irreducibility hypothesis on $\sigma_B$ ensures that the last term of the sequence vanishes. It is worth pointing out  that without this assumption standard transversality may  fail a priori, due to multiple coverings as shown in \cite{M}.  The remaining obstructions are dealt with by perturbing the Hamiltonian connection as in \cite{MS}.  
More generally, $\ov{\pi}$ extends to a map, still denoted $\ov{\pi}$, between the compactifications $\ov{\MMC}_{0,l}(P,\sigma,J_P)$ and $\ov{\MMC}_{0,l}(B,\sigma_B,J_B)$ of  the moduli spaces.
These compactifications are stratified spaces for which the strata can be represented  by \emph{stable stratum data}  $\S$, as pointed out by Kontsevitch, that is, by connected trees with tails together with an effective decomposition of the represented second homology class. We can repeat the arguments above for each stratum $\M_{\S_P}(P)$ mapped to a stratum $\M_{\S_B}(B)$ under $\ov{\pi}$, in order to show that transversality is generically realized for the irreducible elements in $\M_{\S_P}(P)$ whenever $\M_{\S_B}(B)$ does not contain reducible elements. Then condition \eqref{ssp} ensures  that the ``boundary" of the compactified moduli spaces above, given by lower strata,  have  codimension at least 2 with respect to the top stratum consisting of simple maps. 

\begin{rem} At this point, it is worth mentioning  that the restriction to the genus 0 case is not essential. Although we have not treated  the case of higher genus curves, all  the results should still go through  with minor modifications, except regarding  the applications to $c$-splitting and symplectic uniruledness below.

Another  noteworthy observation, is that the restrictions on $\sigma_B$ and the relative semipositivity conditions are only of  technical order. It is believed that those \emph{ad hoc} hypothesis can be avoided by using virtual perturbations (see \cite{CL}, \cite{LT}, \cite{Rvirt}, \cite{Z}), which have been developped in order to deal with transversality issues for general symplectic manifolds. Removing these assumptions, is part of a joint work in progress with Shengda Hu.  
\end{rem}

%  url = {http://www.citebase.org/abstract?id=oai:arXiv.org:math/0610369},
%  url = {http://www.citebase.org/abstract?id=oai:arXiv.org:0704.3899},
% url = {http://www.citebase.org/abstract?id=oai:arXiv.org:math/0610370},

\noindent\textbf{Fibration structure.} It is natural to ask about the structure of the map $$\ov{\pi}:\ov{\M}_{0,l}(P,\sigma,J_P)\sra \ov{\M}_{0,l}(B,\sigma_B,J_B).$$
In particular, it would be interesting to understand when $\ov{\pi}$ is a fibration, at least above the top stratum of the target space. When this is the case, we can recover the  product formula using integration over the  fibers of $\ov{\pi}$ (see Section 5).
Assuming the linearized operators involved in the exact sequence above are all surjective, it  follows that the restriction of $\ov{\pi}$ to $\M_{0,l}(P,\sigma,J_P)$ is a smooth submersion onto $\M_{0,l}(B,\sigma_B,J_B)$. However, this map is not proper. This latter condition is important, as one can easily construct a smooth submersion which is not proper and which does not induce a fibration structure. To solve this problem,  we consider the fiberwise compactification of $\ov{\pi}$.   The properness issue then ``disappears" but at the cost of losing the obvious smooth structure. 
Nevertheless,   Chen and Li recently showed in the general case of a symplectic manifold $(X,\om)$, that  one can define a smooth orbifold atlas on  $\ov{\M}_{0,l}(X,A,J)$, where the  charts   are given by gluing maps \cite{CL}.  
%Precisely they proved the following:
%\begin{theorem*}{\rm(}Chen, Li,  2006{\rm)} Let $(P,\om_P)$ be a symplectic manifold, and let $J_P$ be an $\om_P$-tame structure with respect to which transversality on every stratum of $\ov{\M}_{0,l}(P,\sigma,J_P)$ is realized. Then $\ov{\M}_{0,l}(P,\sigma,J_P)$ is a smooth orbifold.
%\end{theorem*}
There are many variants in the gluing of pseudo-holomorphic   spheres  procedure (see \cite{CL}, \cite{L}, \cite{MS}, \cite{Rvirt}, \cite{Sik},  among others), which appears naturally in Gromov-Witten theory, as well as in Floer theory.
%appears naturally in Gromov-Witten theory, as well as in Floer theory (see \cite{CL}, \cite{L}, \cite{MS}, \cite{Rvirt}, \cite{Sik},  among others), as an inverse procedure to Gromov's convergence.
 The  approach followed in \cite{CL} is to use \emph{balanced curves} in order to define  a natural slice for the action of  the  group of reparametrizations of  $S^2$, $\G$,  reducing the action of this latter non-compact group to that of $S^1$. As a consequence, they obtain gluing maps which are well-defined after quotient by the reparametrizations. 
Adapting their ideas to the  Hamiltonian fibration case, we construct  gluing maps $Gl^P$ and $Gl^B$ satisfying
$$\ov{\pi}\circ Gl^P=Gl^B\circ \ov{\pi}.$$
 Under some transversality assumptions  for a given fibered almost complex structure $J_P$ on $P$, realized  by projective fibrations over  projective space (as  described in \eqref{strongreg}),  we obtain the following:

\begin{fibrationmodule} The moduli spaces $\ov{\M}_{0,l}(P,\sigma, J_P)$ and $\ov{\M}_{0,l}(B,\sigma_B,J_B)$ are smooth orbifolds, and the map $\ov{\pi}$
%$$\F_{\pi}:\ov{\M}(P,\sigma,J_P)\lra \ov{\M}(B,\sigma_B,J_B),$$
restricts to a smooth locally trivial fibration {\rm(}of orbifolds{\rm)} above each stratum of  $\ov{\M}_{0,l}(B,\sigma_B,J_B)$. Moreover, the product formula can be recovered using integration over the fibers of  $\ov{\pi}$ above the top stratum of $\ov{\M}_{0,l}(B,\sigma_B,J_B)$.
\end{fibrationmodule}
 
 \noindent\textbf{Applications.} 
In 1997, Seidel  defined in \cite{Se}, a representation 
of the space of Hamiltonian loops of a given symplectic manifold in the automorphism group of the corresponding quantum homology.
 Lalonde, McDuff and Polterovich  have shown, under the relative semi-positivity assumption,  that the rational  cohomology of the total space splits as modules for any Hamiltonian fibration over $S^2$ \cite{LMP}.  McDuff  removed the semi-positivity assumption using virtual techniques  \cite{Mfib}. 
   In general, we say that a  fibration is \emph{rationally c-split} if $$H^{*}(P;\Q)\cong H^{*}(B;\Q)\otimes H^{*}(F;\Q)$$ as modules.
  This splitting is realized when the inclusion  $\iota_F^P:F\into P$ induces an injection in rational homology, and if, in addition, the second page of  the  Leray-Serre spectral sequence splits.  More generally,  Lalonde and McDuff conjectured that every Hamiltonian fibration verifies the $c$-splitting property \cite{LM}. They showed  that this splitting property holds for a large panel of Hamiltonian fibrations, in particular for Hamiltonian fibrations over $\CP^n$. The proof they give can be seen as a simple consequence of the product formula (when \eqref{ssp} is satisfied), and of  the invertibility of Seidel's morphism \cite{LMP}, \cite{Se}. 
  
  Another  consequence  is the symplectic uniruledness of Hamiltonian fibrations over  rationally connected bases. As defined in \cite{HLR}, a symplectic manifold $(X,\om)$ is   \emph{(symplectically) uniruled} if there is a non vanishing $GW$-invariant  with at least one point as a constraint. In other words if there exists $A\in H_2(X,\Z)$ and homology classes  $c_2,...,c_l\in H_*(X)$ such that:
\begin{equation*}\la pt,c_2,c_3,...,c_l\ra_{0,l,A}^X\neq 0.\end{equation*}
 A symplectic manifold is  \emph{rationally connected} if there is a non-zero $GW$-invariant involving  2 point  insertions \cite{HLR}, i.e the equation above is still true with $c_2=pt$.  In summary:

\begin{csplitting}\label{csp} 
Let $\pi: P\sra B$ be a Hamiltonian fibration. Assume  $B$ is rationally connected with respect to a class,  $\sigma_B\in H_2(B;\Z)$,  verifying the  hypothesis of Theorem B. Then,  $P$  is $c$-split and symplectically uniruled.
\end{csplitting}

\proof Let $C$ be the image of a  map  counted in $\la pt,pt,c_3^B,...,c^B_l\ra_{0,l,\sigma_B}^B\neq 0$. As already mentioned, $P_C$ is a Hamiltonian fibration over $S^2$, and by a result of 
 Lalonde, McDuff and Polterovich, \cite{LMP},   for every  $a\in H_*(F)$ there is an equivalence class $\sigma'$ of section classes in $P_C$, as well as an element $b\in H_*(F)$,  such that:
$$0\neq \la \iota_F^{P_C}(a),\iota_F^{P_C}(b)\ra_{0,2,\sigma'}^{P_C}=\la \iota_F^{P_C}(a),\iota_F^{P_C}(b),\iota_F^{P_C}([F]),...,\iota_F^{P_C}([F])\ra_{0,l,\sigma'}^{P_C},$$
where the last equality is a consequence of the Divisor axiom.
Applying the product formula, as given in \eqref{simpversionPF}, we conclude that:
\begin{equation}\label{preuveunireglage}
\la \iota_F^P(a),\iota_F^P(b),\pi^{-1}(c_3^B),...,\pi^{-1}(c^B_l)\ra_{0,l,\iota^P_{P_C}(\sigma)}^P\neq 0.
\end{equation}
Hence, by taking $a=pt$ we obtain that  $P$ is uniruled. Now, suppose  $\pi$ is not c-split. Then there exists  $a\in H_*(F;\Z)$ in the kernel of $\iota_F^P$. Therefore, the $GW$-invariants having  $\iota_F^P(a)$ as an entry, must vanish. But this would contradict  \eqref{preuveunireglage}. 
\qed\\

The proof of this corollary indicates, that unless we have a good knowledge of Seidel's morphism, the number of point insertions should a priori decrease.  Still, as a result, every Hamiltonian fibration over $(\CP^n,\om_{FS})$ is $c$-split and uniruled, where $\om_{FS}$ is the standard K\"ahler form on the complex projective space. The same applies for Hamiltonian fibrations over $(S^2\times S^2,\om\oplus\om)$  since in that case there is only one curve representing the diagonal and passing through two points.\\

The paper is organized as follows. Section 1, introduces the basic ingredients needed. We also define a particular affine connection on $P$, whose torsion is given by the symplectic curvature associated to the coupling form.  
In Section 2, we describe the linearization of the Cauchy-Riemann problem associated to the fibered almost complex structures. It is shown that the linearization is compatible with the projection $\pi$. 
Then we prove the \emph{structure theorems} in Section 4, ensuring the  $GW$-invariants are well-defined. In Section 5, we give the  proof of the product formula.
The Section $6$,  is devoted to showing the  locally trivial smooth (orbi-)fibration structure of $\ov{\pi}$. We then recover the product formula using integration over the fibers of $\ov{\pi}$. We conclude by  giving an example of a non-trivial induced fibration of moduli spaces.\\

\subsection*{Acknowledgements.} I would like to thank  Fran\c{c}ois Lalonde and Shengda Hu, for being so generous of their time and support,  and for the fruitfull discussions I had with them. I should also thank Dusa McDuff for bringing symplectic  uniruledness to my mind. Also, many thanks to Fr\'ed\'eric Bourgeois and his corrections and suggestions regarding this work. I should also thank R\'emi Leclercq and Liam Watson for their useful comments.

\section{The framework}

In this section we set the basic notions that will be needed in the rest of the paper. 
We recall that symplectic and Hamiltonian fibrations are classified by $\mathrm{BSymp}(F,\omega_F)$ and $\mathrm{BHam}(F,\omega_F)$, respectively. Again, for $b\in B$, let  $\omega_{b}$ be the induced symplectic form on the fiber $F_b:=\pi^{-1}(b)$. For the sake of clarity, we begin by recalling the notions of  Hamiltonian connections and coupling form. This exposition follows  \cite{GLS} and \cite{MS}, where the proofs of all the claims can be found.
 
\subsection{Hamiltonian connections and coupling form.}  Consider the vertical subbundle,  $Vert\subset TP$, over $P$, whose fiber at each point, $p\in P$,  is given by the subspace
  $Vert_p:=\ker d\pi(p)$. A connection on the fibration, $\pi:P\sra B$, is defined  by a  splitting of $T_pP$ for each $p\in P$: 
 $$T_pP=Hor(p)\oplus Vert_p,$$ 
where $Hor(p)$  is called the horizontal plane at $p$.
%Given such a connection one can lift paths from $B$ to $P$ and define parallel transport.
The notations,  $X^h$ and $X^v$, will refer to the horizontal and  vertical parts of a vector field $X$ on $P$, with respect to the above splitting.  Also, given a vector field $X$ on $B$, we will denote by $\ov{X}$ its horizontal lift to $TP$.  
Now, let $R$ denote the  \emph{the symplectic curvature}   associated to the connection. This is the 2-form on $B$ with values in $Vert$, such that  for $v$ and $w$  two vector fields on $B$: 
$$ R(v,w)(p):=[\ov{v},\ov{w}]^{v}(p),\hsp p\in P.$$
Any closed extension $2$-form $\tau$ of $\omega_F$ defines a connection                                                                                                                                                                                                                                                                                                               with horizontal planes:
$$Hor_{\tau}(p):=\{v\in
  T_pP\hspace{0.1cm}|\hspace{0.1cm}\tau(v,w)=0\hspace{0.2cm}\forall
  w\in Vert_p\}. $$
\noindent In fact, any 
form $\tau'\in\Omega^2(P)$ such that $\ker(\tau'-\tau)$ is in $Vert$,  defines the same horizontal plane field. Nevertheless, a specific choice can be made by requiring that:
\begin{equation*}\label{couplingcondition}
  \pi_*\tau^{n+1}=
  \int_F\tau^{n+1}= 0.
\end{equation*}

\noindent Such  $\tau$ is  called \emph{coupling form} associated to the connection, and its values on pairs of horizontal vectors is determined by the curvature of the connection:
\begin{equation*}-d(\tau(p)(v,w)):= \iota(R(d\pi(p)v,d\pi(p)w)(p))\omega_{\pi(p)}(p).
\end{equation*}
\noindent This results from the fact that the holonomy of the connection is Hamiltonian. 

%\begin{rem} We will often write $R(p)(v,w)$ instead of $\tau(p)(v,w)$. 
%\end{rem}

For transversality purposes we will need to allow  the connection to vary. For  $b\in B$, let $C_0^{\infty}(F_b)$ denote  the space of smooth functions on $F_b$ having zero mean value. We will  consider  Hamiltonian deformations of the coupling form $\tau$. By this we mean exact deformations,
 $$\tau_H:=\tau-d\tilde{H},$$
 where  $\wt{H}$ is a section in $C^{\infty}(\pi^*T^*B)$, i.e $\wt{H}_p$ is a cotangent vector in $T_{\pi(p)}B$,  and it satisfies the property that  for fixed  $b\in B$ and $X\in T_bB$ the function $H_{b,X}(p):=\wt{H}_p(X)$ belongs to  $C_0^{\infty}(F_b)$. The subset of $C^{\infty}(\pi^*T^*B)$ having this property will be denoted by  $\H$. By definition,
$\tau_H$ is a closed extension 2-form of $\omega_F$, and one verifies that its associated horizontal distribution is given by: 
$$Hor_{\tau_H}(p) =\{v-X_{\wt{H}_p(v)}(p)\hspace{0.1cm}|\hspace{0.1cm} v\in Hor_{\tau}(p)\},\hspa p\in P,$$
where $X_{\wt{H}_p(v)}$ denotes the (unique) Hamiltonian vector field on $F_{\pi(p)}$ induced by the function $\wt{H}_p(v)$. 

\begin{rem} Since we are working with functions having zero mean value, $\tau_H$ is nothing else but the coupling form associated to $Hor_{\tau_H}$. Also, observe that the symplectic curvature changes under exact deformation.  
\end{rem}

\subsection{Almost complex structures.}
An almost complex structure $J$ on a symplectic manifold, $(F,\omega)$, is a smooth section of  the bundle of endomorphisms of the tangent space of $F$, such that:
$$\forall p\in F,\hsp (J(p))^2=-id_{T_pF}.$$ 
Let  $\mathcal{J}(F)$ denote the set of  almost complex structures on $F$, and let  $\mathcal{J}_F\equiv \mathcal{J}(F,\omega)$  be the subset of $\mathcal{J}(F)$ consisting of \emph{$\omega$-tame} almost complex  structures, i.e such that $\omega(\cdot,J \cdot)$ is positive definite.  A choice of such a $J$, gives $TF$ the structure of a complex fiber bundle. Let  $c_1^{TF}\in H^2(F,\Z)$ denote the corresponding first Chern class. Since  $\mathcal{J}_F$ is contractible and non-empty ~\cite{MS}, this class is well-defined, and does not depend on our choice of compatible almost complex structure.  Let $\pi:P\sra (B,\om_B)$ be a Hamiltonian fibration  with fiber $(F,\om)$, and coupling form $\tau$. 

\begin{defn}\label{fiberedstruct}  An  almost complex structure  $J_P\in\mathcal{J}(P)$ is \emph{compatible with $\pi$ and $\tau$}, or \emph{fibered}, if and only if there exists $J_B\in\mathcal{J}(B,\omega_B)$ such that:
\begin{itemize}
  \item $d\pi \circ J_P=J_B\circ d\pi$,
  \item $J_b:=\left.J_P\right|_{F_b}\in\mathcal{J}(F_b,\omega_{b})$ for all $b\in B$,
  \item $J_P$ preserves the horizontal distribution induced by $\tau$.
\end{itemize}
We denote by $\mathcal{J}(P,\tau,\pi)$ the set of such almost complex structures.
\end{defn}

Each fibered almost complex structure determines  a family $\{J_b\}_{b\in B}$ of $\omega_{b}$-tame almost complex structures.
Let $\J^{V}\equiv \mathcal{J}^{V}(P,\pi,\om)$ denote the set of almost complex structures of the vertical subbundle which are $\omega_{b}$-tame on each fiber $F_b$. Also, let $c^v\in H^2(P,\Z)$ be the corresponding first Chern class.   Conversely, for any $J\in\mathcal{J}^{V}$ and $J_B\in\mathcal{J}_B$,  there exists a unique fibered structure $J_P$ extending $J$ and defined as the horizontal lift of $J_B$ on the horizontal distribution. Furthermore,  when deforming the coupling form using $H\in \H$, the unique extension of $J$ that projects on $J_B$ and that preserves the distribution induced by $\tau_H$ is given  by:
\begin{equation*} J_P^H(p)(v)=J_P(p)(v)+J(p)X_{H_{\pi(p)}(\pi_{*_p}v)}-X_{H_{\pi(p)}(J_B(\pi(p))\pi_{*_p}v)},\hspa p\in P,
 \end{equation*}
where $v\in T_pP$. Let  $\J_P\equiv\J(P,\pi)$ denote the union of the $\J(P,\pi,\tau_H)$ over all deformations $\tau_H$ of a given coupling form $\tau$.  We have:

\begin{lem} 
The space $\J_P$ is parametrized by the product  $\J_B\times\mathcal{J}^{V}\times\H$.
\end{lem} 

The above isomorphism is given by the choice of  $\tau$. In fact, the dependance is on the factor $\H$ corresponding to affine space of Hamiltonian connections. The choice of $\tau$ simply fixes the origin.

%Let $p_{i}$, $i=1,2,3$,  be the projection over the $i^\text{th}$ factor of $\J_P$, and  set $p_{ij}$, $i,j\in\{1,2,3\}$, $i<j$, to be the projections on the product of the $i^\text{th}$ and $j^\text{th}$ factors of $\J_P$. These projections will be used in the structure theorems section.

\begin{rem}
Note that, for any family $\{J_b\}_{b\in B}\in\mathcal{J}^{V}$ and any given $J_B$, we can find a positive $\kappa\in\R$ such that $J_P$ is $\omega_P$-tame for $\omega_P=\tau+\kappa\pi^*\omega_B$. 
\end{rem}
%Indeed, let $v\in TP$, then by definition  we have:
%\begin{equation*}
%\omega_P(v,J_Pv)=\omega(v^v,Jv^v)+R(\pi_*v,J_B\pi_*v)+\kappa\omega_B(\pi_*v,J_B\pi_*v).
%\end{equation*}
%so that $\omega_P$ is already $J_P$-tame in the vertical directions. Thus  it only remains to show that for $v\in Hor$,  the sum of the two last terms in the right member of the equation 
%%$$R(\pi_*v,J_B\pi_*v)+\kappa\omega_B(\pi_*v,J_B\pi_*v)>0$$
%is strictly positive for some positive function $\kappa$. But since  $\omega_B(\pi_*v,J_B\pi_*v)>0$ whenever $\pi_*v\neq0$ we can just set for all $z\in B$
%$$\kappa(z)>\max\Big\{\sup_{\{v\in TB|\|v\|=1\}}\{-R(z)(
%v,J_Bv)\},0\Big\}$$
%where the norm on $TB$ is the one induced by $\omega_B$ and $J_B$.
%\item Given  $J_B$, $J_P$, and the family  $J\in\mathcal{J}^{Vert}$,  we will set 
%$c_1(TB,J_B)\in H^2(B;\Z)$,  $c_1(TP,J_P)\in H^2(P;\Z)$ and $c_1(Vert,J)\equiv c^v\in H^2(P;\Z)$, to be the corresponding first Chern classes.
%\end{itemize}

\subsection{A specific affine connection.} Fix a coupling form $\tau$. We define an affine connection on $TP$,   extending the  \emph{vertical  Levi-Civita (L-C) connection} introduced in \cite{MS} (Chapter 8), and lifting the L-C connection on $TB$.  This construction will be needed in order to relate the linearization of the Cauchy-Riemann associated to  $J_P=(J_B,J,H)\in \J_P$, to the  linearization  of the Cauchy-Riemann operator  associated to $J_B$ (see Section 2).  First, let  $g_{J_P}$ be the hermitian metric on  $P$ defined as
\begin{equation*}\label{productmetric}
 g_{J_P}:=g_J\oplus \pi^*g_{J_B},
\end{equation*}
relatively to the splitting  $TP=Vert\oplus Hor_{\tau_H}$, where  $$g_{J_B}:=\frac{1}{2}(\om_B(.,J_B.)-\om_B(J_B.,.)),$$ and  $g_J:=\{g_{J_b}\}_{b\in B}$ is the analogous family of Hermitian metrics on $Vert$. Let $\nabla^B$ denote the L-C connection on  $TB$ relatively to $g_{J_B}$.  Also, set $\nabla^F$ to be the L-C connection  on  $TF_b$, with  $b\in B$, relatively to   $g_{J_b}$. For any vector fields  $X$ and $Y$ on $P$, set
%$$\nabla_X Y := [X^h, Y^v]^v + \nabla^F_{X_v} Y^v + [X^v, Y^h]^h + (\nabla_{\underline X^h} \underline Y^h)^h.$$
$$\nabla_X Y := [X^h, Y^v]^v + \nabla^F_{X_v} Y^v + [X^v, Y^h]^h + \ov{(\nabla^B_{\pi_* X^h}\pi_*Y^h)}.$$
This operation is clearly bilinear in  $X$ and $Y$. In fact, the sum of the first two terms corresponds to the  \emph{vertical L-C connection}, $\nabla^v$, which is the unique connection on $Vert$ induced by the Hamiltonian connection and which restricts to the L-C connection on $F$. The remaining part is what is needed  to extend this vertical connection to an affine connection on $P$ lifting the L-C connection on $B$, which torsion $T$  is given by the symplectic curvature:
\begin{eqnarray*}T(X,Y)=-R(X,Y)=[X^h,Y^h]^v,\hspa X,Y\in\mathcal{X}(TP).\end{eqnarray*}

We show  below that $\nabla$ is indeed  a connection. Since $\nabla^v$ is   a connection, it suffices to show that for  $f\in C^{\infty}(P)$ and $\xi\in Hor_{\tau}$,  and  any $w\in TP$:
$$\nabla_wf\xi=(w(f))(\xi)+f\nabla_w\xi\hspace{0.2cm},\hspace{0.2cm}\nabla_{fw}\xi=f\nabla_w\xi.$$
 Suppose that $w$ is vertical. Then, by definition: $$\nabla_wf\xi=[w,f\xi]^h=f[w,\xi]^h+w(f)\xi.$$ Analogously, we have   $[fw,\xi]=f[w,\xi]-\xi(f)w$, implying that 
 $$[fw,\xi]^h=f[w,\xi]^h$$
 since $w$ is vertical. Now, suppose that $w$ is horizontal. Let $\alpha_t$, $t\in(-\epsilon,\epsilon)$ be the flow of  $w$ starting at $p$, and  let $P^B_t(\pi_{*_p}w_p)$ denote the parallel transport along the projected curve $\pi(\alpha_t)(p)$. Then,
  \begin{eqnarray*}
 (\nabla_wf\xi)_p
 %&=&\left.\frac{d}{dt}\right|_{t=0}(P^B_{t}(\pi_{*_p}w_p)\pi_{*_{\alpha_t(p)}}f\xi)_p^h\\
% &=&\left.\frac{d}{dt}\right|_{t=0}(P^B_{t}(\pi_{*_p}w_p)f(\alpha_t(p))\pi_{*_{\alpha_t(p)}}\xi)_p^h\\
&=&\left.\frac{d}{dt}\right|_{t=0}(f(\alpha_t(p))P^B_{t}(\pi_{*_p}w_p)\pi_{*_{\alpha_t(p)}}\xi)_p^h=df(p)w\xi+f(p)(\nabla_w\xi)_p,\\
%&=& (w(f))(p)\xi+f(p)\left.\frac{d}{dt}\right|_{t=0}(P^B_{t}(\pi_{*_p}w_p)\pi_{*_{\alpha_t(p)}}\xi)_p^h\\
%&=&  df(p)w\xi+f(p)(\nabla_w\xi)_p,
\end{eqnarray*}
where the first equality is given by linearity of the parallel transport in $B$. We further have that 
$$\nabla_{fw}\xi=\left.\frac{d}{dt}\right|_{t=0}(P^B_{t}(\pi_{*_p}f(p)w_p)\pi_{*_{\alpha_t(p)}}\xi)_p^h=\left.\frac{d}{dt}\right|_{t=0}(P^B_{t}(\pi_{*_p}w_p)\pi_{*_{\alpha_{\psi(t)}(p)}}\xi)_p^h$$ for some reparametrization $\psi(t)$ of the interval, so that finally we have $f\nabla_w\xi$.  

\begin{rem}
Note that $\nabla^v$ may not preserve the metric $g_J$ along horizontal directions, whereas  $\wt{\nabla}^v=\nabla^v-\frac{1}{2}J(\nabla^v J)$ does (\cite{MS}, Lemma 8.3.6).
Furthermore, given any vector fields $w,\xi_1$ and $\xi_2$  in $TP$, one can show that $(\nabla_wg)(\xi_1,\xi_2)$ coincides with $(\nabla_{w^h}g)(\xi^v_1,\xi^v_2)$.
 It follows that  the $J_P$-preserving connection $\wt{\nabla}:=\nabla-\frac{1}{2}J_P(\nabla J_P)$  preserves $g_{J_P}$.
\end{rem}

Let $\exp$ stand for both the exponential maps with respect to $\nabla$ and $\nabla^B$. The following straightforward  identities will be useful in the gluing section:
\begin{lem} For $p\in P$,  $X\in T_pP$ and $q\in P$ in the injective radius of $\exp_p$, we have  
\begin{equation}\label{identiteexpo}\pi(\exp_p X)=\exp_{\pi(p)}\pi_{*_p}X, \hspa\pi_{*_p}\exp_p^{-1}(q)=\exp^{-1}_{\pi(p)}(\pi(q)).
\end{equation}
\end{lem}
%\proof  First observe that $\pi(\exp_p tX)$ is a curve in $B$ which is tangent to $\pi_*(X)$ at time 0. If $X\in \ker d\pi(p)$ then  $c_t(p):=\exp_p tX$ lies in the fiber $F_{\pi(p)}$ at each time $t$, and the first equality is verified. If $X$ is not vertical, extend it to a vector field $X=X^h+X^v$ on the restriction of $P$ to $\pi(\exp_p tX)$ such that $X=\dot{c}=\dot{c}^v+\dot{c}^h$ on the curve $c_t(p)$, and such that $X^h$ is the horizontal lift of $\pi_*\dot{c}$.  Then by definition of the exponential, $\nabla_{\dot{c}}\dot{c}=0$ and furthermore,  $[\dot{c}^v,\dot{c}^h]^h$ vanishes.  Consequently: 
%$$0=\pi_*(\nabla_{\dot{c}}\dot{c})=\pi_*(\nabla_{\dot{c}^h}\dot{c}^h)=\nabla^B_{\pi_*\dot{c}}\pi_*\dot{c},$$
%so that the first equality is verified. We leave the second equality to the reader.
\qed

\subsection{Curve independance} 

From  $(H_1)$ in the characterization of Hamiltonian fibrations we  deduce the following lemma, which plays a crucial part in the proof of the product formula:

\begin{lem} \label{curveindependance}
Let $\pi:P\sra B$ be a Hamiltonian fibration as above. Assume we have  $u_1,u_2\in C^{\infty}(S^2,B)$ such that  $[u_1(S^2)]=[u_2(S^2)]$. Then the restricted bundles $\left.P\right|_{u_1}$ and $\left.P\right|_{u_2}$ are isomorphic as Hamiltonian bundles.
 \end{lem}

\begin{pr} If $B$ is simply connected then it follows directly from Hurwicz isomorphism between $\pi_2(B)$ and $H_2(B,\Z)$.
Assume $B$ is not simply-connected. As $\mathrm{Ham}(F,\om)$ is connected, any classifying map for $P$, say %$f:B\sra B\Ham(F,\om)$
$f$, factorizes up to homotopy through a map  $f':B/B_1\sra \mathrm{BHam}(F,\om)$.   In other words, if $\pi_{B_1}$ denotes the projection from $B$ to $B/B_1$, the maps  $f$ and $f'\circ \pi_{B_1}$ are homotopic. Let $P':=(f')^*E\Ham(F,\om)$ and consider 
$$u'_1:=\pi_{B_1}\circ u_1,\hspa u'_2:=\pi_{B_1}\circ u_2.$$
These two maps represent the same homology class $\pi_{B_{1_*}}(\sigma_B)\in H_2(B/B_1;\Z)$, hence, $u'_1$ and $u'_2$ are homotopic so that:
$$\left.P\right|_{u_1}\cong (u'_1)^* P'\cong (u'_2)^* P'\cong \left.P\right|_{u_2}.$$
\end{pr}\qed

\section{The Cauchy-Riemann Fredholm problem in Hamiltonian fibrations}

Let $j_0$ denote the complex structure on $S^2=\C\cup\{\infty\}$ inherited from the multiplication by $i:=\sqrt{-1}$. In a general symplectic manifold $(X,\om)$ with $\om$-tame almost complex structure $J$,  a rational $J$-holomorphic sphere is a smooth map, $u:S^2\sra X$, satisfying the Cauchy-Riemann equation, 
$$\delbar_{J}u:=\frac{1}{2}(du+J\circ du\circ j_0)=0.$$
The set of all such solutions representing a given class $A\in H_2(X,\Z)$,
\begin{equation*}
 \wt{\M}(X,A,J):=\left\{u\in C^{\infty}(S^2,X)|\delbar_{J}u=0\hs,\hs [u(S^2)]=A \right\},
 \end{equation*} 
  is  the \emph{moduli space of parametrized $J$-holomorphic maps representing $A$}. 
 Consider the Fr\'echet space,
  $$\mathcal{B}_X(A):=\left\{u\in C^{\infty}(S^2, X)|[u(S^2)]=A\right\},$$
 and the space
$$\mathcal{E}_X(A,J)=\bigsqcup_{u\in\mathcal{B}_X(A)}\mathcal{E}_{X,u}(J):=\bigsqcup_{u\in\mathcal{B}_X(A)}C^{\infty}(\Lambda_{J}^{0,1}(S^2,u^*TX)).$$ 
The obvious projection from  $\mathcal{E}_X(A,J)$ to $\B_X(A)$, defines a locally trivial bundle  between Fr\'echet spaces. Then,  $\delbar_{J}$ is  a  section of this bundle, and $ \wt{\M}(X,A,J)$ is the corresponding zero set.  The linearized operator $D^X_u$ of  $\delbar_{J}$ at  $u\in \wt{\M}(X,A,J)$ is defined as the differential of $\delbar_{J}$ at $u$ composed with the projection on the fiber $\mathcal{E}_{X,u}(J)$. To give a meaning to this  vertical projection outside of the zero section, we consider the Hermitian connection on 
$\mathcal{E}_X(A,J)$,
 $$\wt{\nabla}^{X}:=\nabla^{X}-1/2(J(\nabla^{X} J)),$$
where $\nabla^{X}$ is the L-C connection on $X$ with respect to  the metric $g_{J}$.
Then,
 \begin{eqnarray*}
D^X_u: %\mathcal{X}_{B,u}
\mathcal{X}_{X,u}:= C^{\infty}(S^2,u^{*}TX)&\sra& \mathcal{E}_{X,u}(J)\\
\xi&\mapsto& \wt{\nabla}^{X}_{\xi}\ov{\partial}_{J}(u)=(\nabla_{du}^{X})^{0,1}\xi-  \frac{1}{2}J(u)(\nabla^{X}_{\xi}J)(\partial_{J}u), 
\end{eqnarray*}
where $(\nabla^{X})^{0,1}$ is the $J$ anti-linear part of $\nabla^{X}$.   It is well-known that  $D^{X}_u$ is Fredholm  for $u\in \wt{\M}(X,A,J)$. Moreover, if $\delbar_{J}$ is transversal to the zero section, then the moduli space is a smooth finite dimensional oriented manifold \cite{MS}, \cite{RT}.

Now, let $\pi:P\sra B$ be a Hamiltonian fibration with coupling form $\tau$, let $J_P$ be  a fibered complex structure relatively to $\tau$, and let $\sigma\in H_2(P,\Z)$.  The connection $\nabla$  we constructed
induces a splitting of the tangent space of $\E_P(\sigma,J_P)$ at all points, and projection to the fiber direction can again be defined. Hence, the linearization $D^P_u$ of  $\delbar_{J_P}$ can be defined for all $u \in \B_P(\sigma)$. Note that $\nabla$ is not Levi-Civita since  its torsion is  given by the symplectic curvature. This is what gives rise to the extra curvature  term in the following expression for $D^P_u$:

\begin{lem} 
Let $u\in \mathcal{B}_P(\sigma)$ and $\xi\in C^{\infty}(S^2,u^{*}TP)$. Then,
$$  D^P_u\xi=(\nabla_{du})^{0,1}\xi-\frac{1}{2}J_P(u)(\nabla_{\xi}J_P)(\partial_{J_P}u)+R^{0,1}(du^h,\xi^h),$$
where $\nabla^{0,1}$ and $R^{0,1}$, respectively stand for the $J_P$ anti-linear parts of $\nabla$ and $R$.
\end{lem}

%It follows from the general expression of   $D_u$
% that:
%\begin{equation*}
%D_u\xi=\frac{1}{2}(\nabla_{\xi}du+J_P(u)\nabla_{\xi}du\circ j_0)-\frac{1}{4}\left(J_P(u)(\nabla_{\xi}J_P(u)) du-(\nabla_{\xi}J_P(u)) du\circ j_0\right),
%\end{equation*}
%and the lemma follows by definition of $\nabla$ and the fact that its torsion is given by the symplectic curvature.

% \noindent\emph{Remark:} %\begin{itemize}\item 
% A simple computation shows that
%  \begin{eqnarray*}
% \delbar_{J_P}u=0&\Longleftrightarrow&\ov{\partial}_{J_P}u+X^{0,1}_{H(du)}=0,
%\end{eqnarray*}
%so that  a Hamiltonian deformation of the connection induces a deformation of the Cauchy-Riemann operator. \\
% \item Everything that follows below will still remain  valid for any other  choice of  complex structure on $S^2$, since all these structures are equivalent (under conformal equivalence) up to the action of $\G$.\\
%\end{itemize}

\begin{rem} When $\xi$ is vertically valued the curvature term disappears and we recover the vertical linearized operator introduced by McDuff and Salamon (\cite{MS}, Chapter 8). In the rest of the paper we will designate by $D^{v}$ the restriction of $D^P$ to vertically valued vector fields.
\end{rem}

\subsection{Splitting of Fredholm systems} Since we consider fibered almost complex structures on $P$, the projection $\pi$ naturally induces a map,
\begin{equation*}\label{fibstr}
\pi: \wt{\M}(P,\sigma,J_P)\sra \wt{\M}(B,\sigma_B,J_B),\hspa u\mapsto\pi(u),
\end{equation*}
 where $\sigma_B:=\pi_*\sigma$. As we will see,  $\pi$ induces a submersion between the Fredholm system $(\B_P(\sigma),\E_P(\sigma,J_P),\delbar_{J_P})$ and  $(\B_B(\sigma_B),\E_B(\sigma_B,J_B),\delbar_{J_B})$.  We recall the notion of Fredholm system as given in  \cite{CL}.
\begin{defn} A Fredholm system of index $d$ is a triple $(\mathcal{B},\mathcal{E},s)$ such that:
\begin{itemize}
\item[$(F_1)$] $\mathcal{E}$ is the total space of a Banach fibration over a Banach manifold $\B$, with  projection $p$,
\item[$(F_2)$] $s:\mathcal{B}\longrightarrow \mathcal{E}$ is a section such that, for all $x\in s^{-1}(0)$, the linearization $L_x:T_x\mathcal{B}\sra\mathcal{E}_x:=\pi^{-1}(x)$  of $s$ at $x$, is Fredholm of index $d$.
\end{itemize}
When, in addition, $s$ is proper, we say that the system is \emph{compact}.
\end{defn}

We call  the set, $s^{-1}(0)$, the \emph{moduli space} of the system.

\begin{ex} Let  $(X,\om)$ be a general symplectic manifold with $A\in H_2(X,\Z)$ and $\om$-tame almost complex structure $J$. Fix an integer $p>2$.  For  $u\in \mathcal{B}_X(A)$  we respectively equip  $\mathcal{X}_{X,u}$ and $\mathcal{E}_{X,u}(J)$ with  $W^{1,p}$  (Sobolev) and  $L^p$ norms (relatively to the metric  $g_{J}$ and a fixed volume form on $S^2$). 
 Explicitly, given  $\xi\in \mathcal{X}_{X,u}$ and $\eta\in \mathcal{E}_{X,u}(J)$:
\begin{equation*}
\|\xi\|_{1,p}=\left(\int_{S^2}(|\xi|_{g_{J}}^p+|\nabla\xi|_{g_{J}}^p)dvol_{S^2}\right)^{\frac{1}{p}},\hspace{0.3cm}
\|\eta\|_p=\left(\int_{S^2}|\xi|_{g_{J}}^pdvol_{S^2}\right)^{\frac{1}{p}}.
\end{equation*}
We denote by  $\mathcal{X}_{X,u}^{1,p}$ and $ \mathcal{E}_{X,u}^p(J)$  the completed vector spaces with respect to these norms, and by  
$\mathcal{B}^{1,p}_X(A)$ and $\mathcal{E}^p_X(A,J)$ the corresponding completions of $\mathcal{B}_X(A)$ and  $\mathcal{E}_X(A,J)$. 
Under these completions,  $D^X_u$ is Fredholm so that the triple,
$$(\mathcal{B}^{1,p}_X(A),\mathcal{E}^p_X(A,J),\delbar_{J}),$$ satisfies conditions ($F_1$) and ($F_2$). This Fredholm system  will in general not be compact. 
\end{ex}

Next we define the following natural notion of morphisms between two Fredholm systems $(\mathcal{B},\mathcal{E},s)$ and $(\mathcal{B'},\mathcal{E'},s')$.

\begin{defn} A map $\Pi:=(\pi,\ov{\pi}):(\mathcal{B},\mathcal{E},s)\sra (\mathcal{B'},\mathcal{E'},s')$ between Fredholm systems is a   Banach orbifold vector bundle morphism, i.e
$$\pi:\mathcal{B}\to\mathcal{B'},\hspace{0.5cm} \ov{\pi}:\mathcal{E}\to\mathcal{E'},$$
such that 
\begin{itemize}
\item[i)] $s'\circ \pi=\ov{\pi}\circ s$
\item[ii)] $L'_{\pi(x)}\circ d\pi(x)=d\ov{\pi}(x,0)\circ L_{x},$ for every $x\in s^{-1}(0)$.
\end{itemize}
We say that $\Pi$ is  a submersion if, furthermore, $d\pi$ and $d\ov{\pi}$ are surjective.
\end{defn}

When $\Pi$ is a submersion, we directly extract the following exact sequence:
\begin{equation}\label{splittingsequence}
 \ker L''_x  \rightarrowtail  \ker L_x  \rightarrow  \ker L'_{\pi(x)}
  \rightarrow \coker L''_x  \rightarrow \coker L_x \twoheadrightarrow \coker L'_{\pi(x)} . 
\end{equation}
\noindent where $L''_x$ stands for the operator $L_x$ restricted to $\ker d\pi(x)$. 
From exactness of this sequence we deduce splitting at the level of operator indices, i.e for every 
$x\in s^{-1}(0)$:
$$Ind(L_x)=Ind(L_{\pi(x)}')+Ind(L''_x). $$

Obviously, if both $L''_x$ and $L'_{\pi(x)}$ are unobstructed (i.e their cokernels vanish), then $L_x$ is also unobstructed.  As a result,  $\ker L_x$ is isomorphic to $\ker L''_x\oplus\ker L'_{\pi(x)}$ (upto a choice of a section from $\ker L'_{\pi(x)}$ to $\ker L_x$ ).

\begin{defn} A submersion $\Pi$ between Fredholm systems is a \emph{splitting}, if the sequence \ref{splittingsequence} is obstruction free, i.e if for all  $x\in s^{-1}(0)$:
$$ \coker L''_x = \coker L_x= \coker L'_{\pi(x)}=0. $$
\end{defn}

\subsection{Splitting of Fredholm systems in the Hamiltonian fibration case.}
Here, we show that a Hamiltonian fibration $\pi:P\sra B$ together with a coupling form $\tau$, induce a submersion  between the Fredholm systems  $(\mathcal{B}_P(\sigma),\mathcal{E}_P(\sigma,J_P),\delbar_{J_P})$ and $(\mathcal{B}_B(\sigma_B),\mathcal{E}_B(\sigma_B,J_B),\delbar_{J_B})$. In fact, $\pi$ induces the maps:
\begin{equation*}
\pi:\mathcal{B}_P(\sigma) \sra\mathcal{B}_B(\sigma_B), \hspace{0.5cm}f\mapsto \pi\circ f,
\end{equation*}
and
\begin{equation*}
\ov{\pi}:\mathcal{E}_P(\sigma,J_P) \lra\mathcal{E}_B(\sigma_B,J_B), \hspace{0.5cm}\eta\mapsto \pi_*(\eta).
\end{equation*}
That $\pi$  induces a submersion between Fredholm systems can then be summarized as follows:
\begin{lem}\label{proj} The following diagrams are commutative:
\begin{equation*}
\xymatrix{ \mathcal{E}_P(\sigma,J_P) \ar[r]^{\ov{\pi}} & \mathcal{E}_B(\sigma_B,J_B)\\
\mathcal{B}_P(\sigma)\ar[u]^{\delbar_{J_P}}\ar[r]^{\pi}& \mathcal{B}_B(\sigma_B) \ar[u]^{\delbar_{J_B}}}
\hspace{1.0cm}
\xymatrix@C=4.0pc{\mathcal{E}_{P,u}(J_P) \ar[r]^{d\ov{\pi}(u,0)} & \mathcal{E}_{B,\pi(u)}(J_B)\\  
\mathcal{X}_{P,u}\ar[r]^{d\pi(u)}\ar[u]^{D_u} & \mathcal{X}_{B,\pi(u)} \ar[u]^{D_u^B}}
\end{equation*}
\end{lem}

\proof For the left handside diagram, commutativity follows from the fact that  $J_P$ is fibered. Thus, $\pi$ maps  $(\delbar_{J_P})^{-1}(0)$ into $(\delbar_{J_B})^{-1}(0)$.  We prove that the second diagram commutes. Namely, we show  that for every $\xi\in\mathcal{X}(u^*TP)$:
$$ \pi_{*} D_u \xi= D^B_{\pi\circ u}\pi_{*}\xi.$$
This is a consequence of the following three identities:
let $\xi$ and $X$ be vector fields on $P$, and let $p\in P$, then
\begin{itemize}
\item[(i)]$\pi_{*}[\xi^v,J_PX^h]=J_B\pi_{*}[\xi^v,X^h]$,
\item[(ii)]$\pi_{*}(\nabla_{\xi}J_P)_p X=(\nabla^B_{\pi_{*}\xi}J_B)_{\pi(p)}(\pi_{*}X)$,
\item[(iii)]$\pi_{*}(\nabla^{0,1}_{du}\xi)=(\nabla^B_{\pi_{*}du})^{0,1}\pi_{*}\xi-\pi{_*}[(\ov{\partial}_{J_P}(u))^h,\xi^v]^h$.
\end{itemize}
Assuming these are verified  we obtain that:
$$ \pi_{*} D_u \xi =(\nabla^B_{d(\pi\circ u)})^{0,1}\pi_{*}\xi-\pi{_*}[(\ov{\partial}_{J_P}(u))^h,\xi^v]^h-\frac{1}{2}J_B(\pi\circ u)(\nabla^B_{\pi_{*}\xi}J_B)(\pi_{*}\partial_{J_P}u),$$
since  $\nabla^{0,1}_{du}\xi^v$ and $R^{0,1}(du,\xi)$ are vertically valued.  Moreover,  since $\pi_{*}\partial_{J_P}u=\partial_{J_B}u_B$ it only remains to show that:
$$[(\ov{\partial}_{J_P}(u))^h,\xi^v]_{u(z)}^h=0 \hspa \forall z\in S^2.$$
But, if  $(\delbar_{J_P}u)(z)^h\neq 0$ for some $z$,   consider the horizontal lift $\ov{X}$ of $(\delbar_{J_B}(\pi(u)))(z)$, which is defined on $T (\pi^{-1}(u(z)))$, and agrees with $(\delbar_{J_P}u)(z)^h$ at $u(z)$. Then,
$$[(\ov{\partial}_{J_P}(u))^h,\xi^v]_{u(z)}=[\ov{X},\xi^v]_{u(z)},$$
where the right handside vanishes since $\ov{X}$ is constant along vertical directions.

Now,  equality  (i) follows by definition of the bracket, by holomorphicity of the projection, and since the flow of a vertical vector field starting at a point $p$ remains in the fiber above $\pi(p)$. Equality (ii) is just a consequence of (i), since this latter  is equivalent to
 $$[\xi^v,J_PX^h]^h=J_P[\xi^v,X^h]^h,$$
 which by definition of $\nabla$ is the same as:
  $$\nabla_{\xi^v}(J_PX^h)=J_P\nabla_{\xi^v}(X^h).$$
Hence, $ (\nabla_{\xi^v}J_P)(X^h)=0$,
which, combined with the fact that the connection is vertically valued when its two entries are in $Vert$, gives us the second equation.

For   (iii)  we have that:
$$\nabla^{0,1}_{du}\xi=\frac{1}{2}(\nabla_{\xi}du+J_P(u) \nabla_{\xi}(du\circ j_0))-R^{0,1}(du,\xi).$$
Hence, by definition of the connection and because the curvature is vertically valued, we get 
\begin{eqnarray*}
2\pi_*(\nabla^{0,1}_{du}\xi) & = &\pi_*\nabla_{\xi}du+J_B(\pi(u)) \pi_*\nabla_{\xi}(du\circ j_0)\\
 & = &  \nabla^B_{\pi_*\xi}d(\pi\circ u)+\pi_*[\xi^v,(du)^h]^h+J_B(\pi(u)) \nabla^B_{\pi_*\xi}d(\pi\circ u)\circ j_0 \\ & &+ J_B(\pi(u))\pi_*[\xi^v,(du\circ j_0)^h]^h \\
 &=& \nabla^B_{d(\pi\circ u)}\pi_*\xi+J_B(\pi(u))\nabla^{B}_{d(\pi\circ u)\circ j_0}\pi_*\xi+\pi_*[\xi^v,(du)^h+J_P(u)(du\circ j_0)^h]^h \\
 &=& 2(\nabla^B_{d(\pi\circ u)})^{0,1}\pi_*\xi+\pi_*[\xi^v,2(\ov{\partial}_{J_P}(u))^h]^h,
\end{eqnarray*}
where the third equality is due to (i) and the fact that  $\nabla^B$ is torsion free. The last follows since $J_P$ preserves the horizontal distribution and the vertical subbundle.
\qed\\

The symplectic connection on $P$ induces the splittings:
\begin{equation*}
\mathcal{E}_{P,u}=\Gamma(\Lambda_{J_P}^{0,1}(S^2,u^{*}TP^h))\oplus\Gamma(\Lambda_{J_P}^{0,1}(S^2,u^{*}TP^v))=:\mathcal{E}_{P,u}^h(J_P)\oplus\mathcal{E}_{P,u}^v(J_P),
\end{equation*}
and
\begin{equation*}
\mathcal{X}_{P,u}=\Gamma(S^2,u^*TP^h)\oplus\Gamma(S^2,u^*TP^v)=\mathcal{X}_{P,u}^h\oplus\mathcal{X}_{P,u}^v.
\end{equation*}
In this splitting, $D^P_u$ 
takes the following matrix form:
\begin{equation}\label{matricelinearisee}
\left(\begin{array}{cc}
	(D^B_{\pi(u)})^h &  0\\
	L_u & D^{v}_u
\end{array}\right)
\end{equation}
where $L_u$ is    linear and given by:
\begin{equation}\label{connectant}L_u: \mathcal{X}_{P,u}^h\lra\mathcal{E}_{P,u}^v,\hspa \xi\mapsto -\frac{1}{2}J(u)(\nabla_{\xi}J)(\partial_{J_P}u)^v+R^{0,1}(du^h,\xi).
\end{equation}
\noindent Thus, applying the diagram \eqref{splittingsequence}, we obtain the exact sequence:
\begin{equation}\label{splitting}
\ker D_u^{v} \rightarrowtail \ker D^P_u\rightarrow \ker D_{\pi(u)}^B\rightarrow \coker D_u^{v}\rightarrow \coker D^P_u\twoheadrightarrow \coker D_{\pi(u)}^B
\end{equation}
where the connectant is given by the restriction of $L_u$  to the kernel of $D^B_{\pi(u)}$. 

\begin{lem}\label{fibsplitting}
A Hamiltonian fibration induces, up to transversality, a splitting of Fredholm system for $(\mathcal{B}^{1,p}_P(\sigma),\mathcal{E}^p_P(\sigma,J_P),\delbar_{J_P})$. 
\end{lem}

The issues regarding transversality are considered in the next section.

%%%%%%%%%%%%%%%%%%%%%%%%%%%%%%
%%%%%%%%%%%%%%%%%%%%%%%%%%%%%%%
% TRANSVERSALITE
%%%%%%%%%%%%%%%%%%%%%%%%%%%%%%%%
%%%%%%%%%%%%%%%%%%%%%%%%%%%%%%%%

\section{Structure theorems}

Let $\wt{\M}_{0,l}(P,\sigma,J_P)$ be the moduli space of parametrized $J_P$-holomorphic maps with $l$ marked points, representing $\sigma$. This space consists of  uples, 
$$(u,x_1,...,x_l)\in \wt{\M}(P,\sigma,J_P)\times (S^2)^l,$$
where the points $x_1,..., x_l\in S^2$ are pairwise distinct. The group $G:=\G$ of reparametrizations of $S^2$  acts (diagonally) on this moduli space.
The quotient of $\wt{\M}_{0,l}(P,\sigma,J_P)$ under this action is usually not compact. Still, its ``Gromov's compactification", $\ov{\M}_{0,l}(P,\sigma,J_P)$,  is a stratified space  consisting of stable holomorphic maps \cite{Hummel},\cite{MS},\cite{RT}. Concretely, the  stratification is given by the \emph{combinatorial-type} of the \emph{labeled connected trees with tails} modeling the stable maps.   When these strata are automorphism free they can be given a manifold structure for a generic choice of fibered almost complex structure. After giving the description of the stable holomorphic maps in Hamiltonian fibrations, we show that the latter structure theorem holds compatibly with the projection $\pi$. This is the content of Theorem B',  which is a slight extension of Theorem B. We conclude  by  observing that  changing of generic almost complex structure induces a cobordism between the corresponding moduli spaces.

\subsection{Stable holomorphic maps and Hamiltonian fibrations} In order to fix notations and terminology, we begin by  introducing, first, the combinatorics needed to describe stable holomorphic maps, and  second,  the moduli space of stable pseudo-holomorphic maps.  

\subsubsection{Labeled graphs.} This exposition is mostly taken from \cite{Costello}. Nevertheless, emphasis is put on  how the combinatorics change when given a homomorphism of semi-groups. 

Let $\S=(V,\textrm{Fl};\text{pr},\varrho)$ be a graph. Here, $V$ denotes the set of vertices,  $\text{Fl}$ denotes the set of flags. The map $\text{pr}:\text{Fl}\sra V$ assigns to a given flag, its associated vertex, i.e for $v\in V$, the set $\text{pr}^{-1}(v)$ gives the valences at $v$. The notation $|v|$ denotes the number of valences at $v$.  The map, $\varrho: \text{Fl}\sra \text{Fl}$ is an involution. The set $E$ of edges in $\S$, is identified to the 2-elements orbits of $\varrho$, while the fixed points of $\varrho$ gives the set, $T$, of tails of $\S$.  In the rest of the paper, we use  the notation, $v E w$ ($v,w\in V$), to indicate that there is an edge between the $v$ and the $w$ vertices of $\S$.

Denote by $C$ the set of connected components of $\S$. There is an induced map, $\text{pr}_C:T\sra C$, assigning a tail to the connected component it is attached to. The genus, $g(\S)$, of the graph is the Euler number of $\S$, i.e 
$$g(\S)=|V|-|E|+1.$$ 
When all the components of $\S$ have genus 0, we say that it is a \emph{forest}; if, furthermore, $\S$ is connected, we say that it is a \emph{tree}. In the rest of  the paper we  will only consider forests. 

The subgraph, $(v,\text{Fl}_v:=\text{pr}^{-1}(v);\left.\text{pr}\right|_{\text{Fl}_v},id)$, is called the \emph{star} of $\S$ at $v$. Given a subset, $V'\subset V$, the \emph{complement} $\S\bs V'$ of $V'$ in $\S$, is defined as the subgraph given by $\S$ from which we removed all  the stars at points coming from $V'$.

A map, $\beta: V\sra B$, where $B$ is an Abelian semi-group, defines a \emph{B-labeling} of $\S$. We denote by $(\S,\beta)$ the corresponding labeled graph. A $B$-labeling $\beta$ induces a $B$-labeling of $C$:
$$\beta_C: C\sra B,\hspa \beta_C(c):=\sum_{v\in c}\beta(v).$$
An example of labeling, is the \emph{genus-labeling} given by a map, $h:V\sra \N$, assigning to each vertex the corresponding genus. In what follows,   the genus labeling will always be trivial. 

%The genus of $(\S,h)$ is then:
%$$g(\S,h):=g(\S)+ \sum_{v\in V}h(v).$$  
%From now on, we will only consider labeled graphs such that $g(\S,h)=0$.

\textbf{Composition of graphs.} The \emph{composition}, $\S'':=\S'\Rightarrow \S$, of two graphs $\S'$ and $\S$, is the graph obtained by replacing the vertices in $\S$ by the components of $\S'$. More precisely,  $\S''$ is determined by an isomorphism of maps between finite sets:
$$[\text{pr}'_C:T'\sra C']\cong[\text{pr}:\text{Fl}\sra V],$$ 
that is, we have isomorphisms $T'\stackrel{\cong}{\sra}\text{Fl}$, and $C'\stackrel{\cong}{\sra}V$, commuting with $\text{pr}_C'$ and $\text{pr}$. If such isomorphisms exist, we say that $\S'$ and $\S$ are \emph{composable}. In that case, $\S''$ is given as follows:
$$V''=V',\hsp \text{Fl}''= \text{Fl}'\hsp \text{pr}''=\text{pr}',\hsp\text{and}\hsp \varrho''=\begin{cases}\varrho' & \text{on elements not fixed by $\varrho$} \\
\varrho & \text{otherwise.}
\end{cases} $$
This notion of composition defines a partial order $\prec$ on graphs: we set $\S''\prec\S$. More generally, labeled graphs can be composed. Let $(\S,\beta)$ and $(\S',\beta')$ be two labeled graphs such that $\S$ and $\S'$ are composable, and let $\pi:B'\sra B$ be a homomorphism of semi-groups.  We say that the graphs are \emph{$\pi$-composable}, if $\beta$ and $\beta'_C$ commute with respect to $\pi$. The \emph{$\pi$-composition} is the labeled graph, $(\S'',\beta'')$, where  $\beta'':=\pi\circ \beta'  $. This defines a $B$-labeling. The \emph{labeled  composition} of $(\S,\beta)$ with $(\S',\beta')$, is the $B'$-labeled graph $(\S'',\beta')$.

 Note that composition of graphs induces a surjective map on the sets of vertices, $\gamma: V''\sra V$, and of flags, $\gamma_{\text{Fl}}: \text{Fl}''\sra \text{Fl}$.  A \emph{section} of the composition is  a pair of right inverses, $\iota: V\sra V''$ and $\iota_{\text{Fl}}: \text{Fl}\sra \text{Fl}''$, such that $(\text{pr}'')^{-1}(\iota(v))=\iota_{\text{Fl}}(\text{pr}^{-1}(v))$, which further verify $\beta''\circ \iota = \pi\circ \beta'\circ \iota= \beta$ when  we have labelings. In this latter case, we say that $(\iota,\iota_{\text{Fl}})$ is a \emph{labeled section}. Note that a section always exists while a labeled one may not. Let $(\S'',\beta')$ be a labeled composition of $(\S',\beta')$ with $(\S,\beta)$, and let $(\iota,\iota_{\text{Fl}})$ be a labeled section.  The components of  $\S''\bs \iota(V)$, together with the restriction of the labeling $\beta'$,  are called \emph{contracted components}. Note that $(\S,\beta)$ is obtained by replacing the contracted components with either tails (if the contracted component had some) or edges between two vertices of $V$. We  say  that $(\S,\beta)$ is the \emph{contraction} of $(\S'',\beta')$. As an example, consider the case  where  we contract a vertex $v\in V''$ with $|v|\leq2$. First, delete the star at $v$, thus creating either one or two tails. Assume the graph is connected. If it has only one vertex, everything disappears. Suppose it has more than one vertex. If $v$ had a tail, then you don't delete the new tail created on its neighbour.  If $v$ had two edges, we  replace the deleted vertex by an edge obtained by gluing together the two new tails. In the case it only had one edge and no tail, we also delete the new tail. This operation is called \emph{elementary contraction}.     

\textbf{Stability of graphs.} Let $\S$ be  a forest.  A vertex $v\in V$ is said to be \emph{stable} if $|v|\geq 3$. The forest  $\S$ is said to be  \emph{stable}, if all its vertices are stable. If, furthermore, we are given a $B$-labeling $\beta$ on $\S$, we say that a vertex $v\in V$ is ($B$)-\emph{stable} if either $\beta(v)\neq0$ or $v$ is stable. The labeled graph $(\S,\beta)$ is $B$-\emph{stable} when all its vertices are ($B$)-stable. 
%From now on, when considering stable graphs, we will drop the genus labeling from the notations. 
The following lemma,  showing how the combinatorics change under  a homomorphism $\pi: B'\sra B$ of   semi-groups, is standard.

\begin{lem} Suppose $(\S'',\beta')$ is a  stable labeled tree, then we have one of the following:
(1)  there exist stable tree $(\S,\beta)$ and a stable forest $(\S',\beta')$ such that $(\S'',\pi\circ \beta')$ is their $\pi$-composition; (2) $|T''|\leq2$. Moreover, in case (1),  there is a labeled section of the composition.
\end{lem}

\proof Suppose $(\S'',\pi\circ \beta')$ is not stable. This implies that there is at least one vertex of $\S''$ which is not stable. Let $v$ be such a vertex, then $|v|\leq 2$ and $\pi(\beta'(v))=0$. In the case $|v|=0$, we directly have (2). If $|v|>0$, then contract the unstable vertices inductively until we get a stable forest, or we get (2). These elementary contractions don't change the number of tails and preserve connectedness of the graph.  Suppose that the resulting tree, $(\S,\beta)$, is  stable. By construction, we have  natural inclusions $\iota: V\hra V''$ and $\iota_{\text{Fl}}: \text{Fl}\hra \text{Fl}''$. The forest, $(\S',\beta')$,  is then given by $\S'=(V'', \text{Fl}'';\text{pr}'',\varrho')$. We define $\varrho'$. By the definition of composition, each edge of $\S$, that is each 2-element orbit $\{f_1,f_2\}$ of $\varrho$, should correspond to two tails in $\S'$. So, letting $\{f'_1:=\iota_{\text{Fl}}(f_1),f'_2\}$ be the 2-element orbit of $f_1'$ under $\varrho''$, we then set $\varrho'(f_1')=f'_1$ and $\varrho'(f'_2)=f'_2$. Now, for the  flags $f\in \text{Fl}''$ not coming from $\S$, we simply set $\varrho'(f)=\varrho''(f)$. \qed\\

The proof indicates that the tree $(\S,\beta)$ is unique while $(\S',\beta')$ is not in general. The natural section mentioned in the proof will be called \emph{$\pi$-section}.

\begin{defn} Let $\pi:B'\sra B$ be a homomorphism of semi-groups. We say that the tree $(\S,\beta)$ above is the \emph{$\pi$-stabilization} of $(\S'',\beta')$ and we denote it $S_{\pi}(\S'',\beta')$. The vertices in $\S$ are called \emph{$\pi$-stable}, and the contracted components are called \emph{$\pi$-unstable}.
\end{defn}

Let $(\S,\beta)=S_{\pi}(\S'',\beta')$. The vertices in the image of section are said to be \emph{$B$-stable}. Now, let $c$ be one of the contracted components in $\S''$. Then $c$ is a tree with  at most 2 tails. A tail of $c$ is called \emph{exterior} if it is  also a tail of $\S''$. Otherwise, it is called \emph{interior}. In the case  $c$ has 1 interior tail, we say that $c$ is a \emph{($\pi$-)contracted branch}. In the case  $c$ has two interior tails, we say that $c$ is a \emph{($\pi$-)connecting branch}. For a connecting branch there is a unique path of  vertices connecting the 2 tails. This path forms a \emph{connecting chain}.  

\subsubsection{Moduli space of stable maps} The discussion below is taken from \cite{MS} and \cite{RT}, where the proofs of all the claims can be found.
Denote by $\ov{\M}_{g,l}$ the Deligne-Mumford moduli space of stable curves with genus $g$ and $l$-marked points. We shall only consider the case $g=0$. Points in $\ov{\M}_{0,l}$ are given by isomorphism classes of elements, $$\j\equiv(\Sigma,\j,\mathbf{x}:=(x_1,...,x_l)),$$ 
where $(\Sigma,\j)$ is a nodal Riemann surface of arithmetic genus $0$  with no self-intersection, together with $l$ pairwise disjoint marked points on $\Sigma$ denoted  $\mathbf{x}$, which  are  disjoint  from the nodes. Furthermore, each component of $\Sigma$ has at least three special points (i.e  marked points or  nodes).  Two pointed nodal curves $(\Sigma,\j,\mathbf{x})$ and $(\Sigma',\j',\mathbf{x}')$ are isomorphic, if there is  a diffeomorphism $\varphi:\Sigma \sra\Sigma'$ satisfying:
$$\varphi^*\j'=\j,\hsp \text{and}\hsp \varphi(x_i)=x'_i.$$
 We denote by  $\aut(\Sigma,\j,\mathbf{x})$
 % the isomorphism class of  $(\Sigma,\j,\mathbf{x})$, and by $\aut[\Sigma,\j,\mathbf{x}]$ \
the automorphisms of $(\Sigma,\j,\mathbf{x})$, i.e the subset of diffeomorphisms of nodal surfaces, $\varphi:\Sigma\sra \Sigma$, such that $\varphi^*\j=\j$ and $\varphi(\mathbf{x} )=\mathbf{x}$. This group is invariant under isomorphism of nodal curves.  Also,  it is standard that the elements of  $\ov{\M}_{0,l}$ are automorphism free. 
 
Let $\ov{\M}_{0,l}(X,A,J)$ denote the compactified moduli space of stable $J$-(pseudo)-holomorphic maps, from (nodal) curves of genus 0 with $l$ marked points   into the symplectic manifold $X$, representing the class $A\in H_2(X,\Z)$. Points in  $\ov{\M}_{0,l}(X,A,J)$ are given by isomorphism classes of parametrized stable pseudo-holomorphic maps $(\j,u)\equiv((\Sigma,\j,\mathbf{x}), u)$, where $(\Sigma,\j,\mathbf{x})$ is a Riemann nodal curve of genus 0 with $l$ marked points (not necessarily stable),  and, $u:\Sigma\sra X$, is $(\j,J)$- holomorphic and such that each component on which $u$ is constant has at least three special points. We say that $(\j, u)$ is isomorphic to $(\j',u')$, if there is an isomorphism of pointed nodal curve, $\varphi$, between $\j$ and $\j'$, such that $\varphi^*u=u'$. Let $\aut(\j,u)$ denote the corresponding automorphism group. It is well-known that stability implies finiteness of the automorphism groups. 
%Let $[\j,u]$ denote the isomorphism class of $(\j,u)$, and denote by $\aut[\j,u]$ the associated automorphism group. It is well-known that stability implies finiteness of the automorphism groups. 
Moreover, if the map $u$ is \emph{reduced}, or \emph{simple} in the sense that $(\j,u)$ has \emph{no ramified component or any two component having the same  non-constant image in $X$}, then $\aut(\j,u)=id$ (when the map has one only one component, this notion of simple map coincides with the notion of somewhere injective map). It is well known that any stable pseudo-holomorphic map can be reduced to a  simple stable map. The reduction process however changes   the  homology class of the map (see \cite{MS}). 
 
 The moduli spaces $\ov{\M}_{0,l}$ and $\ov{\M}_{0,l}(X,A,J)$ are stratified, with strata labeled by stable labeled trees of genus 0 called \emph{stable stratum data} of the strata, or \emph{combinatorial type}. 
  For a tree $\S=(V,\text{Fl};\text{pr},\varrho)$, the set $V$ corresponds to the components of $(\Sigma,\j)$, while  $\text{Fl}$ corresponds to  the set of special points on the curve. For  $\ov{\M}_{0,l}(X,A,J)$ we have a $H_2(X,\Z)$-labeling giving the homology class represented by the image of each component in $X$. Note that   $\ov{\M}_{0,l}$ and $\ov{\M}_{0,l}(X,A,J)$ coincide when $X=pt$ and $A=0$. 
A strata with stratum data $\S$ will be denoted $\M_{\S}$ or $\M_{\S}(X,J)$, and $\M^*_{\S}(X,J)$ will denote  the subset of $\M_{\S}(X,J)$ consisting of simple stable maps.  Furthermore, the partial
order $\prec$ on labeled graphs induces inclusions of strata: 
 $$\M_{\S'}\subset \ov{\M}_{\S} \Leftrightarrow \S'\prec \S$$
 and similarly for $\M_{\S}(X,J)$.  Note that there are finitely many strata in the compactified moduli space, since,  given $l$ and $A$, the  set $\D_{0,l}^A$ of possible combinatorial types for genus 0 stable maps, with $l$ markings representing $A$, is finite.
 
  As we are considering only genus 0 stable maps, we can in fact fix the complex structure on each component of the nodal surface to be the standard one. Let  $\wt{\M}_{\S}(X,J)$ denote the moduli space of parametrized stable $J$-holomorphic maps $(\Sigma,u,\mathbf{x})$ representing $\S$, and let  $G_{\S}$ be the \emph{reparametrization group} of the domain $(\Sigma, \mathbf{x})$. The stratum $\M_{\S}(X,J)$  is then identified to the quotient of $\wt{\M}_{\S}(X,J)$ under  the (proper) action of $G_{\S}$.  For the sake of the reader, we describe it. First, note that a stable map, $(\Sigma,u,\mathbf{x})\in \wt{\M}_{\S}(X,J)$, is determined by a triple, $(u, \mathbf{y},\mathbf{x})$, where $\mathbf{y}:=\{y_{vv'}\}_{vEv'}$ for $v,v'\in V$, is the data given by the nodal points in $\Sigma$. Let $\Sigma_v$ denote the component of $\Sigma$ corresponding to $v\in V$, and let $u_v$ be the restriction of $u$ to $\Sigma_v$. The group  $G_{\S}$   consists of  pairs $(\{\varphi_v\}_{v\in V},\gamma)$,  where $\gamma \in \aut(\S)$ is a tree-with-tails automorphism, and 
  $\varphi_v:\Sigma_v\sra\Sigma_{\gamma(v)}$ 
 is an element of $\G$.  Then,
 $$(\{\varphi_v\},\gamma)\cdot(u,\mathbf{y},\mathbf{x}):=\left(\{u_v\circ\varphi_v^{-1}\},\{\varphi_v(y_{vv'})\}_{vEv'},\{\varphi_{p(x_k)}(x_k)\}_{k\in\{1,...,l\}}\right),$$
is the considered action of $G_{\S}$ on $\wt{\M}_{\S}(X,J)$.
Before carrying out the description of  stable $J_P$ holomorphic maps in a Hamiltonian fibration $\pi:P\sra B$,  we first make sure that we have the appropriate energy bounds in order to apply Gromov's compactness (\cite{MS}, \cite{RT}). 

\subsubsection{Energy identities} Suppose $J_P$ is a fibered structure obtained from a connection $\tau$, an element $J_B\in\J_B$ and a family $J\in\J^{V}$, and let $g_{J_P}$ be the corresponding split metric on $P$. For a smooth map $u:S^2\sra P$, we define its \emph{total energy} to be its Dirichlet norm with respect to $g_{J_P}$:
$$E(u):=\frac{1}{2}\int_{S^2}\|du\|^2_{g_{J_P}}dvol_{S^2}.$$
Since $g_{J_P}$ is split, $E(u)$ can be written as  the sum, 
$$\frac{1}{2}\int_{S^2}\|d(\pi(u))\|^2_{g_{J_B}}dvol_{S^2}+\frac{1}{2}\int_{S^2}\|(du)^v\|^2_{g_{J}}dvol_{S^2}:=E_B(\pi(u))+E^{vert}(u),$$
where $E_B(\pi(u))$ is the energy of  $u_B:=\pi(u)$ (with respect to $J_B$), and $E^{vert}(u)$ is the \emph{vertical energy}. When $u$ is $J_P$-holomorphic,  it turns out that:
$$E_B(u_B)=\int_{S^2}u_B^*\omega_B,\hsp\text{and}\hsp E^{vert}(u)=\int_{S^2}u^*\tau +\int_{S^2} R(u)dvol_{S^2},$$
where the second identity is obtained since: 
$$\tau(du,J_P du)=\omega(du^v,J du^v)- R(du^h,J_P du^h).$$
Consider the \emph{Hofer norm} of the symplectic curvature,
 $$\|R\|_H:=\int_B \left(\max_{p\in F_b} R(p)-\min_{p\in F_b}R(p)\right)\omega_B^{n_B},$$
which is  bounded by compactness of $P$. Then, we obtain the following upper bound:
\begin{lem}
  For every $J_P$-holomorphic map $u$:
\begin{equation}\label{nrj}
 E(u)\leq\int_{S^2}u^*\tau +\|R\|_H+\int_{S^2}u_B^*\omega_B.
\end{equation}
\end{lem}
Applying Gromov's compactness, we conclude that any sequence  of simple $J_P$-holomorphic map representing $\sigma$,  must converge (up to taking a sub-sequence) to a stable $J_P$-holomorphic map.  From the upper bound \ref{nrj}, and  the stability condition for stable maps we deduce the following standard result \cite{MS}, \cite{RT}:
\begin{lem}\label{strata} We have  $|\mathcal{D}^{\sigma}_{0,l}|<\infty$.
\end{lem}
 
\subsubsection{Forgetful maps and Hamiltonian fibrations.}
 The natural map $\pi_{pt}:X\sra pt$, induces the map: $\pi_{{pt}_*}: H_2(X,\Z)\sra \{0\}$ on the labeling groups. We then obtain the standard \emph{forgetful} map:  
  $$\ov{\pi}_{{pt}}:\ov{\M}_{0,l}(X,A)\sra \ov{\M}_{0,l},\hspa \ov{\pi}_{{pt}}(\j,u)=\j^{st},$$ 
  where  $\j^{st}$ denotes the stabilization of $\j$, that is $\j^{st}$ is obtained by contracting the unstable components recursively. It is well-known that, when restricted to strata,  the forgetful map defines maps:
  $$\ov{\pi}^{\S}_{{pt}}:\M_{\S}(X)\sra \M_{S_{\pi_{{pt}_*}}(\S)}.$$
  
  This procedure can be generalized  to any Hamiltonian fibration  $\pi: P\sra B$ with coupling form $\tau$. Again, the moduli spaces of stable pseudo-holomorphic maps in $P$ and $B$ are stratified, with strata labeled by stable labeled trees. Note that $\pi$ induces a map $\pi_*:H_{2}(P,\Z)\sra H_2(B,\Z)$ between the labeling groups. Let $J_P$ be a $(\pi,\tau)$  compatible almost complex structure on $P$ projecting on $J_B$. Given a stable stratum data $\S_P$ for $J_P$-holomorphic maps in $P$, representing a class $\sigma\in H_2(P,\Z)$, we have  that $S_B:=S_{\pi_*}(\S_P)$ is  a stable stratum data for  $\ov{\M}_{0,l}(B,\pi_*(\sigma),J_B)$, and we have a \emph{$\pi$-forgetful map}:
  $$\ov{\pi}_{\S_P}:\M_{\S_P}(P,J_P)\sra \M_{\S_B}(B,J_B),\hspa (\j,u)\to \left(\j^{st,\pi},u_B:=\pi(u)\right),$$ 
  where $\j^{st,\pi}$ is the Riemann nodal surface consisting of the components determined by the $\pi_*$-section of the $\pi_*$-stabilization of labeled graphs (which has combinatorial type $S_{\pi_*}(\S_P)$), and $u_B$ restricts to the  $\pi_*$-stable components. We further have a reparametrization-group-equivariant map: 
  $$\pi_{\S_P}:\wt{\M}_{\S_P}(P,J_P)\sra \wt{\M}_{\S_B}(P,J_B),\hsp (u,\mathbf{y},\mathbf{x})\mapsto (u_B, \mathbf{y}_B:=\mathbf{y}^{st,\pi},\mathbf{x})$$
lifting  $\ov{\pi}_{\S_P}$.  Note that, a priori,   $\ov{\pi}_{\S_P}$ may not respect simplicity of pseudo-holomorphic stable maps (e.g. by sending a simple $J_P$-holomorphic map to a multiply covered $J_B$-holomorphic map).   This has dramatic effects regarding transversality within the range of fibered almost complex structures. For each stable stratum data $\S_P$ we will consequently restrict our attention to the subset, $\MMC^{**}_{\S_P}(P,J_P)$, of simple stable elements in $\MMC_{\S_P}(P,J_P)$
 lying in the preimage of $\M^*_{\S_B}(B,J_B)$ under $\ov{\pi}_{\S_P}$. In particular,  we  denote by  $(\pi_{\S_P}^{-1}(\ov{u}_B))^*$ the set of simple parametrized stable pseudo-holomorphic maps  lying in the fiber of $\pi_{\S_P}$ above $\ov{u}_B\in\wt{\M}^*_{\S_B}(B,J_B)$, and we   use the notation $(\ov{\pi}_{\S_P}^{-1}(\ov{u}_B))^*$ to denote the corresponding quotient under reparametrizations.

%%%%%%%%%%%%%%%%%%%%%%%%%%%%%%%%%%%%%%
%%%%%%%%%%%%%%%%%%%%%%%%%%%%%%%%%%%%%%%
% TRANSVERSALITY ON EVERY STRATA
%%%%%%%%%%%%%%%%%%%%%%%%%%%%%%%%%%%%
%%%%%%%%%%%%%%%%%%%%%%%%%%%%%%%%%%%%%%%

\subsection{Transversality on every strata} We begin by recalling some standard notations and fact concerning  transversality for a symplectic manifold $(X,\om)$. Then, we  apply it to Hamiltonian fibrations. Finally, we formulate the corresponding cobordism invariance.

\subsubsection{The non-fibered case.}\label{transingeneral} Consider a stable stratum data $\S_X=(V,\text{Fl};\text{pr},\varrho)$ for $\ov{\M}_{0,l}(X,A,J)$, with homological labeling $\beta$. For every $v\in V$, let $\sigma_v:=\beta(v)$ and set 
$$\wt{\MMC}^{*}(X,\beta,\J_X):=\left\{(u:=\{u_v\}_{v\in V},J)| J\in \J_X \hsp\text{and}\hsp u_v\in \wt{\MMC}^{*}(X,\sigma_v,J) \right\}.$$
 This  defines a subset of $$\mathcal{B}^{1,p}_X(\beta,\J_X):=\prod_{v\in V}\mathcal{B}_X^{1,p}(\sigma_v)\times \J_X.$$
 We describe $\wt{\M}_{\S_X}^{*}(X,J)$ as a subset of $\wt{\MMC}^{*}(X,\beta,\J_X)\times I(\S_X)$, where the \emph{incidental subvariety} $I(\S_X)$ is the subset of all uples,
 $$(\mathbf{y},\mathbf{x}):=(\{y_{vv'}\}_{vEv'},x_1,...,x_l)\in (S^2)^{2|E|}\times (S^2)^{|V|},$$    
 such that  for every  $v\in V$, the points $y_{vv'}$ for $vEv'$, and $x_m$ with $\text{pr}(x_m)=v$, are disjoint. 
 First, to each edge of the graph $\S_X$, we associate   a copy of the diagonal $\Delta_X\subset X^2$, and  set \emph{the edge diagonal}, $\Delta_E\subset X^{2|E|}$, to be the product of these diagonals over the set $E$ of all edges. We have a natural map, $$ev_E: \wt{\MMC}^{*}(X,\beta,\J_X)\times I(\S_X)\rightarrow (X)^{2|E|},$$
 called \emph{universal edge evaluation map}, assigning to each pair $(u,J, \mathbf{y},\mathbf{x})$ the corresponding ``evaluation at the nodes": 
 $$u(\mathbf{y}):=\{(u_v(y_{vv'}),u_{v'}(y_{v'v}))\}_{vEv'}.$$ 
  The preimage of $\Delta_E$ under $ev_E$, is the \emph{(parametrized) universal moduli space} denoted $\wt{\MMC}_{\S_X}^{*}(X,\J_X)$. Then, it is easy to see that 
  $$\wt{\M}_{\S_X}^{*}(X,J)=(p^{\J_X})^{-1}(J)\cap \wt{\MMC}_{\mathcal{S}_X}^{*}(X,\J_X),$$
where  $p^{\J_X}$ denotes the projection from $\wt{\MMC}^{*}(X,\beta,\J_X)\times I(\S_X)$ to $\J_X$. The standard transversality theorem asserts the following:
\begin{theorem}(\cite{MS},\cite{RT})\label{standardstructheorem} There exists a subset $\J_{X,reg}(\S_X)\subset \J_X$ of second category, such that for each  $J\in \J_{X,reg}(\S_X)$, the moduli space $\M_{\S_X}^{*}(X,J)$  is a smooth oriented manifold of dimension:
 $$\dim (\M^*_{\S_X}(X,J))=2n+2\sum_{v\in V} c_1^{TX}(\sigma_v)+2l-2|E|-6.$$
\end{theorem}

\noindent The set $\J_{X,reg}(\S_X)$ of \emph{regular $\om$-tame almost complex structure for $\S_X$} is explicitly given by the following conditions: 
\begin{enumerate}[i)]
 \item for every $v\in V$, for every $u\in  \wt{\M}^{*}(X,\sigma_v,J)$,  the linearization  $D^X_u$, of  $\delbar_J$  at $u$,    is surjective.
%$$\mathcal{V}_{reg}(J,\sigma):=\big\{(J_B,H)\in \J_B\times\H\hsp|\hsp \coker D_u^H=0\hsp\forall u\in\wt{\M}^{**}(P,\sigma,J_P)\big\}.$$$(J_B,H)\in \mathcal{V}_{reg}(J,\sigma)$ for every element $\sigma\in \vec{\sigma}$ such that $\pi_*\sigma=\sigma_B\neq 0$,
%\item $J\in \mathcal{J}^{vert}_{reg}(\omega,\sigma)$ for every $\sigma\in H_2(P;\mathbb{Z})$ such that $\pi_*\sigma=0$,
\item the restriction of  $ev_E$ to $\wt{\M}_{\S_X}^{*}(X,J)$ is transversal to $\Delta_E$. 
%$(p^{\J_P})^{-1}(J_P)$ 
\end{enumerate}
%i) for every $i\in V$, for every $u_i\in \wt{\M}^*(X,\sigma_i,J)$, the linearized operator $D_{u_i}$ is surjective; ii) the restriction of  $ev_E$ to $\wt{\M}_{\S}^{*}(X,J)$ is transversal to $\Delta_E$.  
Concretely, by i),  for regular $J$, the moduli spaces  $\wt{\M}^*(X,\sigma_v,J)$, ($v\in V$), are naturally oriented manifolds of real dimension 
$$\ind(D^X)=2n+2c_1^{TX}(\sigma_v).$$
Point ii) then implies that $\wt{\MMC}_{\mathcal{S}_X}^{*}(X,J)$ is a smooth oriented manifold. Since the $6(|E|+1)$-dimensional group,  $G_{\S_X}$, acts freely and properly by orientation preserving diffeomorphisms on this latter manifold, it follows that  $\MMC^{*}_{\mathcal{S}_X}(X,J)$ is a smooth oriented manifold of the stated dimension.

We briefly sketch the idea of proof for the genericity  of  $\J_{X,reg}(\S_X)$. We refer to  \cite{MS} for the details. The main idea,  is to show that the universal moduli space is a separable Banach manifold, and that $p^{\J_X}$ is a Fredholm map between separable Banach manifolds in order to apply  Sard-Smale. Of course, this does not apply straighforwardly here, since, for instance, $\J_X$ is not Banach. Instead, we consider, $\J_X^r$, the set of $\om$-tame almost complex structure of class $C^r$, with $r\geq2$, and let $\mathcal{E}_X^{p,r}(\beta)$  be the disjoint union over $\J_X^r$ of the direct sums: $$\bigoplus_{v\in V} \E^p_{X,u_v}(J)\hsp\text{where}\hsp (u,J)\in\mathcal{B}_X^{1,p}(\beta,\J_X^r).$$
By a standard argument, $\mathcal{E}_X^{p,r}(\beta)$ is a  locally trivial $C^{r-1}$ Banach fibration  over $\mathcal{B}_X^{1,p}(\beta,\J_X^r)$. This Banach  fibration admits the section
$$s_X:(u,J)\mapsto \delbar_{J}(u)$$
and clearly, $\wt{\MMC}^{*}(X,\beta,\J_X^r)$ is a subset of $s_X^{-1}(0)$. Moreover, the linearization, $\wt{D}^X_{(u,J)}$, of $s_X$ at $(u,J)$ is given by:
 \begin{eqnarray}\label{linearizationformula} \wt{D}^X_{(u,J)}:\bigoplus_{v\in V} \X^{1,p}_{X,u_{v}}\oplus T_{J}\mathcal{J} &\lra& \bigoplus_{v\in V}\E^{p}_{X,u_{v}}(J)\\
(\mathbf{\xi},Y) &\mapsto&D^X_{u}\mathbf{\xi}+\frac{1}{2}Y.(du). j_0 \nonumber
\end{eqnarray}
One can show by a standard argument that this operator is onto for every $(u,J)$. It follows that  $\wt{\MMC}^{*}(X,\beta,\J_X^r)$  is a separable Banach manifold. One concludes that $\wt{\M}^*(X,\J^r_X)$ is a Banach manifold by showing that the edge evaluation map $ev_E$ is transversal to $\Delta_E$. This is done recursively on the set of all labeled forest, the induction argument being made on the number of edges of the forests. 

Finally, a simple computation  shows that  $p^{\J^r_X}$ is Fredholm, with the same  kernel and cokernel   as the linearized operator $\wt{D}^X$.  Thus,  by  Sard-Smale (when $r$ is big enough) one obtains a generic subset $\J_{X,reg}^r(\S_X)\subset \J_X^r$,   which in fact coincides with the regularizing set defined above,  but with $C^r$ elements only.  One concludes in the $C^{\infty}$ case by  an argument due to Taubes (see \cite{MS}, Chapter 3).

\subsubsection{The Hamiltonian fibration case} Let $\pi:P\sra B$ be a Hamiltonian fibration with coupling form $\tau$. Before reformulating Theorem B and giving its proof, we fix some notations.

Fix a stable stratum data $\mathcal{S}_P=(V,\text{Fl};\text{pr},\varrho)$ with homological labeling $\beta_P$, and let  $\mathcal{S}_B=S_{\pi_*}(\S_P)=(V_B,\text{Fl}_B;\text{pr}_B,\varrho_B)$, with corresponding labeling $\beta_B$. Let $E$ and $E_B$ denote the corresponding set of edges. We will  assume below that the homomorphism $\beta_B':=\pi_*\circ\beta_P$ is non-zero, thus forcing $\beta_B$ to be non zero. Set $\sigma_v:=\beta_P(v)$. 

Let, $\wt{\MMC}^{**}(P,\beta_P,\J_P)$, be the restriction of $\wt{\MMC}^{*}(P,\beta_P,\J_P)$ to   simple maps having simple projection under $\pi$, and let  $$\wt{\MMC}_{\mathcal{S}_P}^{**}(P,\J_P):=ev_E^{-1}(\Delta_E)$$  be the corresponding universal moduli space. Similarly to the general case, we say that $J_P\in \J_{P,reg}(\S_P)$ if and only if: $i)$  $\hsp\forall v\in V,\hsp \forall u\in  \wt{\M}^{**}(P,\sigma_v,J_P)$,  the operator $D^P_u$  is onto; $ii)$  the restriction of  $ev_E$ to $\wt{\M}_{\S_P}^{**}(P,J)$ is transversal to $\Delta_E$. 
% \begin{enumerate}[i)]
% \item for every $i\in V$, for every $u\in  \wt{\M}^{**}(P,\sigma_i,J_P)$  the operator $D^P_u$  is surjective.
%\item the restriction of  $ev_E$ to $\wt{\M}_{\S_P}^{**}(P,J)$ is transversal to $\Delta_E$. 
%\end{enumerate}

Note that  we have a natural map
\begin{eqnarray*}\Pi_{\S_P}: \wt{\MMC}^{**}(P,\beta_P,\J_P)\times I(\S_P)&\sra&\wt{\MMC}^{*}(B,\beta_B,\J_B)\times I(\S_B)\\ 
 (u,J_P,\mathbf{y},\mathbf{x})&\mapsto& (u_B,J_B,\mathbf{y}_B,\mathbf{x}),
\end{eqnarray*}
induced from $\pi$ and the projection $p_1:\J_P\sra\J_B$. In fact, the  restriction  of $\Pi_{\S_P}$ to $\wt{\MMC}_{\mathcal{S}_P}^{**}(P,\J_P)$ induces a map between universal moduli spaces: 
$$\wt{\pi}_{\S_P}:=\pi_{\S_P}\times p_{1}:\wt{\MMC}_{\mathcal{S}_P}^{**}(P,\J_P)\sra\wt{\MMC}_{\mathcal{S}_B}^{*}(B,\J_B).$$
Since the edge evaluation maps commute with $\wt{\pi}_{\S_P}$ and the projection from $P^{2|E|}$ to $B^{2|E_B|}$, it follows easily that if  $J_P$ is regular for $\S_P$ then $J_B$ is regular for $\S_B$. 

Now, fix a regular $J_B$. Note that 
\begin{equation*}\label{realisationfibre}
(\ov{\pi}_{\S_P}^{-1}(\ov{u}_B))^*=(p_{23}\circ p^{\J_P})^{-1}(J,H)\cap \Pi_{\S_P}^{-1}(u_B,J_B) \cap ev_E^{-1}(\Delta_E\cap F^{2|E|}),
\end{equation*}
where $p_{23}:\J_P\sra \J^{V}\times \H$ is the obvious projection. The set $\J\H_{reg}(\ov{u}_B,J_B,\S_P)$ of \emph{fiber regularizing for $(\ov{\pi}_{\S_P}^{-1}(\ov{u}_B))^*$}, actually consists of pairs  $(J,H)\in \J^{V}\times \H$ that  turn this intersection  into an oriented manifold. More precisely, a pair $(J,H)$ is  \emph{fiber regularizing for $(\ov{\pi}_{\S_P}^{-1}(\ov{u}_B))^*$}  if 
for every    $(u\equiv\{u_v\}_{v\in V},\mathbf{y},\mathbf{x})\in(\pi_{\S_P}^{-1}(\ov{u}_B))^*$,  the following conditions are satisfied:
\begin{enumerate}[a)]
\item $\forall v\in V, \hsp \forall u_v\in\wt{\M}^{**}(P,\sigma_v,J_P)$,  the operator $D^{v}_{(u_v,J_P)}$ is onto;
\item the restriction of  $ev_E$ to $(\pi_{\S_P}^{-1}(\ov{u}_B))^*$ is transversal to $\Delta_E\cap F^{2|E|}$. 
\end{enumerate}
Concretely, it follows from $a)$ and $b)$, and the exact sequence \ref{splitting}, that for a regularizing pair $(J,H)$, the set $\pi_{\S_P}^{-1}(\ov{u}_{B})^*$ is a smooth oriented manifold with dimension:
\begin{equation*}\dim (\pi_{\S_P}^{-1}(\ov{u}_{B})^*)= \dim\wt{\MMC}^{**}_{\mathcal{S}_P}(P,J_P)-\dim\wt{\MMC}^{*}_{\mathcal{S_B}}(B,J_B).
\end{equation*}
Since we are only considering vertical deformations of the maps, the dimension of the reparametrizations is $6(|E|-|E_B|)$, thus giving after quotient:
$$\dim(\ov{\pi}_{\S_P}^{-1}(\ov{u}_{B}))^*=2n_F+2\sum_{v\in V}c^v_1(\sigma_v)-2|E|+2|E_B|.$$
We give the proof of the following extension of Theorem B to any strata, which can be seen as a mild extension to the Hamiltonian fibration case of Theorem  \ref{standardstructheorem}.
\begin{structurethm}\label{transversality2}
Let $\mathcal{S}_P$  and $\mathcal{S}_B$ be as above.
\begin{itemize}
 \item [1)]  $\J_{P,reg}(\mathcal{S}_P)$ is of second category in $\J_P$. 
\item[2)]  For any $J_B\in\mathcal{J}_{B,reg}(\mathcal{S}_B)$  and any $\ov{u}_{B}\in \MMC^{*}_{\mathcal{S}_B}(B,J_B)$, the set of regularizing pairs, 
 $\J\H_{reg}(\ov{u}_{B},J_B,\S_P)$, is of second category in $\J^{V}\times\H$.
\end{itemize} 
\end{structurethm}

\proof  The notations used here are the same as the one introduced in Section \ref{transingeneral}.
Following the guideline given in the preceding section, we show that  
\begin{equation}\label{Banachman}\wt{\pi}_{\S_P}:\wt{\M}_{\S_P}^{**}(P,\J_P^r)\sra\wt{\M}_{\S_B}^{*}(B,\J_B^r),\hsp r\geq 2, 
\end{equation}
 is a submersion between separated $C^{r-1}$ Banach submanifolds of $\mathcal{B}^{1,p}_P(\beta_P,\J_P^r)$ and  $\mathcal{B}^{1,p}_B(\beta_B,\mathcal{J}_B^r)$. By  the discussion in Section \ref{transingeneral} we already know that $\wt{\M}_{\S_B}^{*}(B,\J_B^r)$ is a Banach manifold.  We will then proceed as follows:
 \begin{enumerate}[(I)]
 \item we show that  $\wt{\mathcal{M}}^{**}(P,\beta_P,\J_P^r)$ is a Banach manifold and that the restriction of the natural map,
 \begin{equation*}
p:=\pi_{\S_P}\times p_1:\mathcal{B}_P^{1,p}(\beta_P,\J_P^r) \sra\mathcal{B}_B^{1,p}(\beta_B,\J_B^r),\hspa
(u,J_P)\mapsto(\pi(u),J_B).
\end{equation*}
to this product  moduli space is a smooth submersion onto $\mathcal{M}^{*}(B,\beta_B,\J_B^r)$;
\item  assuming  $ev_E$ is transversal to $\Delta_E$, we show that $\wt{\pi}_{\S_P}$ is a submersion;
\item we show that for every  labeled forest $\S_P$ with $l$-tails, the restriction of  $dev_E$ to $\ker d\Pi_{\S_P}$ is transversal to the subspace $\left.T\Delta_E\right|_{TF^E}$. This implies that $ev_E$ is transversal to $\Delta_E$ for any labeled forest $\S_P$ since  $ev_{E_B}$ is transversal to  $\Delta_{E_B}$ for every labeled forest. Hence, $\wt{\M}_{\S_P}^{**}(P,\J_P^r)$ is a Banach manifold.  
\end{enumerate}
The rest of the proof is verbatim the same as in the non-fibered situation and we omit it. 

\underline{\emph{Proof of (I)}}. First, by a standard argument (\cite{Hthesis}) involving the use of both the connection induced by $\tau_H$ and the L-C connection on $TB$, the sets, $\mathcal{E}_P^{p,r}(\beta_P)$ and $ \mathcal{E}^{p,r}_B(\beta_B)$,  are locally trivial $C^{r-1}$ Banach fibrations, over $\mathcal{B}_P^{1,p}(\beta_P,\J_P^r)$ and $\mathcal{B}_B^{1,p}(\beta_B,\J_B^r)$, which can be locally trivialized compatibly $p$.
 Let $s_P$ and $s_B$ be the corresponding Cauchy-Riemann sections. The  linearization $\wt{D}^B_{(u_B,J_B)}$ of $s_B$ at  $(u_B,J_B)$  is given by \eqref{linearizationformula}, while 
%By definition, the sets  $\wt{\MMC}^{*}(B,\beta_B,\J^r_B)$ and $\wt{\MMC}^{**}(P,\beta_P,\J_P^r)$ are subsets of $s_B^{-1}(0)$ and $s_P^{-1}(0)$, respectively. 
for the linearization  $\wt{D}^P_{(u,J_P)}$ of $s_P$ at $(u, J_P)$:
\begin{eqnarray*}\wt{D}^P_{(u,J_P)}:\bigoplus_{v\in V} \X^{1,p}_{P,u_v}  \oplus  T_{J_P}\J_P &\sra &\bigoplus_{v\in V} \E^{p}_{P,u_v}(J_P)\\
(\xi,Y^v,Y,f) \,\,\,\,\,\,\,\,\, &\mapsto & D_{u}\xi+\frac{1}{2}(Y .d\pi(u). j_0)^{h}  +\frac{1}{2}Y^v.(du)^v. j_0+X_{f(du)}^{0,1}
\end{eqnarray*}
%$$s_B:(u_B,J_B)\mapsto \delbar_{J_B}(u_B)\hspace{0.5cm}\text{and}\hspace{0.5cm} s_P:(u, J_P)\mapsto\delbar_{J_P}(u)$$
 Let $\ov{p}$ be the fibration map corresponding to the projection $p$:
 $$\ov{p}:\mathcal{E}^{p,r}(\beta_P)\lra \mathcal{E}_B^{p,r}(\beta_B),\hspa (\eta,J_P)\mapsto (d\pi(\eta),J_B).$$ By definition, $\ov{p}\circ s_P=s_B\circ p$.  
 %$$\wt{D}^B_{(u_B,J_B)}:\bigoplus_{i\in I^h} \X^{1,p}_{B,u_{B,i}}\oplus T_{J_B}\mathcal{J}_B \lra \bigoplus_{i\in I^h}\E^{p}_{B,u_{B,i}}$$ 
%be the linearization of $s_B$ at $(u_B,J_B)$ so that
 Furthermore, from lemma  \ref{proj}, and since $X_{f(du)}^{0,1}$ and $Y^v.(du)^v. j_0$ are vertically valued, we deduce that:
  $$d\ov{p}\circ \wt{D}^P_{(u,J_P)}= \wt{D}^B_{(\pi(u),J_B)}\circ dp.$$
The maps $dp$ and $d\ov{p}$ being both surjective,  the pair $(p,\ov{p})$ defines a submersion of Fredholm systems, 
%between $(\mathcal{E}^{p,r}(\beta_P),\mathcal{B}^{1,p}_P(\beta_P,\J_P^r),s_P)$ and $(\mathcal{E}_B^{p,r}(\beta_B),\mathcal{B}^{1,p}_B(\beta_B,\J_B^r),s_B)$ 
and we end up having the exact sequence:
\begin{equation*} 0\rightarrow \ker \wt{D}^v \rightarrow \ker \wt{D}^P  \rightarrow \ker \wt{D}^B \rightarrow \coker \wt{D}^v  \rightarrow \coker \wt{D}^P\rightarrow\coker \wt{D}^B\sra 0,
 \end{equation*}
where $ \wt{D}^v$ denotes the vertical operator associated to $\wt{D}^P$. 

To prove the claim, we show  that  the pair $(p,\ov{p})$ defines a splitting when we restrict $s_P$ to  $\wt{\MMC}^{**}(P,\beta_P, \J_P^r)$ and  $s_B$ to $\wt{\MMC}^*(B,\beta_B,\mathcal{J}_B^r)$.  It is enough to prove this when the tree structure of $\S_P$ is preserved under $S_{\pi_*}$.  Let $V_0$ denote the subset of $V$  on which $\beta'_B$ vanishes, and denote by $V_+$ its complement in $V$. Since $\wt{D}^B$ is onto for every $(u_B,J_B)\in\wt{\MMC}^*(B,\beta_B, \mathcal{J}^r_B)$, 
 it suffices  to show that the vertical operator $\wt{D}^v$ is surjective  at every points of $\wt{\MMC}^{**}(P,\beta_P, \J^r_P)$. Notice that $\wt{D}^v$ is closed, and suppose it is not dense. Then, by  Hahn-Banach, we would have a non-zero element: 
$$\{\eta_v\}_{v\in V}\in\bigoplus_{v\in V} L^{q}(\Lambda^{0,1}(S^2,u_{v}^*TP^v)),\hspa \frac{1}{p}+\frac{1}{q}=1$$
 such that each $\eta_v$, is of class $W^{1,p}$, is in the cokernel of $D^v_{u_{v}}$, and is such that:
\begin{eqnarray*}0&=&\int_{S^2}\sum_{v\in V_+}\LA\frac{1}{2}Y^v.(du_v)^v. j_0+X_{f(du_v)}^{0,1},\eta_v \RA +\sum_{v\in V_0}\LA \frac{1}{2}Y^v.(du_v)^v. j_0, \eta_v \RA dvol_{S^2}.
\end{eqnarray*}
%Actually, as $\wt{D}^B_{(\mathbf{\pi(u)},J_B)}$ is onto and since
%$$\pi_* \circ \wt{D}_{(u,J_P)}=\wt{D}^B_{(\mathbf{\pi(u)},J_B)}\circ \pi_*, $$
%we can conclude that every $\eta_i$ must  lie in $L^{q}(\Lambda^{0,1}(S^2,u_{i}^*TP^v))$ and therefore the terms involving the element $Y\in T_{J_B}\mathcal{J}_B$ must vanish in the sum above. 
Next, we show that we can find $Y^v$ and $f$ such that all the components in the sum must be strictly positive unless all the $\eta_v$ are identically zero. Let $Z(u_v)$ denote the set of non-injective points of $u_v$, and consider the subset in $S^2$:
$$X(u_v):=Z(u_v)\cup\bigcup_{v'\in V_+,v'\neq v} u_v^{-1}(u_{v'}(S^2))\cup\bigcup_{v'\in V_0,v'\neq v}u_v^{-1}(u_{v'}(S^2)).$$
Since we consider simple maps, the complement of this set is open dense in $S^2$. Let $x_v$ be a point of the complement  of $X(u_v)$. Then, there is a neighborhood $\mathcal {V}$ of $x_v$ which is embedded via $u_v$ into a neighbourhood $\mathcal{U}_v$ of $u(x_v)$ in $P$, and which does not intersect the image of any other $u_{v'}$. Now, assume $v\in V_+$. From transitivity of the action of Hamiltonian diffeomorphisms on the manifold, we  can find a function $f\in T_H\mathcal{H}=\H$  supported in $\mathcal{U}_v$, such that:
$$\int_{S^2} \LA X_{f(du_v)}^{0,1},\eta_v \RA dvol_{S^2} >0 .$$
When  $v$'s in $V_0$, 
%This time we can find neighbourhoods $\mathcal{U}_i$ in $P$ of $u_i(x_i)$ where $x_i$ is an injective point of $u_i$ in the complement of 
%$$X(u_i):=Z(u_i)\cup\bigcup_{j\in V_+} u_i^{-1}(u_j(S^2))\cup\bigcup_{j\in V_0,j\neq i} u_i^{-1}(u_j(S^2)).$$
 we can also  find  an element $Y^v\in T_J\mathcal{J}^{vert}$ supported in  $\mathcal{U}_v$ and such that
$$\int_{S^2} \LA \frac{1}{2}Y^v. (du_v)^v.j_0,\eta_i \RA dvol_{S^2} >0 .$$
The neighbourhoods $\mathcal{U}_v$  can be chosen small enough so that they are pairwise non intersecting. Set $\left.Y^v\right|_{u_v}=0$  for $v\in V_+$ and $\left.f\right|_{u_v}=0$ for $v\in V_0$. Then $Y^v$ and $f$ are well-defined on the whole manifold  $P$, which implies that all the $\eta_v$'s are vanishing. This ends the proof of the first claim.

\underline{\emph{Proof of (II)}}. Assume $ev_E$ is transversal to $\Delta_{E}$ when $\S_P$ is any labeled forest. Furthermore, suppose that  the graph structure of $\S_P$ is preserved under $S_{\pi_*}$. This is enough since we can always place ourselves in this situation  by adding marked points  in the fiber components so that they are all equipped with at least three special points. This procedure does not alter the transversality for $ev_E$, as the latter does not  depend on the infinitesimal movement of the marked points. Let  $\S_{P}^{+}$ be the stable stratum data resulting from adding the marked points, and consider the map 
$$For^P :\wt{\MMC}_{\mathcal{S}_P^{+}}^{**}(P,\J_P)\sra \wt{\MMC}_{\mathcal{S}_P}^{**}(P,\J_P)$$ that forgets the $k$ added marked points, together with stabilizing the resulting map. Define  $For^B$ in a similar way.
 Then:
$$For^B\circ \wt{\pi}_{\S_P}=\wt{\pi}_{\S_P}\circ For^P.$$
%\begin{equation*}\xymatrix@C=4.0pc{ \wt{\MMC}_{\mathcal{S}(k)}^{**}(P,\J_P)\ar[r]^{\wt{\F}_{\pi}}\ar[d]^{\F_{P,k}} & \wt{\MMC}_{\mathcal{S}_B(k)}^{*}(B,\J_B) \ar[d]^{\F_{B,k}} \\  \wt{\MMC}_{\mathcal{S}}^{**}(P,\J_P)\ar[r]^{\wt{\F}_{\pi}} & \wt{\MMC}_{\mathcal{S}_B}^{*}(B,\J_B)  }
%\end{equation*}
 Clearly  $For^P$ is a submersion. It is not hard to see that $For^B$ is also a submersion. Moreover, using an adaptation to the fibered case of Lemma   3.4.7  in \cite{MS}   one can show that  $\wt{\pi}_{\S_P}$ is also submersion  \cite{Hthesis}. 
 
\underline{\emph{Proof of (III)}}.   The proof proceeds by induction on the number of edges of the  labeled forests $\S_P$. When the forest has no edge, the assertion is vacuous. Suppose it is true for forests with at most $N$ edges, and suppose $\S_P$ is a forest with $N+1$ edges. Pick any edge  given by the pair, $(y_{vv'},y_{v'v})$, cut it out  and replace it by the two new marked points, $y_{vv'}$ and $y_{v'v}$. This procedure gives a new forest $\S_P'$ with two more tails,  which satisfies the induction hypothesis, and such that the sets $I(\S_P)$ and $I(\S_P')$  coincide. Let $E'$ denote the set of edges of $\S_P'$. Then $ev_{E'}$
%$$ev_{E'}:\MMC^{**}(P,\vec{\sigma},\P )\times I(T)\lra P^{E'},$$
is transversal to  $\Delta_{E'}$ so that $\wt{\MMC}^{**}_{\S_P'}(P,\J^r_P)$ is a Banach manifold. Consider the evaluation
\begin{equation*}\label{ev}
ev_{vv'}:\wt{\MMC}^{**}_{\S_P'}(P,\J^r_P)\sra P\times  P,\hspa
    (u,\mathbf{y},\mathbf{x},J_P) \mapsto  (u_v(y_{vv'}),u_{v'}(y_{v'v})).
\end{equation*}
We prove that $ev_{vv'}$ is transversal  to the diagonal $\Delta_P\subset P\times P$. Assume that $\pi$ preserves the tree structure of $\S_P'$ (if not we can add marked points). Let $ev^B_{vv'}$ be the analog of $ev_{vv'}$, but in the case of the base $B$. It is known that $ev^B_{vv'}$ is transversal to $\Delta_B\subset B\times B $ at every point of $\wt{\MMC}^{*}_{S_{\pi_*}(\S_P')}(B,\J^r_B)$ (see \cite{MS}). Furthermore,
 $$ev_{vv'}^B\circ \wt{\pi}_{\S_P'}=(\pi\times\pi)\circ ev_{vv'},$$
Since both $\pi\times\pi$ and  $\wt{\pi}$  are submersions it suffices  to check that: 
\begin{equation*}\label{cokvert}  \forall \wt{u}\in \wt{\MMC}^{**}_{\S_P'}(P,\J_P),\hsp \left.\coker dev_{vv'}(\wt{u})\right|_{\ker d\wt{\pi}(\wt{u})}=0.
 \end{equation*}
By the symmetry  arising from quotienting $TF\oplus TF$ by $T\Delta_F$, it suffices to show that the restriction of
$dev_{vv'}(\wt{u})$ to $V$ surjects onto $T_{u_v(y_{vv'})}F\times\{0\}$. But:
$$W=\Big\{(\{\xi_{v}\}, 0, 0, 0,Y^v,f)\in T_{\wt{u}}\wt{\M}_{\mathcal{S}_P'}^{**}(P,\J_P^r)\,|\, \xi_v\in W^{1,p}(u_v^*TP^v),\hs \forall v\in V \Big\}.
$$
Hence, by definition of $ev_{vv'}$:
$$dev_{vv'}(\wt{u})(\{\xi_v\}, 0, 0, 0,Y^v,f)=(\xi_v(y_{vv'}),\xi_{v'}(y_{v'v})).$$
Now let $(v,0)\in T_{u_v(y_{vv'})}F\times\{0\}$ and suppose the $i$ component is not ghost. Then choose any $\xi_v \in \X^{1,p,v}_{P,u_v}$ such that $\xi_v(y_{vv'})=v$.  Adapting Lemma 3.4.7 in \cite{MS}  to the present situation, 
% which, adapted to our situation 
%says that, given a finite set of pairwise distinct points, $\{z_0,z_1,...,z_r\}\in \Sigma_i$, for every positive $\epsilon$, the set
%$L^p(\Lambda_{J_P}^{0,1}(S^2, u_i^*TP^v))$
% coincides with
% $$\Big\{D_{u_i}^{v,H}\xi_i+\frac{1}{2}Y^v\circ (du_i)^v\circ j \,\,|\,\,\xi_i(z_r)=0,\,\,r\geq 1, \,\,\textrm{and}\,\, \textrm{supp}(Y^v)\subset B_{\epsilon}(u_i(z_0))\Big\},$$
%when $i\in I_0$ and with
%$$\Big\{D_{u_i}^{v,H}\xi_i+X_{f(du_i)}^{0,1} \,\,|\,\,\xi_i(z_r)=0,\,\,r\geq 1, \,\,\textrm{and}\,\, \textrm{supp}(X_f)\subset B_{\epsilon}(u_i(z_0))\Big\},$$
%whenever $i\in I_+$, 
we can find $Y^v$ or $f$ supported in a small enough neighbourhood in $P$ (such that it does not intersect the image of any other component) and a  vector field $\zeta\in W^{1,p}(u_v^*TP^v)$ such that, $\zeta(y_{vv'})=0$ for $vEv'$ and $(\xi_v-\zeta,0,Y^v,0)$, or $(\xi_v-\zeta,0,0,f)$, lies in $\ker \wt{D}_{(u_v,J_P)}$. Then set 
$$\xi_{v'}=0 \hspa\forall v'\neq v.$$
If $u_v$ is ghost, then consider $V_{\textrm{gh}}(v)$ the vertices  of the largest subtree in $\S_P'$, containing $v\in V$ and  consisting only of  ghost components. For all $k\in V_{\textrm{gh}}(v)$ we must have $\xi_v=w$. Consider now all the elements $k\in V\backslash V_{\textrm{gh}}(v)$  such that there exists $v'\in V_{\textrm{gh}}(v)$  for which $kEv'$ and write this set as $K$. All these components have a point in commun in the image of the stable map, namely, $u_v(\Sigma_v)=u_v(y_{vv'})$. Then, for every $m\in K$ choose any $\xi_m \in \X^{1,p,v}_{P,u_m}$ such that $\xi_m(y_{mv'})=w$. Applying the argument used in the non-ghost case to all the components indexed by $K$, we  find vertically valued vector fields $\{\zeta_{m}\}_{m\in K}$ such that 
$$\text{$\zeta(y_{mv'})=0$ when $mEj$, $m\in K$ and $v'\in V_{\text{gh}}(v)$}, $$ 
and such that, for all $m\in K$,  either $(\xi_m-\zeta,0,Y^v,0)$ or $(\xi_m-\zeta,0,0,f)$ lies in the kernel of $\wt{D}_{(u_m,J_P)}$. Finally, set
$\xi_{v'}\equiv 0$ for every  component not indexed by $K\sqcup V_{\textrm{gh}}(v)$.
\qed\\
 
\begin{rem} 
In the proof above,  it is essential that we allow the connection to vary. In particular such perturbations enable us to avoid \emph{horizontal}  $J_P$-holomorphic maps, i.e maps $u$ such that $\textrm{Im}(du)\subset Hor$, for which the index of $D_u^v$ is negative (see \cite{MS}).
\end{rem}

\subsubsection{Cobordisms}
We end this section by stating the invariance of the moduli spaces under changes of the regular structures. Let $\S_P$ be a stable stratum data, and let $\S_B$ be its projection. Given two  regular structures $J_{P}^0$ and $J_{P}^1$ in $\J_{P,reg}(\S_P)$, we designate by  $ \J_P(J_{P}^0,J_P^1)$ the set of paths $\{J_{P}^s\}$ in $\J_P$, $s\in[0,1]$,  with endpoints  $J_{P}^0$ and $J_{P}^1$. Similarly define $\J_B(J_B^0,J_{B}^1)$ for pairs  $J_{B}^0$ and $J_B^1$ in $\J_{B,reg}(\S_B)$. For elements  $\gamma\in  \J_P(J_{P}^0,J_{P}^1)$ and $\gamma_B\in \J_B(J_{B}^0,J_{B}^1)$ we set:
$$\wt{\mathcal{W}}_{\S_P}^{**}(P,\{J_{P}^s\}):=\gamma^*\wt{\MMC}_{\S_P}^{**}(P,\J_P)\hspac \textrm{and}
\hspac \wt{\mathcal{W}}_{\S_B}^{*}(B,\{J_{B}^s\}):=\gamma_B^*\wt{\MMC}_{\S_B}^{*}(B,\J_B).$$

It is not hard to see that if $\gamma$ is transversal to $p^{\J_P}$, and respectively $\gamma_B$ is transversal to $p^{\J_B}$, the quotient under $G_{\S_P}$ and $G_{\S_B}$ of the above pullbacks are then oriented manifolds with boundaries, and with dimensions:
$$\dim(\M_{\S_P}^{**}(P,J_{P}^0))+1\and\dim(\M_{\S_B}^{**}(B,J_{B}^0))+1.$$ In such case we say that the paths are \emph{regular}, and denote by  
$$\J_{P,reg}(\S_P,J_{P}^0,J_{P}^1)\hspac\text{and}\hspac\mathcal{J}_{B,reg}(\S_B,J_{B}^0,J_{B}^1)$$
the set of  such regular paths. Note that a regular path $\gamma$ in $\J_P$ projects to a regular path $\gamma_B:=p_{1}(\gamma)$ in $\J_B$.

For fixed $J_B\in \J_{B,reg}(\S_B)$ and $u_B\in \wt{\MMC}^{*}_{\S_B}(B,J_B)$ set $\J\H((J,H)^0,(J,H)^1)$ to be the set of paths with endpoints  $(J,H)^0$ and $(J,H)^1$ in $\J\H_{reg}(u_B,J_B,\S_P)$, and for any such path $\gamma$ set:
$$\wt{\mathcal{W}}_{\S_P}^{**}(\pi^{-1}(u_B),J_B,\{(J,H)^s\}):=\gamma^*(\wt{\pi})^{-1}(u_B,J_B).$$
If $\gamma$ is transversal to the restriction of $p_{23}\circ p^{\J_P}$  to the fiber $(\wt{\pi})^{-1}(u_B,J_B)$, then the quotient under $G_{\S_P}$ of the corresponding pullback is an oriented manifold with boundary, of dimension $$\dim(\ov{\pi}^{-1}(u_B)\cap \M^{**}(P,J_{P}^0))+1.$$  In such case, we say that $\gamma$ is regular,  and the set of regular paths is denoted by  $\J\H_{reg}(u_B,J_B,\S_P,(J,H)^0,(J,H)^1)$.%

\begin{prop}\label{cobordismgen} The regularizing sets: $\J_{P,reg}(\S_P,J_{P}^0,J_{P}^1)$,  $\mathcal{J}_{B,reg}(\S_B,J_{B}^0,J_{B}^1)$ and $\J\H_{reg}(u_B,J_B,\S_P,(J,H)^0,(J,H)^1)$ are of second category.
 \end{prop}

%%%%%%%%%%%%%%%%%%%%%%%%%%%%%%%%%%%%
%%%%%%%%%%%%%%%%%%%%%%%%%%%%%%%%%%%%

%%%%%%%%%%%%%%%%%%%%%%%%%%%%%%%%%%
%%%%%%%%%%%%%%%%%%%%%%%%%%%%%%%%%%
% Fromule PRODUIT
%%%%%%%%%%%%%%%%%%%%%%%%%%%%%%%%%%
%%%%%%%%%%%%%%%%%%%%%%%%%%%%%%%%%%

\section{The product formula} In this section we establish the product formula. Before doing so, we recall  the definition of Gromov-Witten invariant  for  a semi-positive symplectic manifold $(X,\om)$. For a detailed exposition of the following standard facts, we refer to  \cite{MS} or \cite{RT}.  Let $A\in H_2(X,\Z)$, and  consider the  $l$-pointed  evaluation map:
$$ev^X_{l,J}:\M_{0,l}^*(X,A,J)\sra X^l,\hspace{0.2cm} (u,x_1,\ldots,x_l)\mapsto (u(x_1),\ldots,u(x_l)).$$
This defines a $\dim(\M_{0,l}^{*}(X,A,J))$-pseudocycle of $X^l$ for every $J\in \J_{X,reg}\subset \J_X$, where
\begin{equation}\label{genericsetforGW}\J_{X,reg}:=\bigcap_{\S\in \D^{A}_{0,l}} \J_{X,reg}(\S)
\end{equation}
which is of second category since $\mathcal{D}^{A}_{0,l}$ is finite. Note that $\J_{X,reg}$ depends on $\om$.

We recall  that a  \emph{$d$-dimensional pseudo-cycle} in a manifold $X$, is a pair $(M,f)$, where $M$ is an oriented manifold  $M$ of dimension $d$, and $f:M\to X$ is a smooth map such that the closure,  $\ov{f(M)}$, is compact, and such that its  omega-limit $\Omega_f$ is of codimension at least 2 in $X$. Given classes, $c^X_1,\ldots,c^X_l\in H_{*}(X)$, it is possible to represent them by  pseudo-cycles  $(M_1,f_1),\ldots,(M_l,f_l)$ in $X$, of respective dimensions $\dim M_i:=\deg(c^X_i)$. We can further  assume that these cycles are in general position, and such that
$ev^X_{l,J}$ is strongly transverse to the product cycle, $$\mathcal{C}:=(M_1,f_1)\times...\times(M_l,f_l).$$ Then, the corresponding  \emph{Gromov-Witten invariant} is the algebraic number of isolated points in the  preimage  of  $\mathcal{C}$ under $ev^X_{l,J}$,
\begin{equation*} \la c^X_1,\ldots,c^X_l\ra^{X,J}_{0,l,\sigma}:=  ev^X_{l,J}. \mathcal{C}, 
\end{equation*}
which is  set to be $0$ unless:
\begin{equation*}\label{mathdim} 2n(1-l)+2c_1^{TX}(\sigma)+2l-6+\sum_{i=0}^l \deg(c^X_i)=0.
\end{equation*} 
This number only depends on the  bordism class of the  pseudo-cycles involved. In particular, it does not depend on the  regular almost complex structure, which we will drop from  the notations. 

Now, let $\pi:P\sra B$ be a Hamiltonian fibration with coupling form $\tau$, and let $\iota_F^P$ denote the inclusion of $F$ in $P$. Consider $\sigma\in H_2(P,\Z)$ with $\sigma_B=\pi_*\sigma\neq 0$.  For $(u_B,\mathbf{x})\in \MMC^{*}_{0,l}(B,\sigma_B,J_B)$,
 we have the commutative diagram:
\begin{equation}\label{evdiagram}
\xymatrix{\ov{\pi}^{-1}(u_B,\mathbf{x})/G \ar[r]\ar[d]^{ev_{(u_B,\mathbf{x})}}&\MMC^{**}_{0,l}(P,\sigma,J_P)\ar[d]^{ev^P_{l,J_P}}\ar[r]^{\ov{\pi}}& \MMC^{*}_{0,l}(B,\sigma_B,J_B) \ar[d]^{ev^B_{l,J_B}}\\
F^l\ar[r]^{(\iota_F^P)^l}&P^l\ar[r]^{\pi^l} & B^l}
\end{equation}
where 
$$ev_{(u_B,\mathbf{x})}:\ov{\pi}^{-1}(u_B,\mathbf{x})/G\lra F^l,\hspace{0.5cm} u\mapsto u(\mathbf{x}):=(u(x_1),...,u(x_l))\in\prod_{i=1}^lF_{u_B(x_i)}.$$
The product formula is obtained by considering the (respective) intersections of $ev_{(u_B,\mathbf{x})}$, $ev^B_{l,J_B}$, and $ev^P_{l,J_P}$, with the  product pseudo-cycles: 
\begin{equation*} (\mathcal{C}^F,f^F):=\prod_{i=1}^{l}(M^F_i,f^F_i),\hspa(\mathcal{C}^B,f^B):=\prod_{i=1}^{l}(M^B_i,f^B_i),\hspa (\mathcal{C}^P,f^P):=\prod_{i=1}^{l}(M^P_i,f^P_i),
\end{equation*}
where, $(M^F_i,f^F_i)$, $(M^B_i,f^B_i)$, and $(M^P_i,f^P_i)$,  respectively represent torsion free  homology classes,  $c_i^F$, $c_i^B$, and $c_i^P$,  verifying condition \eqref{cond11}. In particular,
$$\begin{cases}(M^B_i, f^B_i)=(pt, f^B_i) & \text{for  $i=1,...,m$}\\
 (M^F_i,f^F_i)=(F,id_F)& \text{for  $i=m+1,...,l$.}
\end{cases}$$ 
Furthermore,  a $d$-dimensional pseudo-cycle $(M,f)$ in the fiber $F$ of $P$ defines a $d$-dimensional pseudo-cycle $(M,\iota_F^P\circ f)$ in the total space. Similarly, any $d$-dimensional  pseudocycle  $(M,f)$ in $B$  defines a $d+\dim F$ pseudocycle $(f^*P,\ov{f})$ in $P$, where   $\ov{f}$ stands for the bundle map associated to $f$. These operations actually preserve the bordism classes.  We conclude that:
$$(M^P_i,f^P_i)=\begin{cases} (M^F_i,\iota_F^P(f^F_i)) &  \text{if $i=1,...,m$}\\
((f^B_{i})^*P,\ov{f}^B_{i}) & \text{otherwise.}
\end{cases} $$ 
 Regarding orientations of the product pseudo-cycles,  the exact sequence,
\begin{equation*} 0\lra df^F(T\CC_F)\stackrel{(\iota_F^P)^l}{\lra} df^P(T\CC_P)\stackrel{\pi^l}{\lra} df^B(T\CC_B)\lra 0,
\end{equation*}
gives:
\begin{equation*}
 \det df^P(T\CC_P)\cong \det df^B(T\CC_B)\otimes \det df^F(T\CC_F).
\end{equation*}
Therefore, if we choose the cycles $(\CC_B,f^B)$ and $(\CC_F,f^F)$ to be positively oriented, the cycle $(\CC_P,f^P)$ must also be positively oriented. Now, assume the evaluations are pseudo-cycles, and that strong transversality with the product cycles is achieved. Then $(ev_l^B)^{-1}(f^B)$ is a finite set, $\{(u_{B,\alpha},\mathbf{x}_{\alpha})\}_{\alpha\in A}$, of isolated simple, $l$-pointed, $J_B$-holomorphic maps.
For each $\alpha$, let $\iota^P_{F,\alpha}$ denote the embedding of   $F^l$ into $F^l_{u_{B,\alpha}(\mathbf{x}_{\alpha})}$. Also, in order to simplify   notations, set
 $$ev_{\alpha}:=ev_{(u_{B,\alpha},\mathbf{x}_{\alpha})},\hsp f^F_{\alpha}:=\iota^P_{F,\alpha}\circ f^F,$$ 
and write:
$$n_{\alpha}:=ev_{\alpha}. f^F_{\alpha},\hsp n_B:=ev^B_l.f^B,\hsp\text{and}\hsp n_P:= ev^P_l.f^P.$$ 
Notice that the numerical conditions under which $n_{\alpha}$ and $n_B$ are non necessarily zero, provide the condition under which  $n_P$ is possibly non vanishing. Now, in the above notations, the product formula now reads:
\begin{equation}\tag{\textbf{PF}}\label{relationprincipale}\forall \alpha\in A,\hspa n_P=n_{\alpha}n_B.
\end{equation}

We will prove this relation, and then prove the Corollary in the last subsection,  Before doing so, we  make sure that the evaluation maps in \eqref{evdiagram}, are simultaneously   pseudo-cycles  that remain in the same  bordism class under change of regularizing almost complex structure. This will give meaning to all the numbers in \eqref{relationprincipale}.  
 
\subsection{Evaluation maps as pseudo-cycles.}  We begin by fixing some notations. Define $\J_{B,reg}\subset \J_B$ and $\J_{P,reg}\subset \J_P$ as in \eqref{genericsetforGW}. 
%(Note that $\J_{P,reg}$ depends on $\om$, $\om_B$ and $\tau$).
These sets are  of second category.  For  $J_B\in \J_{B,reg}$ and  $u_B\in \MMC^*(B,\sigma_B, J_B)$, we set,
 $$\J\H_{reg}(u_B,J_B):=\bigcap_{\left\{\S_P| S_{\pi_*}(\S_P)=\S^{top}_B\right\}}\J\H_{reg}(u_B,J_B,\S_P),$$
where $\S^{top}_B$ denotes the top stable stratum data for  $\ov{\MMC}(B,\sigma_B, J_B)$, i.e the stratum having only one vertex as a tree structure. By Theorem B', this set is also of second category.  Furthermore, for fixed $J_B$, we say that the class $\sigma_B\in H_2(B,\Z)$ only  admits \emph{irreducible effective decompositions with respect to $J_B$} if  every stratum  $\MMC_{\S_B}(B,J_B)$  is only made of irreducible elements. 
This condition is in particular realized by primitive classes,  for example the class of a line in $\CP^n$, or  the diagonal in  $S^2\times S^2$ with the standard product complex structure.  Let  $\J_{irr}(\sigma_B)$ denote the subset of $\J_B$ with respect to which $\sigma_B$  admits only irreducible effective decompositions. Then  $\J_{irr}(\sigma_B)$ is open in $\J(B,\omega_B)$. Nevertheless, nothing garantees that   it is non-empty.  
In the  theorem below, the restriction to $\J_{irr}(\sigma_B)$ is essential in order to avoid simple stable maps having a reducible projection. 

We  show that all the evaluation map in the above diagram are  pseudo-cycles. Note that   condition \eqref{ssp} is equivalent to: 
\begin{equation}\label{semipositivity}\forall A\in \pi_2(F):\hsp  \omega(A)>0\hsp, c^{TF}_1(A)\geq 3-n_P\hspace{0.5cm}\Longrightarrow c^{TF}_1(A)\geq 0.
\end{equation}
This is  weaker than asking for $P$ to be semi-positive. However, this implies that the fiber is semi-positive.

\begin{theorem}\label{evpc} Assume \eqref{semipositivity} and that $\J_{irr}(\sigma_B)\neq\emptyset$. Then:
 \begin{enumerate}[i)] 
\item For every $J_P\in \J_{P,reg}$ with $J_B\in \J_{irr}(\sigma_B)$, the evaluation maps, $ev^B_{l,J_B}$ and $ev^P_{l,J_P}$, are pseudo-cycles. Moreover, changing of regular structure along regular path induces a bordism between the relevant evaluation maps, as long as  the almost  complex structure on  $B$  varies in a  connected component of $\J_{irr}(\sigma_B)$.
\item Fix a regular structure $J_B$, and let $(u_B,\mathbf{x})\in\MMC_{0,l}^*(B,J_B,\sigma_B)$. Then, for any element in  $\J\H_{reg}(u_B,J_B)$, the couple $(\ov{\pi}^{-1}(u_B,\mathbf{x}),ev_{(u_B,\mathbf{x})})$ is a pseudo-cycle that remains in the same  bordism class, under change of regularizing pair along regular paths.
\end{enumerate}
\end{theorem}

\proof We only show the first statement, the proof of the second being similar. Fix $J_P\equiv(J_B,J,H)\in \J_{P,reg}$. Hence $J_B\in\J_{B,reg}$ and  $ev^B_{l,J_B}$ is then a pseudo-cycle (\cite{MS}). By Lemma \ref{strata}, there are only  finitely many stable stratum datas $\S_P$ representing geometric limits of curves in $\MMC_{0,l}^{**}(P,\sigma,J_P)$. Hence, by Gromov's compactness,
$$\Omega_{ev^P_{l,J_P}}\subset\bigcup_{\S_P}ev^P_{l,J_P}(\MMC^{**}_{\S_P}(P,J_P)),$$
where the union is taken over all reduced stratum datas. Let $\beta^{red}_P$ denote the homological labeling associated  to a reduction $\S_P^{red}$  of $\S_P$.  Let $V^{red}_+$ denote the set of $\pi$-stable components in $\S^{red}_P$, and let $V_0^{red}$ denote the set of $\pi$-unstable components in $\S_P^{red}$. Since, by assumption,  $\sigma_B$ only admits irreducible decompositions, only $\pi$-unstable components of $\S_P$  are contracted in the reduction process. Therefore,  there exist integers $m_v>0$ for all $v\in V^{red}_0$ such that:
$$\sigma=\sum_{v\in V^{red}_+}\beta^{red}_P(v) +\sum_{v\in V^{red}_0}m_v\beta^{red}_P(v).$$
Condition \eqref{semipositivity} further implies  that  $c^v(\beta^{red}_P(v))>0$, for every $v\in V^{red}_0$, and we conclude that
 $$\dim\MMC^{**}_{\S^{red}_P}(P,J_P)\leq\dim \MMC_{0,l}^{**}(P,\sigma,J_P)-2$$
for every stratum reduction. Hence, $ev^P_{l,J_P}$ is a pseudo-cycle. The independance statement is shown as follows. Let $J_{P}^t$ be a regular  path of fibered almost complex structures between regular fibered structures, projecting on a path in $\J_{B,reg}(J_{B}^0,J_{B}^1)\cap\J_{irr}(\sigma_B)$ (with $J_{B}^0,J_{B}^1$ in the same connected component of $\J_{irr}(\sigma_B)$). 
Then, any  Gromov limit of a sequence  in  $\mathcal{W}_{0,l}^{**}(P,\sigma,\{J_{P}^t\})$,
is such that its non-trivial roots are irreducible, while its fiber components may actually be reducible. Argumenting exactly as above, the lower strata in $\ov{\mathcal{W}}_{0,l}^{**}(P,\sigma,\{J_{P}^t\})$ have codimension at least 2 in $P$. Thus,
$$ev^P_{l,\{J_{P}^t\}}:\mathcal{W}_{0,l}^{**}(P,\sigma,\{J_{P}^t\})\sra P^l $$
is a  pseudo-cycle inducing a bordism between  $ev^P_{l,J_{P}^0}$ and $ev^P_{l,J_{P}^1}$.
\qed

\begin{rem} Note that in this context we do not need to impose any semi-positivity assumption on $B$ due to the specific decomposability hypothesis imposed on $\sigma_B$. Note also that if  $\sigma_B$ is undecomposable, we can drop the restriction on $\J_{irr}(\sigma_B)$.
\end{rem}

\subsection{Proof of theorem A} We begin by proving that all the terms in \eqref{relationprincipale} are well-defined.  From Theorem \ref{evpc} and since $\J_{irr}(\sigma_B)\neq\emptyset$,  the evaluation maps  $ev^B_l$ and $ev^P_l$  generically define pseudo-cycles. Choose (generically) the cycles $(M^B_i, f^B_i)$ so that  $ev^B_l$ is strongly tranvserse to $(\CC_B,f^B)$. Then $(ev^B_l)^{-1}(f^B)$ is finite and, as already mentioned, the corresponding GW-invariant, $n_B$, only depends on the bordism classes of $(M^B_i,f^B_i)$, and on the connected components of $\J_{irr}(\sigma_B)\cap\J_{B,reg}$. Thus, 
$$(ev^B_l)^{-1}(f^B)=\{(u_{B,\alpha},\mathbf{x}_{\alpha})|\alpha\in A\},$$ 
is finite,  and for every $\alpha\in A$, the map  $ev_{\alpha}$  also  defines a pseudocycle for generic fiber regularizing pairs. Consequently,  $(F^l,f^F_{\alpha})$ 
  is a pseudo-cycle of $F^l_{u_B(\mathbf{x}_{\alpha})}$ for every $\alpha\in A$.  Since  $A$ is finite, we can furthermore choose the cycles $(M^F_i,f^F_i)$  such that, for every $\alpha\in A$, the  evaluation $ev_{\alpha}$ is transversal to $f^F_{\alpha}$. This, together with the fact that $ev^B_l$ is transversal to $f^B$ implies that $ev^P_l$ is  transversal to $f^P$. The independance of  $n_{\alpha}$ and $n_P$, with respect to the choice of regularizing triple, follows from Theorem  \ref{evpc}.
  
 Next, we prove   \eqref{relationprincipale}. Let $C_{\alpha}\subset B$ denote the image of $u_{B,\alpha}$,  let  $P_{C_{\alpha}}$ be the restriction of $P$ to $C_{\alpha}$, and let, $\iota_{\alpha}:P_{C_{\alpha}}\hookrightarrow P$ and $\iota^{P_{C_{\alpha}}}_{F}:F\hookrightarrow P_{C_{\alpha}} $, denote the natural inclusions.  Consider the subset of section classes in $P_{C_{\alpha}}$,
$$B_{\sigma}^{\alpha}=\{\sigma'\in H_2(P_{C_{\alpha}},\Z)|(\iota_{\alpha})_*\sigma'=\sigma\}.$$ 
Then, 
\begin{equation}\label{relationfibre}n_{\alpha}=\sum_{\sigma'\in B_{\sigma}^{\alpha}}\la \iota_F^{P_{C_{\alpha}}}(c^F_1),...,\iota_F^{P_{C_{\alpha}}}(c^F_l)\ra_{0,l,\sigma'}^{P_{C_{\alpha}}}.
\end{equation}
%actually gives  the algebraic number  of $J_{P_{C_{\alpha}}}$-holomorphic maps intersecting the  cycles $\iota^P_{F,\alpha}\circ f^F_i$ in $F_{u_{B,\alpha}(x_{i,\alpha})}$ at the points $x_{i,\alpha}\in S^2$, and representing the classes $\sigma'_{\alpha,j}$. 
Indeed, let  $J_{P_{C_{\alpha}}}$ denote the restriction of $J_P$ to $P_{C_{\alpha}}$, then  $\iota_{\alpha}$ naturally induces an identification:
$$\ov{\iota}_{\alpha}: \bigsqcup_{\sigma'\in B_{\sigma}^{\alpha}}\ov{\MMC}(P_{C_{\alpha}},J_{P_{C_{\alpha}}},\sigma')\sra \ov{\pi}^{-1}(u_{B,\alpha},\mathbf{x}_{\alpha}),$$
which is an orientation preserving diffeomorphism when restricted to any stratum. Furthermore, by simplicity of   $u_{B,\alpha}$,  the $l$ marked points are naturally identified  to the $2l$-dimensional manifold $\MMC^*_{0,l}(C_{\alpha},[C_{\alpha}])$. We obtain the following  diagram:
\begin{equation*}\xymatrix{\ov{\pi}^{-1}(\mathbf{x}_{\alpha})=\bigsqcup_{\sigma'\in B_{\sigma}^{\alpha}}\MMC^*(P_{C_{\alpha}},\sigma')\ar[r]\ar[d]^{ev_{\mathbf{x}_{\alpha}}^{P_{C_{\alpha}}}}&\bigsqcup_{\sigma'\in B_{\sigma}^{\alpha}}\MMC_{0,l}^*(P_{C_{\alpha}},\sigma')\ar[r]^{\ov{\pi}}\ar[d]^{ev_l^{P_{C_{\alpha}}}}&\MMC^*_{0,l}(C_{\alpha},[C_{\alpha}])\ar[d]^{ev_l^{C_{\alpha}}} \\ F^l\ar[r]^{(\io_{F}^{P_{C_{\alpha}}})^l} & P_{C_{\alpha}}^l\ar[r]^{\pi^l}& C_{\alpha}^l }
\end{equation*}
where the complex structures are dropped in order to simplify notations.  By definition,  $ev^{P_{C_{\alpha}}}_{\mathbf{x}_{\alpha}}$ is the composition of $ev_{\alpha}$ with $\iota_{\alpha}$, hence is  a pseudo-cycle for generic fiber regularizing parameters.  Using the above diagram we  conclude that $ev_l^{P_{C_{\alpha}}}$ is generically a  pseudo-cycle. Then, equation  \eqref{relationfibre} follows  since there is only one positively oriented curve in $\MMC^*_{0,l}(C_{\alpha},[C_{\alpha}])$ intersecting transversally $l$ points at the $l$-marked points (giving $\la pt,...,pt \ra_{0,l,[C_{\alpha}]}^{C_{\alpha}}=1$).

%Once again, notice that the conditions under which $n(u_{B,\alpha},\mathbf{x}_{\alpha})$ and $n^B$ are non necessarily zero, provide the condition under which  $n^P$ is possibly non vanishing. We subsequently show that for every $\alpha\in A$ the following equality holds:
%\begin{equation}\label{relationprincipale}n^P=n(u_{B,\alpha},\mathbf{x}_{\alpha})n^B.
%\end{equation}

Now, consider the sign functions,  $\epsilon_P$, $\epsilon_B$, and $\epsilon_{\alpha}$ respectively associated to the curves counted  in $n_P$, $n_B$ and  $n_{\alpha}$.  We have to make sure that the signs of the counted curves are given  compatibly with $\pi^l$, i.e that:
 $$\forall \alpha\in A,\hspa \epsilon_P=\epsilon_B\times\epsilon_{\alpha}.$$
  But this is the case since
% First note that the orientation of $f^{P}$ is induced by those of   $f^B$ and $f^F_{\mathbf{x}}$. 
%Nous pouvons de plus nous  arranger afin que ces cycles  soient positivement orient\'es. 
for every fiber regularizing pair %element in  $\J\H_{reg}(A,\sigma,J_B)$ 
we have the exact sequence:
\begin{equation*} 0\sra T_{u}\left(\ov{\pi}^{-1}(u_{B,\alpha},\mathbf{x}_{\alpha})\right)
\sra T_{(\io_{\alpha}(u),\mathbf{x}_{\alpha})}\M^{**}_{0,l}(P,\sigma)\sra T_{(u_{B,\alpha},\mathbf{x}_{\alpha})}\MMC^*_{0,l}(B,\sigma_B)\sra 0.
 \end{equation*}
%Each term of this sequence is naturally oriented  so that:
%\begin{equation*}\det  T_{(\io^P_{P_{C_{\alpha}}}(u),\mathbf{x}_{\alpha})}\M^{**}_{0,l}(P,\sigma;J_P) \cong T_{u}\left(\pi^{-1}(u_{B,\alpha},\mathbf{x}_{\alpha})\right)\otimes T_{(u_{B,\alpha},\mathbf{x}_{\alpha})}\MMC^*_{0,l}(B,\sigma_B;J_B).
% \end{equation*}
%Consequently, the relation $\epsilon_P=\epsilon_B\times\epsilon_{P_{C_{\alpha}},\mathbf{x}_{\alpha}}$ is verified for all $\alpha\in A$ and we have: 
Consequently:
\begin{eqnarray*}
n_P&=&\sum_{\alpha\in A}\sum_{\{u\in ev^{-1}_{\alpha}(f^F_{\alpha})\}}
\epsilon_P(\io_{\alpha}(u))\\ 
%&=&\sum_{\alpha\in A }\sum_{\{u\in n_{\alpha}\}} \epsilon_{\alpha}(u)
%\epsilon_B(u_{B,\alpha},\mathbf{x}_{\alpha}) \\
&=&\sum_{\alpha\in A}\left(\epsilon_B(u_{B,\alpha},\mathbf{x}_{\alpha})\sum_{\{u\in ev^{-1}_{\alpha}(f^F_{\alpha})\}} \epsilon_{P_{C_{\alpha}},\mathbf{x}_{\alpha}}(u)\right)
\end{eqnarray*}
which, by Lemma \ref{curveindependance} coincides with $n_{\alpha_0}n_B$ for any $\alpha_0\in A$.
\qed\\

%\subsection{Proof of corollary: $c$-splitting and symplectic uniruledness} 

\section{Gluing and  fibration of moduli spaces}
The aim of this section is to show that, under some circumstances, $\pi$ defines   a locally trivial smooth orbi-fibration above the top stratum of $\ov{\M}_{0,l}(B,\sigma_B,J_B)$, with $\sigma_B\neq 0$.
We start with the following simple observation. From the transversality theorem,  the restriction of $\ov{\pi}_{\S_P}$ to $\M^{**}(P,\sigma,J_P)$ is  generically a smooth  submersion onto  $\M^{*}(B,\sigma_B,J_B)$ over countably many points. It is natural to ask if there exists  a fibered  $J_P$ with respect to which the latter map is everywhere regular.
For fixed $J_B$ and $J$, we say that $H\in \mathcal{H}$ is \emph{parametric} if  $D_u^{v}$ is surjective  for every $u\in \M^{**}(P,\sigma,J_P)$.
From exactness of \eqref{splitting}  $\pi$ is a smooth submersion for parametric $H$. 
As a result, if we  assume that $\sigma$ is an undecomposable effective class projecting onto a non-trivial undecomposable class in $B$ and that $H$ is parametric, then $\ov{\pi}_{\S_P}$ is a (smooth) locally trivial fibration  (Indeed, since $\sigma$ is undecomposable, $\M^{**}(P,\sigma,J_P)$ is compact so that  $\pi_{\S_P}$ is proper, as desired).

\begin{rem}
The set of parametric $H$'s may  be empty. Consider the $(\CP^1,\om_{FS})$-fibration, $$\pi:E:=\mathbb{P}(\mathcal{O}_{\CP^2}(-2)\oplus\mathbb{C})\lra(\CP^2,\om_{FS}).$$
Denote by $L$ the homology class of a line in $\CP^2$, and let  $L_0\in H_2(E,\mathbb{Z})$ be the class such  that  $L=\pi_*(L_0)$ and $L_0\cap [F]=[pt]$, where $[F]$ stands for the class of a fiber in $E$. If  $u$ is a holomorphic curve in  $E$ representing $L_0$, the fibration $u^*TE$ is isomorphic to the direct  sum of $\mathcal{O}_{\CP^1}(2)$, $\mathcal{O}_{\CP^1}(1)$ and $\mathcal{O}_{\CP^1}(-2)$.
A straightforward computation then shows that the index of the vertical linearized operator must be  $-2$, hence there are no parametric $H$.
\end{rem}

We will make the following assumption throughout this section:

\begin{assumption}\label{strongreg}
Let $\S_P$ be  a stable stratum data for pseudo-holomorphic maps in $P$ representing  $\sigma$. We ask that for every  $u\in \wt{\M}_{\S_P}(P,J_P)$:
\begin{itemize}
\item  the operators $D^v_{u}$ and $D^B_{\pi(u)}$ are surjective;
\item  the edge evaluation maps $ev^B$ and  $ev_{\pi_{\S_P}(u)}$ are transversal to the corresponding diagonals.
\end{itemize}
\end{assumption}

Then for every $u\in \M_{\S_P}(P,J_P)$ the operator $D^P_u$ is surjective, and  $ev^P$ is transverse to the associated diagonal (see  Section 3). 
%Furthermore, if  $\S_B$ denotes the image of $\S_P$ under $\pi$,  condition (\ref{strongreg}) ensures that the moduli spaces  $$\wt{\M}_{\S_P}(P,J_P),\hspa \wt{\M}_{\S_B}(B,J_B)\and \wt{\pi^{-1}(u_B)}$$ are naturally oriented manifolds (for every $u_B\in \wt{\M}_{\S_B}(B,J_B)$).
By a  standard argument  the quotient spaces   $\M_{\S_P}(P,J_P)$, $\M_{\S_B}(B,J_B)$, and the fiber of $\ov{\pi}_{\S_P}$ over $\ov{\pi}_{\S_P}(u)$,  are smooth orbifolds. Furthermore, the commutativity  $\pi_*\circ D^P_u=D^B_{\pi_{\S_P}(u)}\circ \pi_*$  implies that these orbifold structures can be chosen compatibly with  $\ov{\pi}_{\S_P}$. 
We will drop the almost complex structures from the notations since  it is understood that we made a choice here.

%%%%%%%%%%%%%%%%%%%%%%%%%%%%%%%%%%%%%%%%%%%%%%%%%%%

                                       % Gluing COURBES NODALES
                                       
    %%%%%%%%%%%%%%%%%%%%%%%%%%%%%%%%%%%%%%%%%%%%%%%%%                                   

\subsection{Gluing in the non-fibered case: $B=pt$.}  We give the gluing procedure in $\ov{\M}_{0,l}(X,A)$ with $(X,\om)$ a general symplectic manifold. This will be our guideline when considering more general base. We start by gluing in $\M_{0,l}$.  We follow standard approaches in the litterature, such as \cite{MS}, \cite{Sik}, \cite{CL} amongst others. 

\subsubsection{Gluing for nodal curves}  
%We refer to ~\cite{HM} for proofs in the following section.
  Let, $\S:=(V,\text{Fl};\text{pr},\varrho)$, be a stable stratum data for $\M_{0,l}$, and let  $\j\equiv (\Sigma,\j,\mathbf{x})\in\M_{\S}$. For $v\in V$, let $\Sigma_v$ denote the corresponding (irreducible) component, and for $f\in \text{Fl}$, denote by $z_f$ the corresponding special point on $\Sigma_v$.  By definition each $\Sigma_{v}$ is stable.  Then, upto isometry, there exists a unique  isometric action of a Fuchsian group $\Gamma$ on the hyperbolic half plane $\mathbb{H}$ with respect to which:
 $$\Sigma_v\setminus \{z_f|\text{pr}(f)=v\}\cong  \left.\mathbb{H}\right/_{\Gamma}.$$
 The induced metric belongs to the conformal class given by the complex structure on $\Sigma_v$.  If  $D\subset \mathbb{H}$ denotes  a Dirichlet region of $\Gamma$, then each $z_f$ corresponds to a  vertex at infinity, and we can choose the fundamental region $D$ such that  $z_f$ corresponds to infinity with edges  $x\equiv 0$ and $x\equiv 1$. It is well known that given a real $b>1$, the  \emph{horocycle} at $z_f$, 
  $$\left.\{x+iy\in \mathbb{H}| 2\pi y> b\}\right/_{\la z\mapsto z+1\ra},$$ 
defines a neighbourhood of $z_f$, which can be identified to a punctured disc $D^*(e^{1-b})\subset \C$ via the map $z\mapsto e^{2\pi i z+1}$. We have such neighbourhoods for each $z_f$, and we denote it $D_f^*(r_f)$, with small $r_f>0$. Now, $D_f^*(r_f)$ is conformally equivalent to a cylinder with \emph{negative end}:
$$(-\infty, \ln r_f]\times \R/2\pi \Z\cong D_f^*(r_f): (s,t)\mapsto e^{s+it}$$
or with \emph{positive end}:
 $$[-\ln r_f,+\infty) \times \R/2\pi \Z\cong D_f^*(r_f): (s,t)\mapsto e^{-s-it}.$$

Let  $e_{f,f'}$ be an edge  of $\S$ between the vertices $v:=\text{pr}(f)$ and $v':=\text{pr}(f')$. For a  complex number  
\begin{equation}\label{gluingparameter}\rho_{vv'}:=e^{-R_{vv'}+i\theta_{vv'}}\equiv r_{vv'}e^{i\theta_{vv'}},\hspa \theta\in [0,2\pi),
\end{equation}
 such that $|\rho_{vv'}|< \min\{r_f,r_{f'}\}$, we can glue the components $\Sigma_{v}$ and $\Sigma_{v'}$ at $z_f$ and $z_{f'}$  as follows. Put positive (resp. negative) cylindrical coordinates on $D_f^*(r_f)$ (resp. $D_{f'}^*(r_{f'})$ and  identify the annuli:
\begin{eqnarray*}[-R_{vv'}/2-1,-R_{vv'}/2+1]\times \R/2\pi\Z &\stackrel{\cong}{\sra}& [R_{vv'}/2-1, R_{vv'}/2+1]\times \R/2\pi\Z, \\
(s,t)&\mapsto& (s+R_{vv'},t+\theta_{vv'})
\end{eqnarray*}
 This gives the patching procedure between  the annuli $\A_f:=D_f^*(r_f)\bs D_f^*(e^{-R_{vv'}})$ and $\A_{f'}:=D_{f'}^*(r_{f'})\bs D_{f'}^*(e^{-R_{vv'}})$. The resulting curve $\Sigma_{\rho_{vv'}}(f,f')$ has  a natural complex structure and is the gluing of $\j$ at $e_{f,f'}$ with parameter $\rho_{vv'}$.  
 
 Note that $\rho_{vv'}$ can be naturally identified to   an element of  $$\C_{vv'}:=T_{y_{vv'}}\Sigma_v\otimes T_{y_{v'v}}\Sigma_{v'}\cong\C,\hspa vEv'.$$
Denoting by   $\C_{\j}$ the direct sum of the $\C_{vv'}$ over all edges in $\S$, and by $B_{\epsilon}$ the ball of radius $\epsilon $ at the origin of  $\C_{\j}$,
then for small enough $\epsilon>0$,  the gluing  gives a map:
\begin{equation*}\label{pointwisegluing}gl_{\j}:B_{\epsilon}\subset \C_{\j}\sra \ov{\M}_{0,l},\hspa \rho \mapsto \j_{\rho}\equiv \Sigma_{\rho}.\end{equation*}
 which coincides with the identity map when  $\rho\equiv 0$. Note that  $gl_{\j}$ is  $\text{Aut}(\j)$-equivariant (for the linear action of $\aut(\j)$ on $\C_{\j}$), and it follows from stability of the curve that  $gl_{\j}$ is  injective. Taking the union over $\M_{\S}$ of   the $\C_{\j}$  actually defines an orbibundle $$p_{\S}:\L_{\S}\sra \M_{\S}.$$
Let $\mathcal{L}_{\S,\epsilon}$ denote the restriction to an $\epsilon$ neighbourhood of the zero section, and let $\L_{\S}^*$ be $\L_{\S}$ with the zero section removed.
Given  a proper open subset   $U\subset \M_{\S}$,   there exists $\epsilon>0$, depending on  $U$, such that the above gluing map extends to a locally diffeomorphic map:

\begin{equation*}\label{gluingdomain}gl_{\S}:\mathcal{L}^*_{\S,\epsilon,U}:=\left.\mathcal{L}^*_{\S,\epsilon} \right|_{U}\sra \ov{\M}_{0,l},\hspa (\j,\rho)\mapsto gl_{\j}(\rho).\end{equation*}

More generally, given two stable  stratum datas such that  $\S\prec\S'$, there exists a subbundle  $\mathcal{L}_{\S,\S'}$ of $\mathcal{L}_{\S}$ with fibers identified to $\C^{|E'|-|E|}$, $|E'|$  being the number of edges in $\S'$, as well as gluing maps (defined on proper open subsets) $$gl_{\S,\S'}:\mathcal{L}^*_{\S,\S',\epsilon}\sra \ov{\M}_{\S}.$$

When $\S$ is unstable, one can still define $\L_{\S}$ over $\M_{\S}$, as well as a gluing map $gl_{\S}$, but $gl_{\S}$ is injective and locally diffeomorphic if and only if $\S$ is stable.

\subsubsection{Gluing stable components} Let $\S_X:=(V,\text{Fl};\text{pr},\varrho)$ be a stable stratum data. Here we  assume that the forgetful map $\S:=S_{\pi_{pt_*}}(\S_X)$ and $\S_X$ have the same tree structure.  Then $\ov{\pi}_{pt}^{\S_X}$ coincides with the \emph{forgetting-the-map map}:  
$$\F_X: \M_{\S_X}(X)\sra \M_{\S},\hsp (u,\j)\mapsto \j.$$
 Let  $\mathcal{L}_{\S_X}$ denote  the orbibundle $\F_X^*\mathcal{L}_{\S}$ over $\M_{\S_X}(X)$.  In local coordinates, an element in  $\mathcal{L}_{\S_X}$ is given by a triple $(u,\j,\rho)$. Consider the bundle map:
  $$\ov{\ov{\pi}_{pt}}^{\S_X}:\mathcal{L}_{\S_X}\sra \L_{\S}, \hsp(u,\j,\rho)\mapsto (\j,\rho).$$ 
   The gluing consists in  constructing, compatibly with  $\ov{\ov{\pi}_{pt}}^{\S_X}$,  a $J$-holomorphic map with domain $gl_{\S}(\j,\rho)$ out of  any $(u,\j,\rho)$. More precisely,  the following holds:

\begin{theorem}\label{gluingtheorem}(\cite{CL}, \cite{MS}) For every proper subset $U_X\subset \M_{\S_X}(X)$, there exists a constant, $\epsilon_X>0$, and a locally  diffeomorphic map:
$$Gl_{\S_X}: \mathcal{L}^*_{\S_X,\epsilon_X,U_X}:=\left.\mathcal{L}^*_{\S_X,\epsilon}\right|_{U_X} \sra \M_{0,l}(X,A),\hsp (u,\j,\rho)\mapsto Gl_{\S_X}(u,\j,\rho), $$
such that   
\begin{equation}\label{commutativitygluing}\ov{\pi}_{pt}\circ Gl_{\S_X}=gl_{\S}\circ \ov{\ov{\pi}_{pt}}^{\S_X}.
\end{equation}
\end{theorem}

\proof We give a sketched proof of this theorem which serves as a guiding principle when gluing pseudo-holomorphic curves. This is done in several steps.\\

\noindent\emph{\underline{Step 1}: pregluing.} Here we  construct an approximatively  $J$-holomorphic smooth map out of    $(u,\j,\rho)\in \L^*_{\S_X}$.  For $v\in V$ let $u_v:\Sigma_v\sra X$ denote the corresponding component of $u$. Also, for an edge $e_{f,f'}$ of $\S_X$ with $\text{pr}(f)=v$ and $\text{pr}(f')=v'$, let $y_{v,v'}$ denote the corresponding node on $\Sigma$. Let $\rho_{vv'}\in\C^*$, defined by \eqref{gluingparameter},
 denote the gluing parameter associated to $y_{v,v'}$. Furthermore, let $\beta:\mathbb{R}\sra [0,1]$ denote a smooth cut-off  function with uniformily bounded derivative, $|\beta'(r)|\leq 2$, such  that $\beta(r)=0$ if  $r\leq1$ and $\beta(r)=1$ for $r\geq2$.
 The \emph{pregluing of $(u,\j)$ with parameter $\rho$}, is the smooth map, 
  $$u_{\rho}\equiv \textrm{pgl}^X(u,\j,\rho):\Sigma_{\rho}\sra X,$$   defined as follows: 
for each $v\in V$, for every $v'\in V$ such that $vEv'$,   
\begin{equation*}u_{\rho}(z):=\begin{cases}
  %  u_v(z)  & \text{if} \hs z\in \Sigma_v-D_{f}(2r_{vv'}^{1/4}), \\
     p_{vv'}:=u_v(y_{v,v'})=u_{v'}(y_{v,v'}) & \text{if} \hs z\in D_{f}(r_{vv'}^{1/4})\bs D_{f}(r_{vv'}^{3/4}),\\
   \exp_{p_{vv'}}\left(\beta(|z|/r_{vv'}^{1/4})(\exp_{p_{vv'}}^{-1}u_v(z))\right)  & \text{if} \hs z\in D_{f}(2r_{vv'}^{1/4})\bs D_{f}(r_{vv'}^{1/4})\\
     \exp_{p_{vv'}}\left(\beta(r_{vv'}^{3/4}/|z|)(\exp_{p_{vv'}}^{-1}u_{v'}(\rho/z))\right)  &\text{if} \hs z\in D_{f}(r_{vv'}^{3/4})\bs D_{f}(r_{vv'}^{3/4}/2)
   \end{cases}
\end{equation*}
and   $u_{\rho}(z)$ coincides with $u_v(z)$ away from the annuli above.  Here we need to assume  $\rho$ small enough, so that the discs $D_f(4r_{vv'}^{1/4})$ are sent under $u_v$ in a normal neighbourhood of $p_{vv'}$.  Note that the  mapping $\textrm{pgl}^X$ is continuous with respect to $\rho$. \\

\noindent\emph{Estimates from the pregluing.} 
The following estimates are all standard and their proofs can be found in \cite{MS}, \cite{CL}.  The first estimate tells that $u_{\rho}$ is approximately $J$-holomorphic. The second gives a quadratic estimate ensuring existence and unicity of the gluing map. Finally, the third is needed to derive that the constructed gluing map is locally diffeomorphic.
\begin{lem}\label{pregluingestimates} Let $p>2$ an integer, and $U_X$ a proper open subset of $\M_{\S_X}(X)$.
There is a uniform (with respect to $\rho$)  constant $c^X$  such that  for every $(u,\j)\in \M_{\S_X}(X)$:
\begin{equation}\label{almostholomorphicity}\|du_{\rho}\|_{L^{\infty}}\leq c^X,\hspa \|\delbar_{J}u_{\rho}\|_{L^p}\leq c^X |\rho|^{1/2p}.
\end{equation}
Consequently, there is  a uniform constant $c_1^X$   such that
\begin{equation}\label{quadraticestimate}\|N_{u_{\rho}}^X(\xi_1)-N_{u_{\rho}}^X(\xi_2)\|_{L^p}\leq c_1^X(\|\xi_1\|_{W^{1,p}}+\|\xi_2\|_{W^{1,p}})\|\xi_1-\xi_2\|_{W^{1,p}},
\end{equation}
where $N_u^X$ denotes the non linear part in the Taylor expansion  of $\delbar_{J}$ at  $u$:
$$N_u^X(\xi)=\delbar_{J}\exp_u\xi-\delbar_{J}u-D^X_u\xi,\hsp \xi\in W^{1,p}(u^*TX).$$
Finally, let  $u_t:=\{(u_{vt},u_{v't})\}_{v E v'}$, $t\in[0,\upsilon)$ be 
a path  in $\B^{1,p}_X$, with 
$$\zeta\equiv\{(\zeta_{v},\zeta_{v'})\}_{v E v' }:=\left.\frac{d}{dt}\right|_{t=0} u_t.$$
Let $u_{\rho,t}$ denote the corresponding  path of preglued curves and set $\zeta_{\rho}:=\dti u_{\rho,t}$. 
There is a uniform constant $\wt{c}^X$  such that
 \begin{equation}\label{estimatediffeo}\|\zeta_{\rho}\|_{L^{1,p}}\leq \wt{c}^X\|\zeta\|_{L^{1,p}}\hspa\text{and}\hspa \|D^X_{u_{\rho,0}}\zeta_{\rho}\|_{L^p}\leq \|D^X_{u_0}\zeta\|+ \wt{c}^X|\rho|^{1/2p}.
 \end{equation}
% $$\left\|\left.\frac{d}{dt}\right|_{t=0}\delbar_{J_P} u_{\rho,t}\right\|_{L^{p}}\leq \left\|\left.\frac{d}{dt}\right|_{t=0}\delbar_{J_P} u_{t}\right\|_{L^{p}} + C |\rho|^{1/2p}\|\zeta\|_{C^1}.$$
 In particular, if $u_t$ is a path of holomorphic curves,  i.e if $\zeta$ is in the kernel of $D^X_{u_0}$, the first term of the right handside of the second inequality  vanishes.
\end{lem}

\noindent\emph{\underline{Step 2}:  Right inverses.} The gluing operation will give  a holomorphic map with domain $\Sigma_{\rho}$ obtained  by perturbing the preglued map in directions that are  transverse to the kernel of  $D^X_{u_{\rho}}$, i.e lying in the image of a uniformily bounded  family $Q^X_{u_{\rho}}$ of right inverses for $D^X_{u_{\rho}}$.  Below we  sketch the proof of the following (for details we refer to \cite{MS}):

\begin{prop}\label{rightinverse}
Let $p> 2$ and $U_{X}$ as before. There exists a uniformily bounded  family $Q^X_{u_{\rho}}$ of right inverses for $D^X_{u_{\rho}}$, i.e  there is a uniform constant  $c^X$  such that for every $(\j,u)\in U_{X}$:
\begin{equation}\label{boundedrightinverse}\|Q^X_{u_{\rho}}\eta\|_{W^{1,p}}\leq c^X\|\eta\|_{L^p}.
\end{equation}
\end{prop}

\proof First, one constructs an interpolation $\mathbf{w}_{\rho}:=\{u_{v,\rho:\Sigma_v\sra X}\}_{\{v\in V\}}$ between $u$ and $u_{\rho}$ as follows: for $v\in V$,  
\begin{equation} u_{v,\rho}:=\begin{cases} u_{\rho}(z)& \text{if $z\in \Sigma_v\bs\bigcup_{\{v'| vEv'\}}D_{f}(r_{vv'}^{1/4})$},\\
 u_v(y_{v,v'})=p_{vv'} & \text{otherwise}.
\end{cases}
\end{equation}
where $f$ is the flag in $\text{pr}^{-1}(v)$ associated to the edge between $v$ and $v'$. 
As as $\rho$ goes to  0, the ``flattened" map  $u_{v,\rho}$
converges to  $u_v$ in $W^{1,p}$ norm. Thus,  $D^X_{u_{v,\rho}}$ converges to  $D^X_{u_v}$ in the operator norm, and since  $D^X_{u_v}$ is surjective then  $D^X_{u_{v,\rho}}$ also is.
%Of course to see this one needs to identify the fibers of the relevant Fredholm systems, namely $\E_{P,u_i}^p$ et $\E_{P,u_{i,\rho}}$,   by use of the parallel transport induced by the connection  $\nabla$ defined in the first chapter. 
 As a result,   $D^X_{\mathbf{w}_{\rho}}$
%\begin{equation*}D^X_{\mathbf{w}_{\rho}}:\{\{\xi_v\}_{v\in V}\in\prod_{v\in V}\X_{X,u_{v,\rho}}^{1,p}\hs |\hs \xi_v(y_{vv'})=\xi_{v'}(y_{v'v})\}\lra \prod_{i\in I}\E_{X,u_{v,\rho}}^p,
%\end{equation*}
is surjective. Therefore, $D^X_{\mathbf{w}_{\rho}}$ has right inverse $Q^X_{\mathbf{w}_{\rho}}$ which can be chosen uniquely by requiring that its image lies in the  $L^2$-orthogonal complement  of $\ker D^X_{\mathbf{w}_{\rho}}$. Moreover, $Q^X_{\mathbf{w}_{\rho}}$ is  uniformly bounded. Next,  we construct a quasi right inverse $R^X_{u_{\rho}}$ for $D^X_{u_{\rho}}$ out of  $Q^X_{\mathbf{w}_{\rho}}$, defined by,
\begin{equation}\label{quasirightinverse}
R^X_{u_{\rho}}:\E_{X,u_{\rho}}^p\sra \X^{1,p}_{X,u_{\rho}},\hsp \eta\mapsto \Gamma \circ Q^X_{w_{\rho}}\circ \Lambda (\eta).
\end{equation}
 
We give the definitions of $\Lambda$ and $\Gamma$. 
For $v Ev'$, with corresponding edge $e_{f,f'}$,  recall that $\Sigma_{\rho}(f,f')$ is obtained by patching together the annuli $\A_f$ and $\A_{f'}$. Consider the circles in $\Sigma_{\rho}$  corresponding to the circles
 $$\{-R_{vv'}/2\}\times \R/2\pi\Z\cong \{R_{vv'}/2\}\times \R/ 2\pi\Z,$$
and let $\mathcal{C}$ be the union of all these circles.
Consider the  biholomorphism (onto its image), $$\pi_{\rho}:\Sigma_{\rho}\bs \mathcal{C}\sra\Sigma,$$
defined as the identity map outside the annuli $\mathcal{A}_{f,f'}:=\mathcal{A}_{f}\sim \mathcal{A}_{f'}$ and given by: 
$$\pi_{\rho}(z_{f},z_{f'}):=\begin{cases} z_{f} & \text{if $|z_{f}|>|z_{f'}|$}, \\
z_{f'} & \text{si $|z_{f'}|>|z_{f}|$}.
              \end{cases}
$$
 Observe that, 
 $$u_{\rho}:=\begin{cases} \mathbf{w}_{\rho}\circ \pi_{\rho}& \text{on $\Sigma_{\rho}\bs\mathcal{C}$},\\
u_v(y_{vv'})=u_{v'}(y_{v'v})=p_{vv'} & \text{on $\mathcal{C}$.}          \end{cases}
$$

\noindent Now set,
\begin{equation*}\Lambda: L^p(\Lambda_{J}^{0,1}(\Sigma_{\rho},u^*_{\rho}TX))\sra L^p(\Lambda_{J}^{0,1}(\Sigma,w_{\rho}^*TX)),\hspa
\eta \mapsto  \begin{cases}(\pi_{\rho}^*)^{-1}\eta &\text{on $\text{Im}(\pi_{\rho})$}, \\
0 & \text{otherwise}.
           \end{cases}
 \end{equation*}
%\begin{eqnarray*}\Lambda: L^p(\Lambda_{J_P}^{0,1}(\Sigma_{\rho},u^*_{\rho}(TP))&\lra & L^p(\Lambda_{J_P}^{0,1}(\Sigma,w_{\rho}^*(TP))\\
%\eta & \mapsto & \begin{cases}(\pi_{\rho}^*)^{-1}\eta &\text{on$\text{Im}(\pi_{\rho})$}, \\
%0 & \text{otherwise}.
%           \end{cases}
% \end{eqnarray*}

\noindent Next, we define $\Gamma$. Set $\xi_{v,v'}:=\xi(y_{v,v'})$ and $\beta_{r_{vv'}}(z):=\beta(4\log |z|/\log r_{vv'})$, where  
$\beta $ is a cut-off function as before. Then $\Gamma$ is an interpolation between  $\xi_v$ and $\xi_{v'}$:
\begin{equation*}
\Gamma: W^{1,p}(w_{\rho}^*TX)\sra  W^{1,p}(u_{\rho}^*TX),\hspa
\xi:=\{\xi_v\}_{v\in V} \mapsto  \xi_{\rho},%$z\in \Sigma_{\rho}\bs\bigsqcup_{\{i,j|iEj, i<j\}}\mathcal{A}_{i,j}$}, 
\end{equation*}
where $$\xi_{\rho}=\begin{cases} \xi_v(z_{f})+\beta_{r_{vv'}}(z_{f})(\xi_{v'}(z_{f'})-\xi_{v,v'})\hsp\text{when}\hsp |z_{f'}|\leq|z_{f}|\leq |r_{vv'}|^{1/4}\\ 
\text{$\xi_{\rho}=\xi(z)$ otherwise.}
\end{cases}$$
%\begin{eqnarray*}
%\Gamma: W^{1,p}(w_{\rho}^*TP)&\lra & W^{1,p}(u_{\rho}^*TP)\\
%\xi:=\{\xi_i\}_{i\in I} &\mapsto & \xi_{\rho}:=\begin{cases} \xi_i(z_{ij})+\beta_{r_{i,j}}(z_{ij})(\xi_j(z_{ji})-\xi_{i,j})& \text{if $|z_{ji}|\leq|z_{ij}|\leq |r_{i,j}|^{1/4}$},\\
%\xi_j(z_{ji})+\beta_{r_{i,j}}(z_{ji})(\xi_i(z_{ij})-\xi_{i,j})& \text{if $|z_{ij}|\leq|z_{ji}|\leq|r_{i,j}|^{1/4}$},\\
% \xi(z) & \text{otherwise} %$z\in \Sigma_{\rho}\bs\bigsqcup_{\{i,j|iEj, i<j\}}\mathcal{A}_{i,j}$}, 
%\end{cases}
%\end{eqnarray*}
%where $\xi_{i,j}:=\xi(y_{i,j})$ and $$\beta_{r_{i,j}}(z):=\beta(\log |z|/\log r^{1/4}_{i,j}),$$ 
%$\beta $ is as before. 

It follows from the estimates \eqref{briformulation} in Lemma \ref{boundedrightinverseestimates} below that the maps are $D^X_{u_{\rho}}R^X_{u_{\rho}}$ is invertible for small enough gluing parameters. Consequently,
 \begin{equation}\label{rightinversedefinition} Q^X_{u_{\rho}}:= R^X_{u_{\rho}}(D^X_{u_{\rho}}R^X_{u_{\rho}})^{-1}
 \end{equation} 
is the unique right inverse for  $D^X_{u_{\rho}}$ having the same image as $R^X_{u_{\rho}}$. \qed \\

\noindent\emph{Right inverses estimates.} 
The first estimate below  gives the invertibility of $D^X_{u_{\rho}}R^X_{u_{\rho}}$. The second   estimates how the right inverses $Q^X_{u_{\rho}}$ vary  along a  path in $\B^{1,p}_X$ of preglued maps, which is needed to show that the gluing map is locally diffeomorphic. For detailed proofs we refer to \cite{MS}. 

\begin{lem}\label{boundedrightinverseestimates} 
Let  $p>2$. The operator $R^X_{u_{\rho}}$ depends smoothly on $(\j,u,\rho)$
 and there are  constants $C^X$ and $\wt{C}^X$ independant of $\rho$,  such that:
 \begin{equation}\label{briformulation}\|D^X_{u_{\rho}}R^X_{u_{\rho}}\eta-\eta\|_{L^p}\leq \frac{C^X}{|\log |\rho||^{1-1/p}}\|\eta\|_{L^p},\hspace{0.5cm} \|R^X_{u_{\rho}}\eta\|_{W^{1,p}}\leq \frac{\wt{C}^X}{2}\|\eta\|_{L^p}.
 \end{equation}
Moreover, let $\{u_t\}_{t\in[0,\upsilon)}$ be a path of stable maps, with $\zeta:=\dti{u_t}$. Let $\{u_{\rho,t}\}_{t\in[0,\upsilon)}$ be the corresponding path of preglued maps, with $u_{\rho,0}:=u_{\rho}$.
 there is a uniform constant $\ov{c}^X$ such that  
\begin{equation}\label{rightinvder} \|\dti Q^X_{u_{\rho,t}}\|\leq \ov{c}^X\|\zeta\|_{W^{1,p}}.
\end{equation}
\end{lem}
\proof We only give the proof of estimate \eqref{rightinvder}.
\begin{eqnarray*}
\dti Q^X_{u_{\rho,t}}&=&\left(\dti R^X_{u_{\rho,t}}\right)(D^X_{u_{\rho}}R^X_{u_{\rho}})^{-1}+ R^X_{u_{\rho}}\dti (D^X_{u_{\rho,t}}R^X_{u_{\rho,t}})^{-1}.
\end{eqnarray*}
On one hand,
$$\left\|\dti R^X_{u_{\rho,t}}\right\|=\left\|\Gamma \left(\dti Q^X_{w_t}\right) \Lambda\right\|\leq M \|\zeta\|_{W^{1,p}}$$
since, by a standard argument, $\|\dti Q^X_{w_{\rho,t}}\|$ is bounded above by   $C\|\zeta\|_{W^{1,p}}$  for some positive constant $C$ (see \cite{MS}). On the other hand,
$$(D^X_{u_{\rho}}R^X_{u_{\rho}})\dti (D^X_{u_{\rho,t}}R^X_{u_{\rho,t}})^{-1}=-(D^X_{u_{\rho}}R^X_{u_{\rho}})^{-1}\left(\dti D^X_{u_{\rho,t}}R^X_{u_{\rho}}+ D^X_{u_{\rho}}\dti R^X_{u_{\rho,t}} \right)$$
But the norm of  derivative of $D^X_{u_{\rho,t}}$ is bounded and the estimate follows.
\qed\\

\noindent\emph{\underline{Step 3}: the gluing map $Gl_{\S_X}$.} Let $U_X$ be a proper open subset of $\M_{\S_X}(X)$.   Let $B_{\delta_X}(0)$ denote the $\delta_X\in \R^+$ neighbourhood around the 0-section in $\E_X^p(J)$. Also let $B_{\delta_X,u_{\rho}}(0)$ denote the corresponding ball in $\E^p_{X,u_{\rho}}(J)$. Then, for each $(u,\j,\rho)$ we can find  a  unique element
\begin{equation}\label{transverseelement}f^X(u,\j,\rho)\in B_{\delta_X,u_{\rho}}(0)\subset L^p(\Lambda_{J}^{0,1}(S^2,u_{\rho}^*TX))
\end{equation}
 verifying:
\begin{equation}\label{equationgluing}\delbar_{J}(\exp_{u_{\rho}}Q^X_{u_{\rho}}(f^X(u,\j,\rho)))=0.
\end{equation}
The \emph{gluing map} is then defined by 
\begin{eqnarray}\label{gluingmap}Gl_{\S_X}:\mathcal{L}^*_{\S_X,\epsilon_X,U_{X}}&\sra& \M_{0,l}(X,A)\\
(u,\j,\rho)&\mapsto & (\j_{\rho},\exp_{u_{\rho}}Q^X_{u_{\rho}}(f^X(u,\j,\rho))).\nonumber
\end{eqnarray}
for some positive constant $\epsilon_X$ given by the Implicit Function Theorem below.   Note that we directly have the commutativity \eqref{commutativitygluing}. Also, note that in reality, the gluing map is defined on local uniformizing system for $\mathcal{L}_{\S_X,\epsilon_X,U_{X}}$, but since $\S_X$ and $\S$ are stable, the gluing is invariant under $\aut(\j)$, hence well-defined after quotient by the reparametrizations.  Existence and unicity of $f^X(u,\j,\rho)$ are given by the following standard parametric version of the implicit function theorem, see  \cite{MS}.

\begin{implicitfunctiontheorem}
Let $p>2$ and  let $c_1^X$ denote the positive constant in \eqref{quadraticestimate}. There is  a constant  $\epsilon_X$  such that, for every
$(u,\j,\rho)\in \L_{\S_X,\epsilon_X,U_{X}}$  we have a uniform positive constants $\delta_X$, $\epsilon_1^X$, and $c_2^X$,
 verifying:
$$ \|\delbar_{J}u_{\rho}\|_{L^p}\leq \epsilon_1^X, \hspa \|Q^X_{u_{\rho}}\|\leq c^X_2,\hspa\epsilon^X_1< \delta_X/4, \hspa  \epsilon^X_1<(8c_1^X(c_2^X)^2)^{-1}$$
and a smooth map,
$$f^X:\L_{\S_X,\epsilon_X,U_{X}}\sra B_{\delta_X}(0),$$
 such that  $f^X(u,\j,\rho)$ is the unique solution to \eqref{equationgluing}.
Furthermore, $$\|f^X(u,\j,\rho)\|_{L^p}<2\epsilon^X_1.$$
 \end{implicitfunctiontheorem}

\proof\hs The constants $\epsilon_1^X$ and $c_1^X$ are given by the estimates in Lemma \ref{pregluingestimates}, while  $c_2^X$ is given by Proposition \ref{boundedrightinverse}.  Consider the Fredholm fibered map:
\begin{equation*}F^X:(\prgl^X)^*T\B_X\sra \E_X^p(J),\hspa ((u,\j,\rho),\xi)\mapsto  F^X_{(u,\j,\rho)}(\xi):=\Phi^{-1}_{X,u_{\rho}}(\xi)\delbar_{\j_{\rho},J}(\exp_{u_{\rho}}\xi),
\end{equation*}
where $\Phi_{X}$ is the parallel transport induced by the hermitian connection induced  from the L-C connection $\nabla^{TX}$ with respect to $g_J$.  When $\xi=0$, this map coincides with $\delbar_{\j_{\rho},J}u_{\rho}$ so that  
$$DF^X_{(u,\j,\rho)}(0)(\xi)=D^X_{u_{\rho}}\xi.$$ 
For every $x:=(u,\j,\rho)\in \L_{\S_X,\epsilon_X,U_{X}}$, we want to find $f^X(x)$ in $\E_{X,x}^p(J)$, such that $F^X(Q^X_x f^X(x))=0$, in other words  such that 
$$0=F^X_x(0)+f^X(x)+N^X_x(Q^X_xf^X(x)),$$
where $N^X_x$ is the non linear term of the expansion of $F^X_x$ around $\xi\equiv 0$.
Consider the family of operators $$H_x: \E_{X,x}^p(J)\sra \E_{X,x}^p(J), \hspa \eta\to -F^X_x(0)- N^X_x(Q^X_x\eta).$$
We have
$$\|H_x(\eta)\|_{L^p}\leq \|F^X_x(0)\|_{L^p}+\|N^X_x(Q^X_x\eta)\|_{L^p}\leq \epsilon_1^X+c^X_1(c^X_2)^2\|\eta\|^2_{L^p} $$
so that $H$ maps the ball $B_{\delta_X}(0)$ into itself whenever $\epsilon_1^X +c^X_1(c^X_2)^2 \delta_X^2\leq\delta_X$. This is realized when $\epsilon_1^X <\delta_X/4$ and $4c^X_1(c^X_2)^2\epsilon_1^X<1/2$. Furthermore, 
 \begin{eqnarray*}\|H_x(\eta_1)-H_x(\eta_2)\|_{L^p}&=&\|N^X_x(Q^X_x\eta_1)-N^X_x(Q^X_x\eta_2)\|_{L^p} \\
 &\leq & 2c^X_1(c_2^X)^2\delta_X \|\eta_1-\eta_2\|_{L^p}\\
 &<& \|\eta_1-\eta_2\|_{L^p}.
 \end{eqnarray*}
Thus,  $H$ is  a contraction map so that existence and uniqueness follow. One can further show  that the map $H$ defines a contraction from  $B_{2\epsilon_1^X}(0)$ to itself, when $4c^X_1(c_2^X)^2\epsilon_1^X<1/2$. This gives the estimate $\|f^X(x)\|_{L^p}<2\epsilon^X_1$.  Smoothness of $f^X$ follows from implicit function theorem since  $DF^X_x(Q^X_xf^P(x))$ is an isomorphism from  $Q^X(\E_{X}^p(J))$ to $\E_{X}^p(J)$. 
 \qed\\

Next, we prove that $Gl_{\S_X}$ is locally diffeomorphic:  we show that  there is a proper open subset $U'_{X}\subset U_{X}$   and a positive constant  $\delta_X'<\delta_X$ with respect to  which the map $\left.Gl_X\right|_{U'_{X}}$ is a diffeomorphism onto its image.
  
%\begin{theorem} \label{commutgluing} Let $\S_X$, $U_{X}$  and $\delta_X$, as in the above proposition. Then there is a proper open subset $U'_{X}\subset U_{X}$  such that $\pi(U'_{\S_P})=U'_{\S_B}$, and a positive constant  $\delta_X'<\delta_X$ with respect to  which the map $\left.Gl_X\right|_{U'_{X}}$ is a diffeomorphism onto its image.  Furthermore let $\epsilon_P$ and $\epsilon_B$ as above then
%$$\pi\circ Gl_{\S_P}=Gl_{\S_B}\circ \pi.$$
%\end{theorem}
   Fix $(u_0,\j,\rho_0)$ in $\mathcal{L}^*_{\S_X,\epsilon_X,U_{X}}$ and set $u_{\rho_0}$ to be the corresponding pregluing. Let $W$ be an open neighbourhood of 0 in $W^{1,p}(S^2,u_{\rho_0}^*TX)$. Furthermore, let $U_1$ be the image under $\prgl^X$ of a neighbourhood $U_0$ of $(u_0,\j,\rho_0)$. We decompose $Gl_{\S_X}$ as follows:
\begin{eqnarray*}
U_0\stackrel{\prgl^X}{\lra} U_1\stackrel{1\times f^X}{\lra} U_1\times B_{\delta_X}(0) \stackrel{\Phi}{\lra} W \stackrel{\exp_{u_{\rho_0}}}{\lra} \M_{0,l}(P) 
\end{eqnarray*} 
where $\Phi(u_{\rho},\eta):=\xi_{{\rho}}+Q^X_{u_{\rho}}\eta$, and $\exp_{u_{\rho_{0}}}\xi_{{\rho}}=u_{\rho}$.  Note that the pregluing map is  locally  diffeomorphic. Also, for any path $u_{\rho_t}$ starting at $u_{\rho_0}$, with derivative $\xi$ at $t=0$, differentiating  the equation 
$$0=\delbar_{J}u_{\rho_t}+f^X(u_{\rho_t})+N^X_{u_{\rho_t}}(Q^X_{u_{\rho_t}}f^X(u_{\rho_t})),$$
and using the estimates proved  so far, one obtains the estimate:
$$\|\frac{df^X}{d\xi}\|_{L^p}< Cr^{1/2p}\|\xi\|_{W^{1,p}}.$$
This ensures that for small enough gluing parameter, the differential of the gluing map is well-defined.  Therefore it suffices to check that  $\Phi$ is a diffeomorphism. To prove  this  we identify $W^{1,p}(S^2,u_{\rho_0}^*TX)$ with $\ker D^X_{u_{\rho_0}}\oplus L^p(S^2,u_{\rho_0}^*TX)$ via the map 
$$\xi\mapsto ((Id-Q^X_{u_{\rho_0}}D^X_{u_{\rho_0}})\xi, D^X_{u_{\rho_0}}\xi),$$ and then rewrite $\Phi$ as 
$$\Phi(u_{\rho},\eta)=((Id-Q^X_{u_{\rho_0}}D^X_{u_{\rho_0}})(\xi_{u_{\rho}}+Q^X_{u_{\rho}}\eta),  D^X_{u_{\rho_0}}(\xi_{\rho}+Q^X_{u_{\rho}}\eta)).$$
Then, $D^X\Phi_{(u_{\rho},\eta)}=Id+K_{(u_{\rho},\eta)}$ where
$$K_{(u_{\rho},\eta)}(\xi,\zeta)=
\begin{pmatrix}
   (Id-Q^X_{u_{\rho_0}}D^X_{u_{\rho_0}})(\xi+\frac{dQ^X_{u_{\rho}}}{d\xi}\eta)-\xi   &   (Id-Q^X_{u_{\rho_0}}D^X_{u_{\rho_0}})Q^X_{u_{\rho}}\zeta  \\
D^X_{u_{\rho_0}}(\xi+\frac{dQ^X_{u_{\rho}}}{d\xi}\eta)      &(Q^X_{u_{\rho_0}}D^X_{u_{\rho_0}} - Id)\zeta 
\end{pmatrix}$$ 
Hence we deduce from corollary \ref{rightinverse} and lemma \ref{rightinvder} that for proper susbet $U'\subset U$ and $\delta_1< \delta$, the operator $D^X\Phi_{(u_{\rho},\eta)}$ is invertible for any $(u_{\rho},\eta)\in U'\times B_{\delta_1}(0)$, and that
$$\|D^X\Phi_{(u_{\rho},\eta)}\|\leq 2.$$
For the injectivity of $\Phi$,  let $N(\xi_{\rho},\eta)$ be the non linear part in the expansion of $\Phi(u_{\rho},\eta)$ around $(u_{\rho_0},0)$.  Then,  
$$N^X(\xi_{{\rho}},\eta)=((Id-Q^X_{u_{\rho_0}}D^X_{u_{\rho_0}})Q^X_{u_{\rho}}\eta, D^X_{u_{\rho_0}}Q^X_{u_{\rho}}\eta-\eta).$$
By standard arguments, for $(\xi_{\rho_1},\eta_1)$ and $(\xi_{\rho_2},\eta_2)$, the following estimate holds,
$$\|N^X(\xi_{{\rho_1}},\eta_2)-N^X(\xi_{{\rho_2}},\eta_2)\|_{p}\leq C (\|(\xi_{\rho_1},\eta_1)\|_{1,p}+\|(\xi_{\rho_2},\eta_2)\|_{1,p})(\|(\xi_{\rho_1}-\xi_{\rho_2},\eta_1-\eta_2)\|_{1,p})$$
Assume that $\Phi(u_{\rho_1},\eta_1)=\Phi(u_{\rho_2},\eta_2)$, for $(u_{\rho_1},\eta_1)\neq (u_{\rho_2},\eta_2)$. Then,
\begin{eqnarray*}\|D^X\Phi_{(u_{\rho_0},0)}(\xi_{\rho_1}-\xi_{\rho_2},\eta_1-\eta_2)\|_{L^p}&=&
\|N^X(\xi_{{\rho_1}},\eta_2)-N^X(\xi_{{\rho_2}},\eta_2)\|_{L^p}.
\end{eqnarray*}
Since   $\|(\xi_{\rho_1}-\xi_{\rho_2},\eta_1-\eta_2)\|_{W^{1,p}}>0$ by assumption,  we have,  
$$0<\|D^X\Phi_{(u_{\rho_0}},0)\|_{L^p}\leq C(  \|(\xi_{\rho_1},\eta_1)\|_{1,p}+\|(\xi_{\rho_2},\eta_2)\|_{1,p}),$$
which is impossible for small enough $(\xi_{\rho_i},\eta_i)$, $i=1,2$.  
\qed\\

As such,  the construction above still applies when the domain is not stable, nevertheless:
\begin{enumerate}[$\bullet$]
\item the obtained gluing maps are only defined at the level of parametrized maps, i.e before quotienting by  $\aut(\j)$ even though these maps are  $\aut(u,\j)$ equivariant. 
\item so far we have parametrized our gluing maps according to the gluing for the domains, but the gluing for nodal surfaces is neither injective, nor   locally diffeomorphic in the unstable case. Therefore, we can't treat the gluing parameters $\rho$ as parameters.
\end{enumerate}
The problem is to find  a slice for the  action of the group of reparametrizations of the domain.  

%%%%%%%%%%%%%%%%% Gluing unstable components %%%%%%%%%%%%%%%%%%%%%%%%%%%

%%%%%%%%%%%%%%%%%%%%%%%%%%%%%%%%%%%%%%%%%%%%%%%%%%%%%%%%%%%

\subsubsection{Gluing unstable components}\label{sectiongluingunstable}
Let $(u,\j)\in \M_{\S_X}(X)$, and let $v\in V$ be a $\pi_{pt_*}$-unstable component. Note that $|v|>2$ is possible. Following Chen-Li \cite{CL},
we describe the notion of \emph{balanced component} when $|v|\leq2$. For the commodity of the readers, we furnish the details below, as it will be generalized to the case of general base $B$. Let $\textbf{\text{tr}}\cong\C\subset G$ denote the  subgroup of translations  and let  $\m\cong\C^*\subset G$ denote the subgroup   acting on $S^2=\C\cup\{\infty\}$ by complex multiplication. Then  the semi direct product $\mathcal{G}:=\t \ltimes \m$,  acts on $\C$ by:  $$(m,t)\in \mathcal{G},\hspa (t,m)\cdot z:=m(z-t).$$

\noindent 1)\underline{\emph{Balanced maps.}} First assume $|v|=2$. Upto the action of $\m$, the component $v$ can be parametrized by $\CP^1$ with special points $0=[1:0]$ and $\infty=[0:1]$. Identify $\Sigma_v\bs\{\infty\}$ with $\C$. The parametrization is \emph{balanced} if:
\begin{equation}\label{balancedconditionenergy}
\frac{1}{2}\int_{|z|\leq 1}\|du\|^2dvolS^2= \hbar
\end{equation} 
where $\hbar$ denotes the minimal energy of a non-constant pseudo-holomorphic map in
$X$.  Next, when $|v|=1$, upto the action of $\mathcal{G}$, the component    $v$ can be parametrized by $\CP^1$ such that the special point is $\infty$. Identify $\Sigma_v\bs\{\infty\}$ with $\C$. The parametrization is \emph{balanced} if \eqref{balancedconditionenergy} holds and if  the center of energy  of $u$ is 0:
\begin{equation}\label{balancedconditioncenterofenergy} \int z \|du\|^2dvolS^2=0,\hspa\text{where}\hspa z\in \C.
\end{equation}
We will say that $u$ is \emph{centered} if it is the case.
The $\pi_{pt_*}$-unstable component $\Sigma_v$ with balanced parametrization is called a \emph{balanced component}. Recall that the reparametrizations of a $\pi_{pt_*}$-stable component is of finite order. Here the reparamatrizations of a balanced component is given by $S^1$. Note that the neighbourhoods of the special points in a balanced parametrization for $\Sigma_v$ can be put  in 
standard cylindrical coordinates: for $\infty$,
$$[0,\infty)\times \R/2\pi\Z \cong D^*(r_{\infty})\,\,:\,\, (s,t)\mapsto e^{s+it}$$
while for 0,
$$(-\infty,0]\times \R/2\pi\Z \cong D^*(r_{0})\,\,:\,\, (s,t)\mapsto e^{s+it}.$$

To each unbalanced map $u_v:\Sigma_v\sra X$ smoothly corresponds a unique element $\phi_X^b(u_v)\in \mathcal{G}$ consisting of the pair of translation and real dilation such that:  the center of energy of 
\begin{equation}\label{associatebalancedmap}u^b_v:=u_v\circ(\phi_X^b(u_v))
\end{equation}
 is zero;  half the total energy of $u^b_v$ lies in the unit disc around zero. 
Let  $\wt{\M}^{b}_{0,i}(X,A_v)$, $i=1,2$, denote the sets of   balanced $J$-holomorphic maps representing $A_v$ with one and two marked points. The map $u_v\mapsto u^b_v$ sends an orbit of  $\mathcal{G}$ to an $S^1$  orbit where $S^1$ acts by rotations around the origin, hence we have the following natural identifications:
$$\M_{0,i}(X,A)\cong\left.\wt{\M}^{b}_{0,i}(X,A)\right/S^1,\hspa i=1,2.$$

More generally, we say that  $(u,\j)\in\wt{\M}_{\S_X}(X)$ is \emph{balanced} if each of $\pi_{pt_*}$-unstable component $v$ with $|v|\leq2$ is  balanced. Let $\wt{\M}^b_{\S_X}(X)$ denote the subset of all balanced stable maps, and call it  \emph{moduli space of balanced $J$-holomorphic maps for  $\S_X$}. On $\wt{\M}^b_{\S_X}(X)$ the action of the reparametrizations reduces to the action of  $\aut_{red}^X$:
\begin{equation}\label{reducedgroup}\aut_{red}^X\cong (S^1)^{|V_b|}\times \aut(\S_X).
\end{equation}
where $V_b\subset V$ denotes the subset of $\pi_{pt_*}$-unstable components $v$ with $|v|\leq 2$. From the discussion above we have the identification:
\begin{equation}\label{definitionmodulispacebalancedmaps}\M_{\S_X}(X)=\left.\wt{\M}^b_{\S_X}(X)\right/\aut_{red}^X.
\end{equation}
Next we define the gluing maps for balanced curves.\\

\noindent 2) \underline{\emph{Gluing balanced maps.}} Let $\S:=S_{\pi_{pt_*}}(\S_X)$. Set  $\S^u:=\F_X(\S_X)$ and let $\wt{\M}_{\S^u}$ denote the set of parametrized nodal curves having $\S^u$ as stratum data. Let $\wt{\L}_{\S^u}$ be the corresponding fiber bundle of  gluing parameters. Also, let 
$\wt{\L}_{\S_X}$ denote the  bundle $\F_X^*\wt{\L}_{\S^u}$ over $\wt{\M}_{\S_X}(X)$, and let $\wt{\L}^b_{\S_X}$ denote its restriction to $\wt{\M}^b_{\S_X}(X)$. The forgetful map $\ov{\pi}_{pt}^{\S_P}$ induces a map:
\begin{equation}\label{stabilisationgluingparameters}\ov{\ov{\pi}_{pt}}^{\S_X}: \wt{\L}^b_{\S_X}\sra\wt{\L}_{\S},\hspa (u,\j,\rho)\mapsto (\ov{\pi}^{\S_X}_{pt}(u,\j), \rho^{\text{st}}),
\end{equation}
where  $\rho^{\text{st}}$ leaves unchanged the  gluing parameters between $\pi_{{pt}_*}$-stable components,  forgets about all the gluing parameters of  components lying in contracted branch and sends the gluing parameters of a connecting branch to the product of the parameters of  the corresponding connecting chain. 
%The forgetting-the-map map  induces a bundle map:
%$$\ov{\F}_X:\wt{\L}_{\S_X}\sra\wt{\L}_{\S^u},\hspa (u,\j,\rho)\mapsto (\j, \rho).$$

Note that the group  $\aut_{red}^X$ acts naturally  by rotations on $\wt{\L}^b_{\S_X}$.   Since we consider balanced parametrizations,   $\rho^{st}$ is invariant under the reparametrizations of $\j$ and $\ov{\ov{\pi}_{pt}}^{\S_X}$ is well-defined after  quotient:
\begin{equation*}\ov{\ov{\pi}_{pt}}^{\S_X}: \L^b_{\S_X}:= \left.\wt{\L}_{\S_X}\right/\aut^X_{red}\sra \L_{\S}. 
\end{equation*}

Now, by hypothesis \eqref{strongreg}, for any element of $\wt{\M}^b_{\S_X}(X)$, the linearization of the  Cauchy-Riemann operator is surjective.  It follows from Theorem \ref{gluingtheorem}  that for any proper open subset $U_X$ of $\wt{\M}^b_{\S_X}(X)$, that we may choose to be  $\aut_{red}^X$-invariant, there exists $\epsilon_X>0$ and a map 
$$\wt{Gl}^b_{\S_X}:\wt{\L}^{b*}_{\S_X,\epsilon_X,U_X}\sra \wt{\M}_{0,l}(X,A).$$
This map is $\aut_{red}^X$-equivariant, thus the gluing map is well-defined after the quotient by the action:
$$Gl^b_{\S_X}:\left.\wt{\L}^{b*}_{\S_X,\epsilon_X,U_X}\right/\aut_{red}^X\sra \wt{\M}_{0,l}(X,A).$$

Note that the domain of $Gl^b_{\S_X}(u,\j,\rho)$ is  $gl_{\S^u}(\j,\rho)$. Since $\S^u$ is not stable, $gl_{\S^u}$ is not injective nor it is locally diffeomorphic and $\rho$ cannot be treated as a parameter anymore.  
Moreover, note that for $l\geq 3$, this gluing map takes value in $\M_{0,l}(X,A)$, while for  $l<3$, one needs to make sure that the image of the gluing gives a slice for the action of the  automorphisms of $S^2$ with less than three marked points. %There are two issues arising:   and the compatibility with respect to $\ov{\pi}_{pt}^{\S_X}$. The second one
%First note that a contracted branch is connected to a unique $\pi_*$-stable component. The subgraph is a tree with a distinguished  \emph{root}  that is attached to the $\pi_*$-stable component, and each vertex has two edges or less. We make the following convention: when we parametrize a component $\Sigma_v$ of the branch by $\CP^1$, the special point that is the closest to the  root will be parametrized as $\infty=[0:1]$.

\begin{theorem}(Chen-Li \cite{CL})\label{TheoremChenLi} $Gl^b_{\S_X}$ is locally  diffeomorphic. Furthermore,
\begin{equation}\label{gluingcompatibilitystandardcase} \ov{\pi}_{pt}^{\S_X}\circ Gl^b_{\S_X}=gl_{\S}\circ \ov{\ov{\pi}_{pt}}^{\S_X}.
\end{equation}
\end{theorem}

\proof  The gluing among $\pi_*$-unstable components is divided into two cases: (1) the gluing between a balanced component and a stable component; (2) the gluing between two balanced components.   Let $\Sigma_v$ and $\Sigma_{v'}$ be the two components to be glued at the edge $e_{f,f'}$. We may write
 $$\Sigma_v=(S^2,z_f\equiv0,\{x_k\}_{k=1,...,m})\hspa\text{and}\hspa\Sigma_{v'}=(S^2,z_{f'}\equiv\infty, \{x'_k\}_{k=1,...,m'}).$$  
Then (1) and (2) can be deduced from: (a)
$\Sigma_v$ is $\pi_*$-stable while $\Sigma_{v'}$ is not,  $m\geq3$ and $m'=0$;  (b) both components are $\pi_*$-unstable and $m=m'=0$.

%Again this defines a new domain $\Sigma_{\rho}(f,f')$. Doing this at every edge gives a new domain $\Sigma_{\rho}$ for $\rho\in \C_{\j}$, also  denoted by $gl_{\rho}(\j)$.

\underline{\emph{Case (a)}}.   For simplicity we forget about the marked points on $\Sigma_v$. Let $V_{z_{f}}$ denotes a neighbourhood of $z_{f}$, then   $gl_{\S^u}$ sends  the neighbourhood $V_{z_f}\times \C^*$ of $\j$  to $(S^2,\infty)$.    Let $\wt{N}_{u_0}\times V_{\rho_0}$ denote a neighbourhood  of  $(u_0,\rho_0)$ in $\wt{\M}_{\S_X}(X)\times \C^*_{\j}$, 
in local coordinates for $\wt{\L}_{\S_X}$. We want to define a gluing map:
$$Gl:\frac{\wt{N}_{u_0}\times V_{\rho_0}}{\mathcal{G}}\sra \M_{0,l}(X,A).$$
By choosing a proper slice for the action of $\mathcal{G}$ we can construct a well-defined gluing map, namely  $Gl^b_{\S_X}$,  locally given by:
$$Gl^b_{\S_X}: N_{u_0,\j}\times V_{z_{f}}\times V_{\rho_0}\sra \M_{0,l}(X,A),$$
where $N_{u_0,\j}$ stands for an $S^1$ slice in 
$\wt{\M}^b_{\S_X}(X)\cap\F_{X}^{-1}(\j)$ around $u_0$. 
To see that this map is  locally diffeomorphic, we compare it to a gluing map already encountered. 
To do this, use the identification of $\m$ with $V_{\rho_0}$ to obtain a new map: 
$$Gl_1:\frac{\wt{N}_{u_0}}{\t} \times\{\rho_0\}\sra \M_{0,l}(X,A).$$
Next, add  two marked points on the second component, $\{0\}$ and $\{1\}$,  in order to stabilize, and let $\S_X(2)$ denote the corresponding stratum data.  Then  $\wt{N}_{u_0}$ can be written as a product  of $V_{z_{f}}$ with a neighbourhood  $N_{u_0,z_{f}}$ of $u_0$ in 
$\F^{-1}_X(\j)\cap\wt{\M}_{\S_X(2)}(X)$. Then we use the natural identification between $\t$ and $V_{z_{f}}$ in order to obtain a gluing map $Gl'_1$  defined on $N_{u_0,z_f}\times\{z_f\} \times\{\rho_0\}$, which we know is diffeomorphic. This map is in fact  $C^1$-close to $Gl^b_{\S_X}$ so that  $Gl^b_{\S_X}$ is locally  diffeomorphic.

%The idea is that  each  balanced map  can be obtained (upto an $S^1$ action)  by adding two marked points, $\{0\}$ and $\{1\}$, on the root. But this choice fixes, upto rotations,  $\rho$, and $z_f$ via the gluing of nodal surface: 

\underline{\emph{Case (b)}}. We may assume that $z_f=0$ and $z_{f'}=\infty$ in the balanced parametrizations.   This time  $\Sigma_v=(S^2,0)$ and $\Sigma_{v'}=(S^2,\infty)$.  Furthermore,
$$\aut(\Sigma)=\cal{G}_1\times\cal{G}_2=(\t_1\ltimes\m_1)\times(\t_2\ltimes\m_2).$$
%For a gluing parameter $\rho_{vv'}\in \C^*_{\epsilon_X}$ with $\rho_{vv'}=e^{-R_{vv'}+2\pi i\theta_{vv'}}$, $\theta_{vv'} \in [0,1)$, and $R_{vv'}>0$, we glue the two components by identifying the following annuli:
%$$[-R_{vv'},0]\times \R/2\pi\Z \cong [0,R_{vv'}]\times \R/2\pi\Z\,\,:\,\, (s,t)\mapsto (s+R_{vv'},t+\theta_{vv'}).$$
Again $\aut(\Sigma)$  acts on  $\C^*_{\epsilon_X}$, and we set $\aut_{vv'}(\Sigma)$ to be the normal subgroup that fixes the gluing parameters under this action. Hence,   $\aut_{vv'}(\Sigma)$ is isomorphic to $\t_1\times\t_2\times \C_1^*$,  where  $\C^*_1:=\Delta^{-1}(1)$ and 
$$\Delta:\m_1\times\m_2\sra \m,\hspa (m_1,m_2)\mapsto m_1m_2.$$
 The complementary of  $\C^*_1$ in  $\m_1\times\m_2$ is  denoted by $\C^*_2$ and is naturally identified with $\C^*_{\epsilon_X}$.
 The map  $Gl^b_{\S_X}$ comes from a map on  $\wt{\M}_{\S_X}(X)$ locally given by:
$$\wt{Gl}_{\S_X}:U_{u_0}\times V_{\rho_0}\sra \wt{\M}_{0,0}(X,A),$$
where  $U_{u_0}$ denotes  a  $\cal{G}_1\times\cal{G}_2$-invariant neighbourhood of $u_0\in\wt{\M}_{\S_X}(X)$, and $V_{\rho_0}$ is a neighbourhood of $\rho_0\in\C^*_{\epsilon_X}$. The map  $Gl^b_{\S_X}$ is  obtained by choosing an appropriate slice for the action of  $\aut(\Sigma)$ once we restrict ourselves to balanced maps. Quotienting by the automorphism group we get a map 
$$Gl_{\S_X}:\frac{U_{u_0}\times V_{\rho_0}}{\cal{G}_1\times \cal{G}_2}\sra \wt{\M}_{0,0}(X,A),$$
which we would like to take values in  $\wt{\M}_{0,0}(X,A)/\aut(S^2)$. Using the identification between $V_{\rho_0}$ and  a neighbourhood of the identity in  $\C^*_2$ we get a new map:
$$Gl_{\S_X,1}:\frac{U_{u_0}}{\aut_{vv'}(\Sigma)}\times \{\rho_0\}\sra\frac{\wt{\M}_{0,0}(X,A)}{\aut(S^2)},$$
which  is well-defined, locally diffeomorphic, and close to $Gl^b_{\S_X}$. That this map is indeed well-defined follows since $\aut_{vv'}(\Sigma)$ and $\aut(S^2)$ are locally diffeomorphic around the identity.

\underline{Proof of \eqref{gluingcompatibilitystandardcase}.}  
%Again, for $(u,\j)\in\wt{\M}^b_{\S_X}(X)$, the correpson denote the image of $\wt{\M}^b_{\S_P}(P)$ under $\pi$: $(u,\j)\mapsto (\pi(u),\j)$. Note that the parametrization of the domain $\pi(u)$ comes from the parametrization of the domain of  $u$ up to an  $S^1$-action.
 First consider a contracted branch. It is connected to a unique $\pi_*$-stable component. The subgraph is a tree with a distinguished  \emph{root}  that is attached to the $\pi_*$-stable component. 
 We can parametrize each  component $\Sigma_v$ of the branch by $\CP^1$, such that  the special point  closest to the  root is  given by  $\infty=[0:1]$. By gluing from the farthest component to the closest,  it suffices to consider  the case of only one component attached to a root. But in that case, the glued surface  $\j_{\rho}$ is isomorphic to  the domain of the root component for every small enough $\rho$, so that  the forgetful map takes the corresponding glued maps to the same point. 
% Since  $u$ is balanced, $\aut(\j)$ is identified to $S^1$. Now, each unbalanced map the balanced    maps  can be obtained (upto an $S^1$ action)  by adding two marked points $\{0\}$ and $\{1\}$ on the root. But this corresponds to fixing, upto rotations,  $\rho$, and the nodal point on the  $\pi_*$-stable component.  

 Now consider a connecting branch. Observe that we can treat the components that are not in the connecting chain in the same way as the the components of a  contracted branch. Therefore we only consider the case where the connecting branch coincides with the connecting chain.  This chain  is connected to  exactly two $\pi_*$-stable components. The subgraph is a tree with two distinguished components, a \emph{root} and a \emph{top}  that are   attached to $\pi_*$-stable components. 
 We can parametrize each  component $\Sigma_v$ of the branch by $\CP^1$, such that  the special point that is the closest to the root is  given by  $\infty=[0:1]$, and the point that is farthest  is given by $0=[1:0]$. 
Let $k$ denote the number of components of the connecting chain. By adding one marked point, say $\{1\}$, on every component of the chain, the resulting nodal Riemmann surface $\j'$  becomes stable. Then for every $\rho$, the gluing $\Sigma_{\rho}$ is obtained from $\Sigma'_{\rho}$ by forgetting the added marked points.  Note that each such choice of point $\{1\}$ fixes one gluing parameter on each component of the chain. For $\rho=(\rho_1,...,\rho_{k+1})$, starting from the top, we can fix all gluing parameter a fixed  small $\rho_0$, except for the gluing parameter associated to the root and the $\pi_*$-stable component, which is then given $\wt{\rho}$ the product of all the $\rho_i$'s.
Now, $\Sigma'_{\rho}$ and $\Sigma^{st}_{\wt{\rho}}$
are isomorphic since they are both realized by gluing on a cylinder of length $\sum\log|\rho_i|$.
 \qed\\

%\subsubsection{Gluing unstable components II: the constracted branches} 

%\begin{theorem} 
%Let $\S_P$ be a stable stratum data and $\S_B$ its   projection under $S_{\pi_*}$. 
%Let $U_{P}$ be a proper open subset of $\M_{\S_P}(P)$ projecting over a proper open subset $U_{B}$ of $\M_{\S_B}(B)$. 
%Let  $Gl_{\S_P}$ and $Gl_{\S_B}$ be the gluing maps above these open sets.  Then
%$$\ov{\pi}\circ Gl_{\S_P}=Gl_{\S_B}\circ \ov{\pi}_{\S_P}.$$
%\end{theorem}

%%%%%%%%%%%%%%%%%%%%%%%% PARTIE A UTILISER DANS CETTE SECTION 

%\subsubsection{Gluing unstable components III: the connecting branches}

%%%%%%%%%%%%%%%%%%%%%%%%%%%%%%%%%%%%%%%%%%%%%%%%%%%

                                       % GLUING NON STABLE 
                                       
    %%%%%%%%%%%%%%%%%%%%%%%%%%%%%%%%%%%%%%%%%%%%%%%%%                                   

%%%%%%%%%%%%%%%%%FIN DE CETTE PARTIE A UTILISER

%%%%%%%%%%%%%%%%% FIN DE LA PARTIE REVISEE %%%%%%%%%%%%%%%%%%%%%%%%%%%

%%%%%%%%%%%%%%%%%%%%%%%%%%%%%%%%%%%%%%%%%%%%%%%%%%%%%%%%%%%

%%%%%%%%%%%%%%%%% Partie avec B general %%%%%%%%%%%%%%%%%%%%%%%%%%%

%%%%%%%%%%%%%%%%%%%%%%%%%%%%%%%%%%%%%%%%%%%%%%%%%%%%%%%%%%%

\subsection{Gluing for general $B$} It is completely parallel to the special case $B=pt$ treated above. We mainly  point out the differences.
Let $\S_P:=(V_P,\text{Fl}_P;\text{pr}_P,\varrho_P)$ be a stable stratum data for $\ov{\M}_{0,l}(P,\sigma)$, and let $\S_B:=(V_B,\text{Fl}_B;\text{pr}_B,\varrho_B)$ be its image under $S_{\pi_*}$. Also let $\S$ denote the image of $\S_P$ (or $\S_B$) under the forgetful map.

\subsubsection{Pregluing} 
Let $(u,\j)$ be a  $J_P$-holomorphic stable map in  $P$, representing the stratum data $\S_P$. We show  that the pregluing of $(u,\j)$  projects under $\pi$ to the pregluing of $(\pi(u),\j)$ with same gluing parameter:

\begin{lem}\label{projpregl} For every $(u,\j)$ stable map, and gluing parameter $\rho$:
\begin{equation*}  \pi(u_{\rho})\equiv\pi(\prgl^P(u,\j))= \prgl^B(\pi(u),\j)\equiv \pi(u)_{\rho}
\end{equation*}
\end{lem}

\proof\hs Assume for simplicity in the notations that $|V|=2$ with elements $v$ and $v'$ and let $e_{f,f'}$ be the corresponding edge. Set 
$$\xi_v(z):=\exp_{p_{vv'}}^{-1}u_v(z),\hspa \xi_{v'}(z):=\exp_{p_{vv'}}^{-1}u_{v'}(\rho_{vv'}/z),$$ 
and let $\beta^+$ and $\beta^-$ respectively denote the functions $\beta(|z|/r_{vv'}^{1/4})$ and $\beta(r_{vv'}^{3/4}/|z|)$.
From  \eqref{identiteexpo}  we deduce that on $\Sigma_v$,
\begin{equation*}
\pi(u_{\rho})=\begin{cases}
   \pi(u_v) & \text{if $z\in \Sigma_v\bs D_{f}(2r_{vv'}^{1/4})$},\\
    \pi(p_{vv'})=\pi(u_v(y_{v,v'})) =\pi(u_{v'}(y_{v,v'}))&\text{if $z\in D_{f}(r_{vv'}^{1/4})\bs D_{f}(r_{vv'}^{3/4})$}\\
    \exp_{\pi(p_{vv'})}\left(\beta^+\pi_{*_{p_{vv'}}}\xi_v(z)+\beta^-\pi_{*_{p_{vv'}}}\xi_{v'}(z)\right) & \text{otherwise}.
\end{cases}
\end{equation*}
We can rewrite this last expression as follows:
$$\exp_{\pi(p_{vv'})}\left(\beta^+\exp^{-1}_{\pi(p_{vv'})}(\pi(u_v(z)))+\beta^-\exp^{-1}_{\pi(p_{vv'})}(\pi(u_{v'}(\rho_{vv'}/z)))\right),$$
so that  $\pi(u_{\rho})$ coincides  with  the pregluing  $\pi(u)_{\rho}$.
\qed\\

\begin{rem} The map obtained is not the pregluing of the stabilized map $\ov{\pi}(u,\j)$. If $u$ is only made of $\pi_*$-stable  components, these two  pregluings coincide. 
%Also note that if every component of  $u$ lie in the same fiber  $F_{\pi(u)}$ of $P$,  then the pregluing of  $u$ must stay in $F_{\pi(u)}$, for $\xi_i,\xi_j$ are then elements of $T_{p_{i,j}}P^v$. 
Furthermore, $u_{\rho}$ may not  necessarily  lie in the restriction $\left.P\right|_{\pi(u)}$, e.g if  $u$ has only one $\pi_*$-stable component $u_0$. Nevertheless, we will see that the \emph{glued map} projects to  $\pi(u_0)$.
\end{rem}

Let $p>2$ and let $(\j,u)\in \wt{\M}_{\S_P}(P)$. Here are some estimates that follow directly from Lemma \ref{pregluingestimates} in the $B=pt$ case.
% and $U_{B}\subset\wt{\M}_{\S_B}(B)$ be proper open subsets such that $\pi(U_{P})=U_{B}$. 
From \eqref{almostholomorphicity}  there are   uniform positive constants $c^B$ and $c^v$ such that: $\|du^v_{\rho}\|_{L^{\infty}}\leq c^v$, $\|d\pi(u)_{\rho}\|_{L^{\infty}}\leq c^B$, and:
$$ \|(\delbar_{J_P}u_{\rho})^v\|_{L^p}\leq c^v |\rho|^{1/2p}, \hspa \|\delbar_{J_B}\pi(u)_{\rho}\|_{L^p}\leq c^B |\rho|^{1/2p}.$$
Moreover, by definition of  $g_{J_P}$:
$$\|du_{\rho}\|_{L^{p}}\leq \|du_{B,\rho}\|_{L^{p}}+ \|du^v_{\rho}\|_{L^{p}},$$
hence there is a positive uniform constant  $c^P$ such that:
\begin{equation}\label{almostholomorphicityfibered}\|du_{\rho}\|_{L^{\infty}}\leq c^P,\hspa \|\delbar_{J_P}u_{\rho}\|_{L^p}\leq c^P |\rho|^{1/2p}.
\end{equation}
Also, from \eqref{quadraticestimate} there is a uniform  positive constant $c_1^P$ such that:
\begin{equation} \label{quadraticestimatefibered}\|N_{u_{\rho}}^P(\xi_1)-N_{u_{\rho}}^P(\xi_2)\|_{L^p}\leq c_1^P(\|\xi_1\|_{W^{1,p}}+\|\xi_2\|_{W^{1,p}})\|\xi_1-\xi_2\|_{W^{1,p}}.
\end{equation}
Furthermore, from \eqref{estimatediffeo},  if   $u_t:=\{(u_{vt},u_{v't})\}_{v E v'}$, $t\in[0,\upsilon)$ is a 
a path  in $\B^{1,p}_P$, with $\zeta:=\dti{ u_t}$, 
and if $u_{\rho,t}$ is the corresponding  path of preglued  with $\zeta_{\rho}:=\dti{u_{\rho,t}}$,  
there are uniform constants $\wt{c}^P$ and $\wt{c}^B$ 
such that $\wt{c}^B\leq \wt{c}^P$ and 
 \begin{equation}  \label{estimatediffeofibered}\|\zeta_{\rho}\|_{L^{1,p}}\leq \wt{c}^P\|\zeta\|_{L^{1,p}},\hspace{0.5cm} \|D^P_{u_{\rho,0}}\zeta_{\rho}\|_{L^p}\leq \|D^P_{u_0}\zeta\|+ \wt{c}^P|\rho|^{1/2p},
 \end{equation}
 and
% $$\left\|\left.\frac{d}{dt}\right|_{t=0}\delbar_{J_P} u_{\rho,t}\right\|_{L^{p}}\leq \left\|\left.\frac{d}{dt}\right|_{t=0}\delbar_{J_P} u_{t}\right\|_{L^{p}} + C |\rho|^{1/2p}\|\zeta\|_{C^1}.$$
\begin{equation}\label{estimatediffeobase}\|\pi_* \zeta_{\rho}\|_{L^{1,p}}\leq \wt{c}^B\| \pi_*\zeta \|_{L^{1,p}},
\hspa \|D^B_{\pi(u_{\rho,0})}\pi_*\zeta_{\rho}\|_{L^p}\leq \|D^B_{\pi(u_0)}\pi_*\zeta\|+ \wt{c}^B|\rho|^{1/2p}.
\end{equation}

%%%%%%%%%%%%%%%%%%%%%%%%%%%%%%%%%%%%%%%%%%%%%%%%%%%

                                       % INVERSE A DROITE 
                                       
    %%%%%%%%%%%%%%%%%%%%%%%%%%%%%%%%%%%%%%%%%%%%%%%%%                                   

\subsubsection{Right inverses}
We give the description of right inverses for $D^P_u$  which are induced by right inverses for $D^v_u$ and right inverses for  $D^B_{\pi(u)}$. By assumption $D^B_{\pi(u)}$ and $D^{v}_u$ are surjective and we can therefore consider their unique $L^2$-orthogonal right inverses  $Q^B_{\pi(u)}$ and $Q^v_u$,  with respect to  $g_{J_B}$ and $g_{J}$. Set
\begin{equation}\label{matrixQ}
Q^P_u:=\left(\begin{array}{cc}
	(Q^B_{\pi(u)})^h &  0\\
	L'_u & Q^{v}_u
\end{array}\right)
.
\end{equation}
From the matrix expression \eqref{matricelinearisee} for $D^P_u$ we get 
\begin{equation}\label{inversematrice}
D^P_u\circ Q^P_u:=D^P_u\circ\left(\begin{array}{cc}
	(Q^B_{\pi(u)})^h &  0\\
	L'_u & Q^{v}_u
\end{array}\right) = \left(\begin{array}{cc}
	Id &  0\\
	L_u\circ(Q^B_{\pi(u)})^h+ D^{v}_u\circ L'_u & Id
\end{array}\right),
\end{equation}
where $L$ is given by \eqref{connectant}. Thus $Q^P_u$ is a right inverse, if $L':(\mathcal{E}_{P})^h_{u}\sra (\mathcal{X}_{P})^v_{u}$ verifies
\begin{equation*}\label{relationL'}L_u\circ(Q^B_{\pi(u)})^h+ D^{v}_u \circ L'_u=0.
\end{equation*}
A natural choice for $L'$ is  $L'=-Q^v\circ L\circ (Q^B)^h$.

%Consider  the matrix representation of $D_u$  given in equation \eqref{matricelinearisee}.
%\begin{equation}
%\left(\begin{array}{cc}
%	(D^B_{\pi(u)})^h &  0\\
%	L_u & D^{v}_u,
%\end{array}\right)
%\end{equation}
%where $L_u$ stands for the linear operator:
%$$L_u: \mathcal{X}_{P,u}^h\lra\mathcal{E}_{P,u}^v,\hspa \xi\mapsto -\frac{1}{2}J(u)(\nabla_{\xi}J)(\partial_{J_P}u)^v+R^{0,1}(du^h,\xi).$$
\begin{rem}\label{estimateL}
By definition, $L_u$ is  bounded.  Namely, $\|L_u\|\leq C''$ for some constant $C''$ depending on $\|J\|_{C^1},\|du\|_{L^{\infty}}$, and $\|R\|_H$.
\end{rem}

%\begin{equation}\label{L2ortho}\text{Im}(Q^B_{\pi(u)})=(\ker D^B_{\pi(u)})^{\bot_{L^2}}\hspace{0.5cm}\text{Im}(Q^v_u)=(\ker D^v_u)^{\bot_{L^2}},
%\end{equation} 
%o\`u les compl\'ements $L^2$-orthogonaux sont consid\'er\'es par rapport aux normes $g_{J_B}$ et $g_{J}$. 
%\begin{equation}\label{naturalchoice}
%L':=-Q^v\circ L\circ (Q^B)^h. 
%\end{equation}
The following is immediate.

\begin{lem}\label{straightestimate} 
Let $Q^P_u$ be a right inverse for  $D^P_u$ and suppose  $Q_{\pi(u)}^B$ and  $Q_u^v$  are as above . 
Then $L_u'$ is uniquely determined by  \eqref{relationL'} and the requirement that   $Q^P_u$ has for image the $L^2$-orthogonal complement of $\ker D^P_u$. 
In this case, 
$$L'=-Q_u^v\circ L\circ (Q_{\pi(u)}^B)^h.$$ 
In particular if $\|Q_{\pi(u)}^B\|<C^h$ and $\|Q_u^v\|<C^v$, for some  positive constants $C^v$ and $C^h$ depending on   $\|du\|_{L^{\infty}}$, then $\|L'\|<C^hC''C^v$. 
%where $C''$ is the constant in the remark {\rm \ref{estimateL}}.
\end{lem}

\begin{rem} The  $L^2$-orthogonal complementarity  condition is a commodity assumption. The lemma above still holds for different choices of  right inverses as long as we ask that the image of $Q^P$ is given by the images of  $Q^B$ and $Q^v$. 
%Then we again have that the image of $L'$ must lie in the image of $Q^v$.
\end{rem}

Let $p> 2$, and let $(\j,u)\in \wt{\M}_{\S_P}(P)$. From Assumption \ref{strongreg}, Lemma \ref{boundedrightinverseestimates} and Proposition \ref{rightinverse} we have uniform constants, $c^B$ and $c^v$, and   right inverses $Q^B_{\pi(u)_{\rho}}$ and $Q^v_{u_{\rho}}$ for $D^B_{\pi(u)_{\rho}}$ and $D^v_{u_{\rho}}$ such that:
\begin{equation}\label{bricomponents}\|Q^B_{\pi(u)_{\rho}}\eta\|_{W^{1,p}}\leq c^B\|\eta\|_{L^p},\hspa \|Q^v_{u_{\rho}}\eta\|_{W^{1,p}}\leq c^v\|\eta\|_{L^p}.
\end{equation}
 Similarly to Proposition \ref{rightinverse} we have:

\begin{prop}\label{rightinversefibered}
 There exists a   constant  $c^P$   independant of $\rho$ and right inverse $Q^P$ for $D^P$ such that for $(u,\j)\in \M_{\S_P}(P)$:
 \begin{equation}Q^P_{u_{\rho}}:=\left(\begin{array}{cc}
	(Q^B_{\pi(u)_{\rho}})^h &  0\\
	-Q^{v}_{u_{\rho}}\circ L_{u_{\rho}}\circ (Q^B_{\pi(u)_{\rho}})^h & Q^{v}_{u_{\rho}}
\end{array}\right)
\end{equation}
and such that 
\begin{equation}\label{briformulationfibered}\|Q^P_{u_{\rho}}\eta\|_{W^{1,p}}\leq c^P\|\eta\|_{L^p}.
\end{equation}
\end{prop}

\proof The right inverses $Q^B_{\pi(u)_{\rho}}$ and $Q^v_{u_{\rho}}$ are obtained from quasi-inverses $R_{\pi(u)_{\rho}}^B$ for $D_{\pi(u)_{\rho}}^B$, and $R^v_{u_{\rho}}$  for $D^v_{u_{\rho}}$, constructed as in \eqref{quasirightinverse}. From these we construct a quasi-inverse for $D^P_{u_{\rho}}$:
\begin{eqnarray}R^P_{u_{\rho}}:\E^p_{P,u_{\rho}}\equiv \E_{P,u_{\rho}}^{p,h}\oplus \E_{P,u_{\rho}}^{p,v}&\sra&\X^{1,p}_{P,u_{\rho}}\equiv \X^{1,p,h}_{P,u_{\rho}}\oplus \X^{1,p,v}_{P,u_{\rho}}\\
(\eta^h,\eta^v)&\mapsto& ((R^B_{\pi(u)_{\rho}}\pi_*\eta^h)^h+L^R_{u_{\rho}}\eta^h, R_{u_{\rho}}^v\eta^v ).\nonumber
\end{eqnarray}
In fact, 
 \begin{equation}\label{estimationsespacetotalquasiinverse}R^P_{u_{\rho}}:=\Gamma \circ Q^P_{\mathbf{w}_{\rho}}\circ \Lambda.
 \end{equation}
where $\mathbf{w}_{\rho} $, $\Lambda$ and $\Gamma$ are defined as in the $B=pt$ case. Note that $\Gamma$ and $\Lambda$ preserve the splitting induced by the Hamiltonian connection on $TP$,  hence they have the following matrix representation 
\begin{equation*}
\Lambda\equiv\left(\begin{array}{cc}
	\Lambda^h &  0\\
	0 & \Lambda^v
\end{array}\right)
\hspace{2cm}\Gamma\equiv\left(\begin{array}{cc}
	\Gamma^h &  0\\
	0 & \Gamma^v
\end{array}\right).
\end{equation*}
It follows from the matrix form  of $Q^P_{w_{\rho}}$ that
\begin{equation*}R^P_{u_{\rho}}=\left(\begin{array}{cc}
	\Gamma^h(Q^B_{\pi(\mathbf{w}_{\rho})})^h\Lambda^h &  0\\
	-\Gamma^v Q^v_{\mathbf{w}_{\rho}} L_{\mathbf{w}_{\rho}}(Q^B_{\pi(\mathbf{w}_{\rho})})^h\Lambda^h & \Gamma^vQ^{v}_{\mathbf{w}_{\rho}}\Lambda^v
\end{array}\right)\equiv  
\left(\begin{array}{cc}
	R_{u_{\rho}}^h &  0\\
	L^R_{u_{\rho}} & R_{u_{\rho}}^v
\end{array}\right)
\end{equation*}
  we end up with the desired expression for $R^P_{u_{\rho}}$. Note that $R^B\equiv d\pi \circ R^h$. 
  
  We show that $R^P$ is bounded and that  $D^P_{u_{\rho}}R^P_{u_{\rho}}$ is invertible for small enough gluing parameters:
    \begin{equation*}\|D^P_{u_{\rho}}R^P_{u_{\rho}}\eta-\eta\|_{L^p}\leq \frac{C^P}{|\log |\rho||^{1-1/p}}\|\eta\|_{L^p},\hspace{0.5cm}\|R^P_{u_{\rho}}\eta\|_{W^{1,p}}\leq \frac{\wt{C}^P}{2}\|\eta\|_{L^p},
\end{equation*}
for uniform constants $C^P$ and $\wt{C}^P$. But from Lemma \ref{boundedrightinverseestimates} there are  unuiform constants,
$C^B$,  $\wt{C}^B$, $C^v$ and $\wt{C}^v$, such that:
\begin{equation*}\|D^B_{u_{B,\rho}}R^B_{u_{B,\rho}}\eta-\eta\|_{L^p}\leq \frac{C^B}{|\log |\rho||^{1-1/p}}\|\eta\|_{L^p},\hspace{0.5cm} \|R^B_{u_{B,\rho}}\eta\|_{W^{1,p}}\leq \frac{\wt{C}^B}{2}\|\eta\|_{L^p}
\end{equation*}
and
\begin{equation*}\|D^v_{u_{\rho}}R^v_{u_{\rho}}\eta-\eta\|_{L^p}\leq \frac{C^v}{|\log |\rho||^{1-1/p}}\|\eta\|_{L^p},\hspace{0.5cm} \|R^v_{u_{\rho}}\eta\|_{W^{1,p}}\leq \frac{\wt{C}^v}{2}\|\eta\|_{L^p}.
\end{equation*}
Set   $\xi_{\rho}=R^P_{u_{\rho}}\eta$ and  suppose without loss of generality that $|V|=2$, and let $V=\{v,v'\}$ with edge $e_{f,f'}$. Outside the patched  annuli  $D_{f}(r_{vv'}^{1/4})\bs D_{f}(r_{vv'}^{3/4})$, 
%and $D_{f'}(r_{vv'}^{1/4})\bs D_{f'}(r_{vv'}^{3/4})$, 
we have that $\xi_{\rho}=\xi= Q^P_{w_{\rho}}\eta$ and $u_{\rho}=w_{\rho}$ which implies that  $D^P_{u_{\rho}}\xi_{\rho}=\eta$. Therefore, the desired estimate for  $R^P_{u_{\rho}}$ is trivially realized on this part of the curve. Then, it suffices to understand what happens on   $D_{f}(r_{vv'}^{1/4})\bs D_{f}(r_{vv'}^{1/2})$. In that region $u_{\rho}$ and $\mathbf{w}_{\rho}$ are constant with value  $p_{vv'}$. Hence, $D^P_{u_{\rho}}$, $D^P_{u_{v,\rho}}$ and $D^P_{u_{v',\rho}}$ coincide with the standard Cauchy-Riemann operator:
\begin{equation*}
\left(\begin{array}{cc}
	(\delbar_{J_B(\pi(p_{vv'}))})^h &  0\\
	0 & \delbar_{J_{\pi(p_{vv'})}}
\end{array}\right).
\end{equation*}
Furthermore,  $L^R_{u_{\rho}}$ must vanish since $L_{\mathbf{w}_{\rho}}$ vanishes pointwise. Then, the result follows from the estimations of  $R^B$ and $R^v$: the first estimate in \eqref{estimationsespacetotalquasiinverse} is obtained by choosing  $C^P\geq \max (C^B,C^v)$, and the second by taking $\wt{C}^P\geq \max(\wt{C}^B,\wt{C}^v)$.

Hence,   $D^P_{u_{\rho}}R^P_{u_{\rho}}$ is invertible for small enough gluing parameters and we  set:
 $$Q^P_{u_{\rho}}:= R^P_{u_{\rho}}(D^P_{u_{\rho}}R^P_{u_{\rho}})^{-1}.$$ 
Now $D^P_{u_{\rho}}R^P_{u_{\rho}}$ is of the following form:
 \begin{equation*}DR:=\left(\begin{array}{cc}
	D^hR^h &  0\\
	L_{D^PR^P}& D^vR^v
\end{array}\right),
\end{equation*}
where
\begin{equation*}
L_{D^PR^P,u_{\rho}}=L_{u_{\rho}}R^h_{u_{\rho}}-D^v\Gamma^v Q^v_{w_{\rho}} L_{w_{\rho}}(Q^B_{\pi(w)_{\rho}})^h \Lambda^h.
\end{equation*}
Furthermore,
 \begin{equation*}(D^PR^P)^{-1}:=\left(\begin{array}{cc}
	(D^hR^h)^{-1} &  0\\
	L_{(D^PR^P)^{-1}}:=(D^vR^v)^{-1}L_{DR}(D^hR^h)^{-1} & (D^vR^v)^{-1}
\end{array}\right).
\end{equation*}
 Since all the operators involved are lower triangular we must have that
\begin{equation*}Q^P_{u_{\rho}}:=\left(\begin{array}{cc}
	(Q^B_{\pi(u)_{\rho}})^h &  0\\
	L''_{u_{\rho}} & Q^{v}_{u_{\rho}}
\end{array}\right).
\end{equation*}
We identify $L''_{u_{\rho}}$. To simplify notations we will omit the  $u_{\rho}$ subscripts and we will set  $(D^B)^h=D^h$. Again, $L_{\mathbf{w}_{\rho}}$ vanishes in the region  $D_{f}(r_{vv'}^{1/4})\bs D_{f}(r_{vv'}^{3/4})$, 
so that the image of $L_{\mathbf{w}_{\rho}}$ must lie in the image of $\Lambda^v$. Now this latter map is injective therefore
$$L^R=-R^v(\Lambda^v)^{-1}L_{\mathbf{w}_{\rho}}(Q^B_{\pi(\mathbf{w})_{\rho}})^h\Lambda^h.$$
%Set $\wt{L}:=(\Lambda^v)^{-1}L_{w_{\rho}}(Q^B_{\pi(w)_{\rho}})^h\Lambda^h$. By definition
A simple computation then gives:
\begin{equation*}L''=L^R(D^hR^h)^{-1}+R^v L_{(D^PR^P)^{-1}}= -Q^vLR^h(D^hR^h)^{-1}= -Q^vL(Q_B)^h.
\end{equation*}
\qed\\

Finally, let $\{u_t\}_{t\in[0,\upsilon)}$ be a path in $\B^{1,p}_{P}$, where $\zeta:=\dti{u_t}$. Let $\{u_{\rho,t}\}_{t\in[0,\upsilon)}$ be the corresponding path of preglued curves, with $u_{\rho,0}=:u_{\rho}$.
There are uniform constants $\ov{c}^B$ and $\ov{c}^P$ such that:  
\begin{equation} \label{rightinvderfibered}\|\dti Q^B_{\pi(u_{\rho,t})}\|\leq \ov{c}^B\|\pi_*\zeta \|_{W^{1,p}},\hspa \|\dti Q^P_{u_{\rho,t}}\|\leq \ov{c}^P\|\zeta\|_{W^{1,p}}.\end{equation}

%%%%%%%%%%%%%%%%%%%%%%%%%%%%%%%%%%%%%%%%%%%%%%%%%%%

                                       % GLUING
                                       
    %%%%%%%%%%%%%%%%%%%%%%%%%%%%%%%%%%%%%%%%%%%%%%%%%                                   

\subsubsection{Gluing stable components} We  assume here that both $S_{\pi_*}$ and  $\ov{\pi}_{pt}^{\S_P}$  preserves the tree structure of $\S_P$. 
Then $\ov{\pi}_{pt}^{\S_P}\equiv \F_P$, $\ov{\pi}_{pt}^{\S_B}\equiv \F_B$ and $$\ov{\pi}_{\S_P}: \M_{\S_P}(P)\sra \M_{\S_B}(B),\hsp (u,\j)\mapsto (u_B:=\pi(u),\j).$$
This induces a map between  orbibundles:
  $$\ov{\ov{\pi}_{\S_P}}:\mathcal{L}_{\S_P}:=\F_P^*\mathcal{L}_{\S}\sra  \mathcal{L}_{\S_B}:=\F_B^*\mathcal{L}_{\S} \hspa (u,\j,\rho)\mapsto (\pi(u),\j,\rho).$$
  Let $U_P\subset \M_{\S_P}(P)$ and $U_B\subset \M_{\S_B}(B)$ be proper open subsets such that  $\ov{\pi}_{\S_P}(U_P)=U_B$. 
To simplify the exposition, we assume that the maps involved do not have automorphism. If not so,  the gluing maps $Gl_{\S_P}$ and $Gl_{\S_B}$ obtained below,  are  actually defined on local  uniformizing systems for neighbourhoods $W^P$ and  $W^B$ around the stable (holomorphic) maps $(u,\j)$ and $(\pi(u),\j)$,  which are compatible with the projection $\ov{\pi}_{\S_P}$.
%More precisely, the gluing maps are defined on local systems
%$$(\wt{W}^P,\text{Aut}(u,\j),p_{W^P})\and(\wt{W}^B,\text{Aut}(\pi(u),\j),p_{W^B}),$$
%with that the first one projects to the second one using $\pi$. Hence the gluing domains of the gluing maps   will be given by 
%$$\F^*_Pp^*_{W^P}\mathcal{L}_{\S}=\left.\mathcal{L}_{\S_P}\right|_{\wt{W}^P}\and\F^*_Bp^*_{W^B}\mathcal{L}_{\S}=\left.\mathcal{L}_{\S_B}\right|_{\wt{W}^B}.$$
%Nevertheless, we directly see in their  construction, that these gluing maps are  respectively $\text{Aut}(u,\j)$ and $\aut(\pi(u),\j)$ equivariant since the maps are assumed to be stable, so that  they are actually well-defined on the bundles $$\left.\mathcal{L}_{\S_P}\right|_{W^P}\hsp\text{and}\hsp\left.\mathcal{L}_{\S_B}\right|_{W^B}.$$
% 

  From Theorem
\ref{gluingtheorem} there exists a positive constant $\epsilon_B$, and a diffeomorphism
$$Gl_{\S_B}: \mathcal{L}^*_{\S_B,\epsilon_B,U_B} \sra \M_{0,l}(B,\sigma_B),\hsp (u_B,\j,\rho)\mapsto Gl_{\S_P}(u_B,\j,\rho).$$
Similarly to Theorem \ref{gluingtheorem} we have:
%This map is a diffeomorphism obtained obtained by perturbing the pregluing $\text{pgl}^B$

\begin{theorem} For small enough positive constant $\epsilon_P$, there is  a diffeomorphism:
$$Gl_{\S_P}: \mathcal{L}^*_{\S_P,\epsilon_P,U_P}  \sra \M_{0,l}(P,\sigma),\hsp (u,\j,\rho)\mapsto Gl_{\S_P}(u,\j,\rho), $$
such that: 
\begin{equation}\label{commutativitygluingfibered}\ov{\pi}\circ Gl_{\S_P}=Gl_{\S_B}\circ \ov{\ov{\pi}_{\S_P}}.
\end{equation}
\end{theorem}

%%%%%%%%%%%%%%%%%%%%%%%%%%%%%%%%%%%%%%%%%%%%%%%%%%%

                                       % GLUING STABLE
%%%%%%%%%%%%%%%%%%%%%%%%%%%%%%%%%%%%%%%%%%%%%%%%%                                   

\proof First, recall  that  for  $(\pi(u),\j,\rho)\in \L^*_{\S_B,\epsilon_B,U_B}$, 
$$Gl_{\S_B}(\pi(u),\j,\rho)=(\j_{\rho},\exp_{\pi(u)_{\rho}}Q^B_{\pi(u)_{\rho}}(f^B(\pi(u),\j,\rho))),$$ 
where $f^B$ is as in \eqref{transverseelement}:
\begin{equation*}f^B(\pi(u),\j,\rho)\in B_{\delta_B,\pi(u)_{\rho}}(0)\subset L^p(\Lambda_{J}^{0,1}(S^2,\pi(u)_{\rho}^*TB))
\end{equation*}
for a positive (uniform) constant $\delta_B$ given by the Implicit Function Theorem. From the estimates \eqref{almostholomorphicityfibered}, \eqref{quadraticestimatefibered}, \eqref{briformulationfibered},  Implicit Function Theorem applies here, and there are  uniform constants  $\epsilon_P$ and $\delta_P$, and  a smooth map 
$$f^P:\L_{\S_P,\epsilon_P,U_P}\sra B_{\delta_P}(0)$$
such that
\begin{equation*}f^P(u,\j,\rho)\in B_{\delta_P,u_{\rho}}(0)\subset L^p(\Lambda_{J}^{0,1}(S^2,u_{\rho}^*TP))
\end{equation*}
is the unique solution to $\delbar_{J_P}(\exp_{u_{\rho}}Q^P_{u_{\rho}}(f^P(u,\j,\rho)))=0$. Then the gluing map $Gl_{\S_P}$ is defined by;
$$Gl_{\S_P}(u,\j,\rho)=(\j_{\rho},\exp_{u_{\rho}}Q^P_{u_{\rho}}(f^P(u,\j,\rho))).$$ 
The proof that this is a locally diffeomorphic is verbatim  the proof  of Theorem \ref{gluingtheorem} using the estimates \eqref{estimatediffeofibered} and \eqref{rightinvderfibered}.

Next we show \eqref{commutativitygluingfibered}. We begin by showing that
$$d\pi \circ f^P=f^B\circ \ov{\pi}_{\S_P}.$$
By Lemma \ref{projpregl}, $\pi(u_{\rho})=\pi(u)_{\rho}$. 
Let $\xi^P$ denote  $Q^P_{u_{\rho}}f^P(u,\j,\rho)$. Then,
$$\pi(\exp_{u_{\rho}}(\xi^P))=\exp_{\pi(u)_{\rho}}d\pi(\xi^P).$$
Since   $\xi^P$ is the unique solution to  $\delbar_{J_P}\exp_{u_{\rho}}(\xi^P)=0$,
 we  deduce that
$$\delbar_{J_B}\exp_{\pi(u)_{\rho}}d\pi(\xi^P)=0.$$
Moreover,
$$\xi^B:=d\pi(\xi^P)=d\pi(Q^P_{u_{\rho}}f^P(u,\j,\rho))=Q^B_{\pi(u)_{\rho}}d\pi(f^P(u,\j,\rho)),$$
implying that  $\xi^B$  is in the image of  $Q^B_{u_{B,\rho}}$, so  that  $d\pi(f^P(u,\j,\rho))$ is in the image of $f^B$ (by the implicit function theorem and choosing $\delta_P$ smaller than $\delta_B$).
Finally,
 \begin{eqnarray*} 
\ov{\pi}(Gl_{\S_P}(u,\j,\rho))&=&\pi(\exp_{u_{\rho}}Q^P_{u_{\rho}}f^P(u,\j,\rho))\\
 &=&\exp_{\pi(u)_{\rho}}(d\pi (Q_{u_{\rho}}f^P(u,\j,\rho)))\\
 &=& \exp_{\pi(u)_{\rho}}(Q^B_{\pi(u)_{\rho}}f^B(\pi(u),\j,\rho))\\
 &=& Gl_{\S_B}(\pi(u),\j,\rho).
 \end{eqnarray*}
\qed\\

\subsubsection{Gluing: the unstable case.} Let $(u,\j)\in \M_{\S_P}(P)$, and let $v\in V_P$ be a $\pi_{*}$-unstable component. Again, $|v|>2$ is possible. Generalizing the approach in \cite{CL}, we   describe the notion of \emph{balanced component} when $|v|\leq2$.\\

\noindent1) \underline{\emph{ Balanced maps.}}  First assume $|v|=2$. In this case $u_v:\Sigma_v\sra P$ can be parametrized  by $\CP^1$ with special points $0=[1:0]$ and $\infty=[0:1]$. Identify $\Sigma_v\bs\{\infty\}$ with $\C$. If $\pi(u_v)$ is non-constant, we say that $u_v$ is  balanced if it is  \emph{horizontally balanced}, i.e if $\pi(u_v)$ is balanced in the sense of \eqref{balancedconditionenergy}. On the other hand, if $\pi(u_v)$ is constant, we say that  $u_v$ is balanced if it is \emph{vertically balanced}, i.e if  half the vertical energy ($\hbar^v$) of $u_v$ is contained in the unit disc around 0.  
  
  Next, when $|v|=1$, $\Sigma_v$ can be parametrized  by $\CP^1$ with special $\infty=[0:1]$.  Identify $\Sigma_v\bs\{\infty\}$ with $\C$. 
 If  $\pi(u_v)$ is not constant,  $u_v$ is  \emph{balanced} if it is \emph{horizontally balanced} 
    as in the $|v|=2$ case, and if $u_v$ is \emph{horizontally centered}, i.e  if the center of energy of  $\pi(u_v)$ is zero. If $\pi(u_v)$ is constant, we say that  $u_v$ is balanced if it is \emph{vertically balanced} as above, and if  \emph{vertically centered}, i.e if  the mean  value for  $\|du_v^v\|_{g_{J_P}}^2$ is zero. 
    %(see \eqref{balancedconditioncenterofenergy}). 

 For $i=1,2$,  let $\wt{\M}^{b,h}_{0,i}(P,\sigma_v)$, resp.  $\wt{\M}^{b,v}_{0,i}(P,\sigma_v)$, denote   the set of horizontally, resp. vertically,  balanced  $J_P$-holomorphic maps with $i=1$ or 2  marked point representing $\sigma_v$. In order to simplify notations, $\wt{\M}^b_{0,i}(P,\sigma_v)$ will designate the set of  (vertically and horizontally) balanced  $J_P$-holomorphic maps with $i$ marked point. 
When $\pi_*\sigma_v\neq 0$, the projection $\pi$ naturally induces a map:
\begin{equation}\label{mapbetweenbalancedmoduli} \pi:\wt{\M}^{b}_{0,1}(P,\sigma_v)=\wt{\M}^{b,h}_{0,i}(P,\sigma_v)\sra \wt{\M}^{b}_{0,i}(B,\pi_*\sigma_v).
\end{equation}
As in the $B=pt$ case, there is a natural smooth surjective map between the moduli space of holomorphic maps to  the moduli space of balanced maps,  $u_v\mapsto u_v^b$, where $u_v^b$ is given by \eqref{associatebalancedmap},  and it is not hard to see that 
\begin{equation}\label{compatibilityassociatebalancedmap}\pi(u_v)^b=\pi(u_v^b)\hsp \text{if} \hsp \pi_*\sigma_v\neq0.
\end{equation}
We conclude that
%$$\phi_P^b: \wt{\M}_{\j_1}(P,\sigma)\lra \wt{\M}^b_{0,1}(P,\sigma),\hspa u\mapsto u\circ(\ov{\phi}_P^b(u))$$
$$\M_{0,i}(P,\sigma_v)\cong\left.\wt{\M}^{b}_{0,i}(P,\sigma_v)\right/S^1,\hsp i=1,2,$$
and when $\pi_*\sigma\neq 0$, the map $\pi$ in \eqref{mapbetweenbalancedmoduli} descends to the (expected) map:
$$\ov{\pi}:\M_{0,i}(P,\sigma_v)\sra \M_{0,i}(B,\pi_*\sigma_v), \hsp i=1,2.$$ 
Finally,  we say that $(u,\j)\in\wt{\M}_{\S_P}(P)$ is balanced if each of its $\pi_*$-unstable components $u_v$ with $|v|\leq2$ is balanced. The corresonding moduli space of balanced $J_P$-holomorphic maps  is denoted $\wt{\M}^b_{\S_P}(P)$. Similarly to  \eqref{definitionmodulispacebalancedmaps}, we have a natural identification,
\begin{equation*}\label{definitionmodulispacebalancedmapsfibered} \M_{\S_P}(P)=\left.\wt{\M}^b_{\S_P}(P)\right/\aut^P_{red},
\end{equation*}
where  $\aut_{red}^P$ denotes the group of reparametrizations acting on  $\wt{\M}^b_{\S_P}(P)$ (see \eqref{reducedgroup}). 
% a product of copies of $S^1$, one for each $\pi_*$-unstable components in $\S_P$,  with the action of $\S_P$.   
%the reparametrization group reduces to that of  $(S^1)^{I^u}$ combined with the action of $\aut(\S_P)$, where $I^u$ is the subset of $I$ corresponding to the unstable components of the domain. Since $\aut(\j)$ fibers over  $\aut(\S_P)$, which is finite,  the group acting is  given by $(S^1)^{I^u}\times\aut(\S_P)$ which we denote by  $\aut_{red,P}$. 
% where $\aut_{red,B}$ corresponds to $(S^1)^{I^u\cap I^h}\times\aut(\S_B)$.
%\begin{equation*}\aut_{red,P}:=(S^1)^{I^u}\times\aut(\S_P)\and\aut_{red,B}:=(S^1)^{I^u\cap I^h}\times\aut(\S_B).
%\end{equation*}
Observe  that  the projection $\pi_{\S_P}$  restricts to an $\aut^P_{red}$, $\aut^B_{red}$ equivariant map $$\pi_{\S_P}:\wt{\M}^b_{\S_P}(P)\sra \wt{\M}^b_{\S_B}(B).$$
%\begin{equation*}\F_{\pi}:\wt{\M}^b_{\S_P}(P)\lra \wt{\M}^b_{\S_B}(B),
%\end{equation*}
and descends to $\ov{\pi}_{\S_P}:\M_{\S_P}(P)\sra \M_{\S_B}(B)$ after quotienting by  the reparametrizations. 
%we recover the map
%%$(S^1)^{I^u}$ et $\aut(\S_P)$. 
%$$\F_{\pi}:\M_{\S_P}(P)\lra \M_{\S_B}(B).$$

\begin{rem}   Let $P$ denote  the bundle $(S^2\times S^2,\omega_0+\omega_0)$ over the base   $B:=(S^2,\omega_0)$, where $\om_0$ denotes the Fubini-Study form on $S^2$.  Here $\pi$ represents the projection to the first factor, $J_P$ is the product complex structure. Consider the holomorphic section $u(z)=(z, z+b)$.  
%We have seen that $C^B(\pi(u))=0$ and $E_B(\pi(u), B_1(0))=\pi/2=\hbar^h$,  i.e the projection 
A simple computation shows that $\pi(u)$ is balanced  (the map $z\mapsto az+b$  is balanced if and only if  $|a|=1$ and $b=0$).  In addition, 
%we can see that $C^v(u)=b$ and $E^{vert}(u, B_1(b))=\pi/2=\hbar^v$, so that 
one sees that $u(z-b)$ is vertically balanced. 
Note that if we had adopted the  definition for balanced maps in $P$  with respect to the energy density $\|du\|^2_{g_{J_P}}$, the mapping \eqref{mapbetweenbalancedmoduli} would not necessarily exist. Indeed,  the map 
%$$C= \frac{1}{E(u)}\left(E^{vert}(u)C^v(u)+E_B(u_B)C^B(u_B)\right)=\frac{1}{2\pi}(\pi b)= \frac{b}{2},$$
%and we easily see that $E(u,B^1(b/2))= \hbar$. Therefore 
$u(z-b/2)$ is balanced in this sense, but projects to $z-b/2$ which is not balanced.

Finally, had we defined a balanced map in  $P$ as being both horizontally and vertically balanced, the compatibility \eqref{compatibilityassociatebalancedmap} may not be realized. For example,   consider the maps  $u(z)=(z,az)$, for $a\in \R^+$, with $a\neq1$. Then the center of energy of  $u$ is  0,  so that it is horizontally balanced. However, the map is not vertically balanced, and if $\hbar^h$ denotes the energy of the projection $\pi(u)$, then  the energy of $u$ in the unit disc  around 0 is given by $\hbar^h+\pi(1-\frac{1}{1+a^2})$,
%$$E(u,B_{1}(0))=E_B(u_B,B_{1}(0))+E^{vert}(u,B_1(0))=\hbar^h+\pi\left(1-\frac{1}{1+a^2}\right),$$
which gives  $\hbar$ if and only if $a=1$. 
\end{rem}

%More explicitly let $\{\j_i\}_{i\in I}$ be the normalization of the nodal domain $\j$,
%%set
%%$$\j:=\{\j_i\}_{i\in I}/\sim:=\bigsqcup_{i\in I}(\Sigma_i, \{y_{ij}\}_{iEj}, \{x_m\}_{m\in D^{-1}(i)})/_{y_{ij}\sim y_{ji}},$$
%and let $I^u$ be the subset of $I$ corresponding to the unstable components of the domain.
%% $$I^u:=\left\{i\in I\,\,|\,\, \j_i \hs\text{is stable}\right\},$$
%Let $I^s$ be  the complementary of  $I^u$  in $I$. Recall that   $\wt{\M}_{\S_P}(P,J_P)$ is a submanifold of
%$$\prod_{i\in I}\widetilde{\M}_{\j_i}(P,\sigma_i)=\prod_{i\in I^s}\widetilde{\M}_{\j_i}(P,\sigma_i)\times\prod_{i\in I^u\cap I^h}\widetilde{\M}_{\j_i}(P,\sigma_i)\times\prod_{i\in I^u\cap I^v}\widetilde{\M}_{\j_i}(P,\sigma_i),$$
%verifying the appropriate incidence  relations. In these notations, the map  $\phi^b_P$ induces a map:
%\begin{eqnarray*} \phi^b_P:\prod_{i\in I}\widetilde{\M}_{\j_i}(P,\sigma_i)&\sra&\prod_{i\in I^s}\widetilde{\M}_{\j_i}(P,\sigma_i)\times\prod_{i\in I^u\cap I^h}\wt{\M}^{b,h}_{\j_i}(P,\sigma_i)\times\prod_{i\in I^u\cap I^v}\wt{\M}^{b,v}_{\j_i}(P,\sigma_i)\\
%\{u_i\} &\mapsto& (\{u_i\}_{i\in I^s}, \{\phi^b_P(u_i)\}_{i\in I^u\cap I^h},\{\phi^b_P(u_i)\}_{i\in I^u\cap I^v})
%\end{eqnarray*}
%
%\begin{defn}\label{balmapstrat}  The set  $\wt{\M}^b_{\S_P}(P,J_P):=\phi^b_P(\wt{\M}_{\S_P}(P))$ is called moduli space of balanced $J_P$-holomorphic maps for the stratum data $\S_P$.
%\end{defn}

%%%%%%%%%%%%%%%%%%%%%%%%%%%%%%%%%%%%%%%%%%%%%%%%%%%

                                       % GLUING NON STABLE 
                                       
    %%%%%%%%%%%%%%%%%%%%%%%%%%%%%%%%%%%%%%%%%%%%%%%%%                                   

\noindent 2) \underline{\emph{ Gluing of balanced maps and compatibility.}}  Let  $\S_B^u$ and $\S^u$ respectively denote the projections of $\S_P$ under $\pi_*$ and  $\F_P$. These stratum data are not stable here. Note that  $\S_B^u$ and  $\S_P$ have the same tree structure but the homological data for $\S_B^u$ is the projection of the homological data for $\S_P$ under $\pi_*$. 

Let $\wt{\L}_{\S_P}$ and $\wt{\L}_{\S_B}$ denote the bundles of gluing parameters over the balanced moduli spaces $\wt{\M}_{\S_P}^b(P)$ and $\wt{\M}^b_{\S_B}(B)$  obtained by pull-backing $\wt{\L}_{\S^u}$ under $\F_P$ and $\F_B$ respectively. Consider the bundle map lifting $\ov{\pi}_{S_P}$:
$$\ov{ \ov{\pi}_{\S_P}}: \wt{\L}_{\S_P}\sra\wt{\L}_{\S_B},\hspa (u,\j,\rho)\mapsto (\ov{\pi}_{\S_P}(u,\j), \rho^{st}),$$
where $\rho^{st}$ denotes the stabilization of $\rho$ with respect to $\pi_*$ this time (see \eqref{stabilisationgluingparameters}). 
This map descends  to a  well-defined map:
$$\ov{ \ov{\pi}_{\S_P}}:\L_{\S_P}:=\left.\wt{\L}_{\S_P}\right/\aut_{red}^P\sra \L_{\S_B}:=\left.\wt{\L}_{\S_B}\right/\aut_{red}^B.$$

Let $U_P\subset \M_{\S_P}(P)$ and $U_{B}\subset \M_{\S_B}(B)$ be proper open subsets such that $\ov{\pi}_{\S_P}(U_P)=U_B$.
From the discussion in Section \ref{sectiongluingunstable} we have a gluing map,
$$Gl^b_{\S_P}:\left.\wt{\L}^*_{\S_P,\epsilon_P,U_P}\right/\aut_{red}^P\sra \M_{0,l}(P,\sigma),$$
and a  gluing map $Gl^b_{\S_B}$ above  $U_B$ defined similarly. We prove the following generalization of  Theorem \ref{TheoremChenLi}:

%The forgetting-the-map map  induces a bundle map:
%$$\ov{\F}_X:\wt{\L}_{\S_X}\sra\wt{\L}_{\S},\hspa (u,\j,\rho)\mapsto (\pi(u),\j, \rho).$$
%Let $\wt{\L}^b_{\S_X}$ denote its restriction to $\wt{\M}^b_{\S_X}(X)$.
%The group  $\aut^X_{red}$ acts on  each fiber of $\wt{\L}_{\S_X}$ by rotations,
% therefore   $\aut_{red}^X$ acts naturally  on $\wt{\L}^b_{\S_B}$. 
%By hypothesis \eqref{strongreg}, for any element of $\wt{\M}^b_{\S_X}(X)$, the linearization of the  Cauchy-Riemann operator is surjective.  It follows from Theorem \ref{gluingtheorem}  that for any proper open subset $U_X$ of $\wt{\M}^b_{\S_X}(X)$, that we may choose to be  $\aut_{red}^X$-invariant, there exists $\epsilon_X>0$ and a gluing map 
%$$\wt{Gl}^b_{\S_X}:\wt{\L}_{\S_X,\epsilon_X,U_X}\sra \wt{\M}_{0,l}(X,A).$$
%This map is $\aut_{red}^X$-equivariant, thus the gluing map is well-defined after the quotient by the action:
%$$Gl^b_{\S_X}:\left.\wt{\L}_{\S_X,\epsilon_X,U_X}\right/\aut_{red}^X\sra \wt{\M}_{0,l}(X,A).$$
% Let  $\wt{\M}_{\S_B^u}(B)$ denote the set  of  parametrized nodal $J_B$-holomorphic maps having $\S_B^u$ as stratum data and let $\wt{\L}_{\S_B^u}$ be the corresponding fiber bundle of  gluing parameters: $\wt{\L}_{\S_B^u}$ is the pull-back  of $\wt{\L}_{\S^u}$ under the forgetting-the-map map from $\wt{\M}_{\S_B^u}$ to $\wt{\M}_{\S^u}$. 
 
%   Let $(u,\j) \in\wt{\M}^b_{\S_P}(P)$,   $(u_B,\j_B):=\pi_{\S_P}(u,\j)\in \wt{\M}^b_{\S_B}(B)$, and  $(\pi(u),\j) \in\wt{\M}^b_{\S^u_B}(B)$.  

\begin{theorem} 
The gluing maps  $Gl^b_{\S_P}$ and $Gl^b_{\S_B}$ are locally diffeomorphic and such that 
\begin{equation}\label{compatibilityfibered}\ov{\pi}\circ Gl^b_{\S_P}=Gl^b_{\S_B}\circ \ov{\ov{\pi}_{\S_P}}.
\end{equation}
\end{theorem}

\proof  The diffeomorphic issue follows directly from  Theorem \ref{TheoremChenLi}. The compatibility \eqref{compatibilityfibered}  is obtained as follows.  Let  $\wt{\M}^b_{\S_B^u}(B)$ denote the image of $\wt{\M}^b_{\S_P}(P)$ under $\pi$: $(u,\j)\mapsto (\pi(u),\j)$. Note that the parametrization of the domain $\pi(u)$ comes from the parametrization of the domain of  $u$, which is fixed up to an  $S^1$-action. Now, let $\wt{\L}^b_{\S_B^u}$ denote the pull-back of  $\wt{\L}_{\S^u}$ under the forgetting-the-map map. Then,  by assumption $\eqref{strongreg}$, we can  construct a gluing map $\wt{Gl}^b_{\S^u_B}$ with value in $\wt{\M}_{0,l}(B,J_B)$, and domain a neighbourhood of the zero section in $\wt{\L}_{\S_B^u}$ restricted to $\pi(U_P)$. This map is equivariant with respect to the balanced maps automorphisms, but since $\S^u_B$  is not  stable, it is not locally diffeomorphic nor  injective. Nevertheless,  $\wt{Gl}^b_{\S^u_B}$ still verifies the compatibility:
  $$\pi\circ \wt{Gl}^b_{\S_P}=\wt{Gl}^b_{\S^u_B}\circ \pi.$$
 Again, the problem is localized in the contracted branches and the connecting branches.  %\begin{enumerate}[i)]
% \item[1)]  $\j$ is only composed by one principal component, $\Sigma_0$ and one branch $Br_0$ that we order in an obvious way starting from the root.
% \item[2)]  $\j$ is made of two principal components $\Sigma_1$, $\Sigma_2$ and a connecting branch  $C^*_{1,2}$ between those two components.
% \end{enumerate}
 
 \underline{\emph{Contracted branch}}.  Set $(u_B,\j_B):=\pi_{\S_P}(u,\j)$.
As in Theorem \eqref{TheoremChenLi},  for every $\rho$ the glued surface  $\j_{\rho}$ is isomorphic  to   $\j_B$, i.e  to the domain of the $\pi_*$-stable component. Therefore, $u_B$ and $\wt{Gl}^b_{\S_P}(u,\j,\rho)$ have the same domain. Fix  $\rho_0$ with small radius. We show that
$$ \wt{Gl}^b_{\S^u_B}(\pi(u),\j,\rho_0)=(u_B,\j_B).$$
%We already have that $\pi(\wt{Gl}_{\S_P}(u,\j,\rho_0))= \wt{Gl}_{\S^u_B}(\pi(u),\j,\rho_0)$. 
By definition, the pregluing of $(\pi(u),\j,\rho_0)$ coincides with $\pi(u)$ except on disc with radius determined by $\rho_0$ on which
$$\pi(u)_{\rho_0}(z)\equiv\prgl_B(\pi(u),\j,\rho_0)(z)=\exp_{\pi(p_{\infty})}\left(\beta(z/\rho_0^{1/4})\exp^{-1}_{\pi(p_{\infty})}(\pi(u_0(z))\right),$$
 $p_{\infty}$ being the image under $u$ of the point of the root identified to $\infty$.  Set
$$\xi:=\exp^{-1}_{\pi(p_{\infty})}(\pi(u_0(z)))\and \beta_{\rho_0}:=\beta(z/\rho_0^{1/4}).$$
Then,  $\wt{Gl}_{\S_B^u}(\pi(u),\j,\rho_0)$ coincides with $u_B\equiv\exp_{\pi(u)_{\rho_0}}\((1-\beta_{\rho_0})\xi\)$.
  % Note that $\pi(u)_{\rho_0}(z)$ coincides with $\pi(u)_{0,\rho_0}$, therefore for  $p>2$:
%$\|\delbar_{J_B}(\pi(u)_{\rho_0})\|_{L^p}\leq c|\rho_0|^{1/2p}=O(|\rho_0|^{1/2p}),$$
%where $c$ is  independant from $|\rho_0|$ (cf \ref{estimatepregl}). 
%and we wish to show that it is in the image of the section $f^B$ of the implicit function theorem.  
%Furthermore, $\pi(u)_{0,\rho_0}$ is $W^{1,p}$ close to $\pi(u_0)$ and (cf \cite{MS} lemma 10.4.3):   
%$$\|\xi-\beta_{\rho_0}\xi\|_{W^{1,p}(D(1))}\leq c'\|\xi\|_{W^{1,p}(D(2|\rho_0|^{1/4}))}=O(|\rho_0|^{1/4})$$
%for  $c'$ not depending on $\rho_0$.  Therefore,   for  sufficiently small $|\rho_0|$, one has that  
%$$|\rho_0|^{1/4}<< |\rho_0|^{1/2p}$$
%so that  $\|\xi-\beta_{\rho_0}\xi\|_{W^{1,p}}< 2 \epsilon_1^B c_1^B$, where as 
%$$\frac{2c_1^B c |\rho_0|^{1/2p}}{2 c' |\rho_0|^{1/2p}}=O(|\rho_0|^{\frac{2-p}{4p}})>> 1$$   \\

\underline{\emph{Connecting branch}}.  As in Theorem \eqref{TheoremChenLi},  we can assume the contracting branch is a connecting chain.
%Again $\pi(u)$ is balanced and here
%$$\aut(\j)=\prod_{k\in C^*_{1,2}}S^1.$$  
 By adding $k$-marked point, one for  each component of the connecting branch,  $\pi(u)$  becomes stable. Let $\j'$ denote the conformal structure resulting from this operation. 
Let  $\wt{\M}^b_{\S^u_B(k)}(B)$  denote the stabilization by adding $k$-marked points. 
  We have a gluing map
% $$Gl_{\S^u_B}: N_{\S^u_B(k),\pi(u)}\times V_{\rho_0}\lra \M_{0,l+k}(B,\sigma_B)$$
$$\wt{Gl}_{\S^u_B(k)}: \wt{\L}^*_{\S^u_B(k)}\equiv\F_B^*\wt{\L}_{\F_{B}(\S^u_B(k))}\sra \wt{\M}_{0,l+k}(B,\sigma_B).$$
By definition of the gluing, and since the components of  $\pi(u)$ coming from the connecting  branch are constant, the maps
$$\wt{Gl}_{\S^u_B(k)}(\pi(u),\j',\rho)\and\wt{Gl}^b_{\S^u_B}(\pi(u),\j,\rho)$$
 coincide for every $\rho$, and any choice of balanced parametrization on the domain $\Sigma$ of $\j$ (here the $\aut(\j)$-equivariance of $ Gl_{\S^u_B}$ is needed). Note that  $\j_{\rho}$ is obtained from $\j'_{\rho}$ by forgetting the added marked points.  We  compare the gluing of  $(\pi(u),\j')$ with that of $(u_B,\j_B):=\pi_{\S_P}(u,\j)$.  
Note that $(u_B,\j_B)$ is the image of  $(\pi(u),\j')$ via the map that forgets the added marked points and stabilizes.
% $$\F_{\S_B^u(k)}:\wt{\M}^b_{\S^u_B(k)}(B)\sra \wt{\M}^b_{\S_B}(B).$$
The fiber over $(u_B,\j_B)$ corresponds to the set of all possible reparametrizations for $\j'$, since  $u_B$ and $\pi(u)$ have same image. 
Thus, it suffices to identify the gluings of the corresponding domains but this gollows fromTheorem 
\eqref{TheoremChenLi}.  Finally, 
\begin{eqnarray*}\pi\circ \wt{Gl}_{\S_P}(u,\j,\rho)&=&\wt{Gl}_{\S_B^u}(\pi(u),\j,\rho)\\
&=&\wt{Gl}_{\S^u_B(k)}(\pi(u),\j',\rho)=\wt{Gl}_{\S_B}(u_B,\j_B,\rho^{st})\\
&=&\wt{Gl}_{\S_B}\circ \ov{\ov{\pi}_{\S_P}}(u,\j,\rho).
\end{eqnarray*}
\qed

From now on, we will always assume  that  gluing maps are obtained  by considering balanced maps, and we will drop the $b$ indices from the  notations for the gluing maps.  

%\begin{theorem}{\rm (}B. Chen, A-M. Li{\rm )}  The map $Gl^b_{\S_B}$ (resp. $Gl^b_{\S_P}$) is a local 
%diffeomorphism. 
%\end{theorem}

\subsubsection{Gluing maps between stratas}
One can generalize  the preceding discussions and introduce gluing maps between different stable strata.  Let $\S_X$ and $\S_X'$ be stable stratum data for $\ov{\M}_{0,l}(X,A)$,   such that $\S_X\prec\S_X'$. Let $\S$ and $\S'$ denote their projections under $\F_X$, and consider the bundles $\L_{\S,\S'}$ and $\L_{\S_X,\S'_X}:=\F_X^*\L_{\S,\S'}$. For an open proper subset $U_X$ of $\M_{\S_X}(X)$ there exists  a positive constant $\epsilon_X$ and locally diffeomorphic map:
\begin{equation*}\label{deuxgluing}Gl_{\S_X,\S'_X}:\L^*_{\S_X,\S'_X,\epsilon_X,U_X}\sra\M_{\S'_X}(X),
\end{equation*}
which coincides with the identity on the zero section. Also, from the definition  of $\L_{\S_X,\S'_X}$, a point of $\L_{\S_X}$  is locally given by a tuple $(u,\j,\rho_1,\rho_2)$ 
where $(u,\j,\rho_1)\in\L_{\S_X,\S'_X} $ and where $\rho_2$  accounts for the remaining gluing parameters. 
Therefore,  $Gl_{\S_X,\S_X'}$ induces a map
\begin{equation*} 
\L_{\S_X}\sra\L_{\S_X'},\hspa(u,\j,\rho_1,\rho_2)\mapsto (Gl_{\S_X,\S'_X}(u,\j,\rho_1), \rho_2)
\end{equation*} 
It follows that  $\L_{\S_X}$  coincides with the pullbacks $Gl^*_{\S_X,\S'_X}\L_{\S'_X}$.   Suppose now we have a third stratum data $\S_X''$ such that $\S_X'\prec\S_X''$. Since $Gl_{\S_X,\S_X'}$ is locally diffeomorphic, we can define a new gluing map:
\begin{equation*}\label{troisgluing}Gl'_{\S'_X,\S''_X}:= Gl_{\S_X,\S_X''}\circ Gl^{-1}_{\S_X,\S_X'}:\L^*_{\S'_X,\S''_X,\epsilon_X,\text{Im}(Gl_{\S_X,\S_X'})}\sra\M_{\S'_X}(X),
\end{equation*}
extending the identity map on the zero section. This new gluing does not necessarily coincide with $Gl_{\S'_X,\S''_X}$. The equality  would mean that the gluing procedure is associative, which  is a priori not   true due to the numerous choices made along the gluing construction (in particular the independance with respect to the choice of right inverses). Nevertheless we can see that these maps are close, in the  $C^{\infty}$ sense, which is enough to give  the moduli spaces the structure of smooth orbifolds, as we will see in the next Section. 

Now consider the Hamiltonian Fibration case $\pi:P\sra B$.  Let $\S_P$ and $\S_P'$ be stable stratum data for $\ov{\M}_{0,l}(P,\sigma)$ such that $\S_P \prec \S_P'$,  and let $\S_B$ and $\S_B'$ be their corresponding projections via $S_{\pi_*}$. 
We see that  $\S_B\prec\S_B'$.  
% and we will assume that $\S^u=\S$ so that: 
%$$(\S')^u=\S',\hsp   \S_B=\S_B^u,\hsp\S'_B=(\S'_B)^u.$$
By the discussion above,  for  $U_P$ and $U_B$ be  proper open subsets of  $\M_{\S_P}(P)$ and $\M_{\S_B}(B)$ such that $\ov{\pi}_{\S_P}(U_P)=U_B$ we do have gluing maps $Gl_{\S_B,\S'_B}$ and $Gl_{\S_P,\S'_P}$ such that:
  $$Gl_{\S_B,\S'_B}\circ \left.\ov{\ov{\pi}_{\S_P}}\right|_{\L_{\S_P,\S_P'}}=\ov{\pi}_{\S_P'}\circ Gl_{\S_P,\S'_P}.$$
   Suppose now we have a third stratum data $\S_P''$ projecting on  $\S''_B$ and such that $\S_P'\prec\S_P''$. Then we also have the commutativity: 
   $$Gl'_{\S'_B,\S''_B}\circ \left.\ov{\ov{\pi}_{\S'_P}}\right|_{\L_{\S_P',\S_P''}}=\ov{\pi}_{\S''_P}\circ Gl'_{\S'_P,\S''_P}.$$

\subsection{Fibration of moduli spaces} In this section we prove that  a Hamiltonian Fibration structure  induces a fibration structure between the appropriate compactified moduli spaces: 

\begin{theorem} Under hypothesis  {\rm(\ref{strongreg})}, the moduli spaces  $\ov{\M}_{0,l}(P,\sigma)$ and $\ov{\M}_{0,l}(B,\sigma_B)$ are smooth orbifolds, and the maps $\ov{\pi}_{\S_P}$ extend to a map
$$\ov{\pi}:\ov{\M}_{0,l}(P,\sigma)\sra \ov{\M}_{0,l}(B,\sigma_B),$$
which restricts to a smooth locally trivial fibration (of orbifolds) above each strata of  $\ov{\M}_{0,l}(B,\sigma_B)$.
\end{theorem}
 
 Regarding the fibration structure, it suffices to show this above a proper open subset $U_B$. But each $\ov{\pi}_{\S_P}$ is a smooth submersion. Moreover, the fibers of $\ov{\pi}$  are  compact, hence $\ov{\pi}$ is proper which ends the proof of the fibration statement. 
 
Following \cite{CL}, we show  below how the  compactified moduli spaces can be given the structure of smooth orbifolds compatibly with the $\ov{\pi}_{\S_P}$.

\subsubsection{Charts data and admissible gluing maps}

Consider a Fredholm system  $(\B,\E,s)$ with moduli space $\M=s^{-1}(0)$, modeled on maps. Assume that the linearization $L_{x_0}$ at $x_0\in\M$ is surjective. Then, a standard construction gives a local coordinate chart  around $x_0$. Such a chart is given by a triple $(U,\phi, f)$ where:
\begin{enumerate}[(i)]
\item $U$ is a submanifold of  a neighbourhood $V_{x_0}$ of zero in $T_{x_0}\B$ (which  we identify to a neighbourhood $V_{x_0}$ of $x_0$ in $\B$ via the exponential map);
\item $\phi:U\times B_{\delta}\sra V'_{x_0}$ is  a diffeomorphism where $B_{\delta}\subset \E_{x_0}$ is an open ball, $V'_{x_0}$ is a neighbourhood of 0 in $V_{x_0}$;
\item $f$ is a smooth section $f:U\sra B_{\delta}$;
\end{enumerate}
with the property that
$$\Psi:U\stackrel{1\times f}{\lra}U\times B_{\delta}\stackrel{\phi}{\lra} V'_{x_0},$$
is a diffeomorphism from  $U$ onto $V'_{x_0}\cap \M$. Here $x_0$ serves as a reference point.  More generally, fix $x_0\in \B$, which may not belong to $\M$,  and let $V_{x_0}\subset \B$ be  a neighbourhood of $x_0$. 
  
  \begin{defn} A triple
$(U,\phi, \Psi)$, or $(U,\phi,f)$, verifying conditions i), ii) and iii) above is a called \emph{chart data}  for $\M$.
\end{defn}
 
%\begin{prop} Two charts coming from different chart datas are  $C^{\infty}$ compatible.
%\end{prop}
\begin{rem} It follows  immediately from the definition,  that  the triples $(U_P,f^P,Gl_{\S_P})$ and $(U_B,f^B,Gl_{\S_B})$  give charts datas for $\M_{0,l}(P)$ and $\M_{0,l}(B)$. Moreover, these are compatible with $\ov{\pi}$.  
\end{rem}

In fact, one can construct  chart data for  $\M_{\S_P}(P)$ from pairs, $\mathcal{Q}_P:=(U_P,Q^P)$, where $U_P$ is a smooth submanifold of $\B^{1,p}_{\S_P}$, and  where
$$Q^P:=\{Q^P_{u}|u\in U_P\}$$
is  a smooth $U_P$-family of right inverses for $D^P_u$. In order to do so, we assume the following conditions (that we actually met when constructing the gluing maps):
\begin{assumption}\label{hyp}
\begin{enumerate}[$\bullet$]
\item  for every $u\in U_P$,  $$\|du\|_{L^p}\leq C_P\hspa\text{and}\hspa \|\delbar_{J_P}u\|_{L^p}\leq\epsilon_P,$$
\item for all $\xi\in T_{u}U_P$, 
$$\|\frac{d}{d\xi}\delbar_{J_P}u\|_{L^p}\leq\epsilon_P \|\xi\|_{W^{1,p}},$$
\item the family $Q^P_u$ is  Lipschitz continuous for  the  constant $C_P$ and for all  $u\in U_P$:
$$\|Q^P_u\|\leq C_P,$$
\end{enumerate}
 the constants $C_P$ and $\epsilon_P$ being such that $C_P\epsilon_P <<1$, with $\epsilon_P$ small. 
\end{assumption}
 
The chart data is then given as follows. 
%Under these assumptions gluing maps are defined and yield local diffeomorphisms.
Fix $u_0\in U_P$, and let $W$ be  a neighbourhood of 0 in  $T_{u_0}\B_{\S_P}^{1,p}$. 
Denote by  $U_P$, the lift  (around $u_0$) of $U_P$ in  $W$ via $\exp^{-1}_{u_0}$. 
 Then set
 $$\phi_P:U_P\times \E^p_{P,u_0}(\S_P)\sra \X_{P,u_0}^{1,p},\hspa (\xi,\eta)\mapsto \xi+Q^P_{u_0}\eta.$$
From the assumptions above,  there is a  unique smooth map  $f^P: U_P\sra B_{\delta}$, around $u_0$, such that 
$$\delbar_{J_P}\exp_{u_0}\phi_P(\xi,f^P(\xi))=0.$$
By reducing $\delta$ and the neighbourhoods involved, we further have that $\phi_P$ is 
diffeomorphic, hence $(U_P,\phi_P,f^P)$ is a chart data for $\M_{\S_P}(P)$.

\begin{rem} Regarding compatibility of the coordinate charts one needs to be carefull, as pointed out in \cite{Rvirt}  Section 3. In fact, the $C^{\infty}$ compatibility is ensured if we restrict our attention to smooth stable maps, which is sufficient to study pseudo-holomorphic stable maps  (by elliptic regularity). 
\end{rem}

From the data $\mathcal{Q}_P$, we can furthermore define another type of gluing map. We explain this. Suppose that $(U_P,\phi_P,f^P)$ is a chart data for a proper subset $\wt{U}_P$ of $\M_{\S_P}(P)$, i.e $\wt{U}_P$ is the image of the diffeomorphism: $\Psi_P:=\phi_P\circ (1\times f^P)$. By further reducing $U_P$ if necessary, we can find a pair as above
$$\mathcal{Q'}_P=(U'_P:=\prgl(\Psi_P^{-1}(\wt{U}_P)), Q')$$
where $Q'$ is a family of right inverses for the elements in $U'_P$, which is constructed from the original family of right inverses $Q^P$, 
giving a  gluing map defined by the composition:
$$Gl_{\mathcal{Q}_P}:\left.\mathcal{L}^*_{\S_P}\right|_{\wt{U}_P}\lra \Psi_P^*\left.\mathcal{L}^*_{\S_P}\right|_{U_P}\stackrel{Gl}{\lra} \M_{0,l}(P,\sigma). $$
%$$Gl_{\mathcal{Q}_P}:\left.\mathcal{L}_{\S_P}\right|_{\wt{U}_P}\sra \M_{0,l}(P,\sigma).$$
%More explicitly $Gl_{\mathcal{Q}_P}$ is given by the composition:
%$$\left.\mathcal{L}_{\S_P}\right|_{\wt{U}_P}\lra \Psi_P^*\left.\mathcal{L}_{\S_P}\right|_{U_P}\stackrel{Gl}{\lra} \M_{0,l}(P,\sigma).$$
Following \cite{CL},    we say that  $Gl_{\mathcal{Q}_P}$  is  \emph{admissible}. Moreover, if $U_P\subset\M_{\S_P}(P)$, we say that $Gl_{\mathcal{Q}_P}$ is of \emph{type-{\rm 1}},
 otherwise we say that it is of  \emph{type-{\rm 2}}. In particular, the gluing map $Gl_{\S_P}$  constructed directly from $\wt{U}_P$ is admissible and of type-1.

 Using gluing maps we introduce a topology basis on $\ov{\M}_{0,l}(P,\sigma)$ as follows:  an open neighbourhood of $(u,\j)\in \ov{\M}_{0,l}(P,\sigma)$ will be the image of some gluing map $Gl_{\S_P}$ previously constructed. Hence,  a neighbourhood is given by charts data of the type  $(\L_{\S_P,\epsilon_P,U_P},Gl_{\S_P})$.
 
  A standard argument shows that these charts are $C^0$ compatible \cite{CL}, \cite{MS}. This can be proved by comparing any admissible gluing map  $Gl_{\mathcal{Q}_P}$ arising from a chart data $\mathcal{Q}_P$ for  a proper subset $\wt{U}_P\subset\M_{\S_P}(P)$, with the type-1 gluing map  $Gl_{\S_P}$ on $\wt{U}_P$. Concretely one shows that  for small enough $\rho$, the map $(Gl_{\S_P})^{-1}Gl_{\mathcal{Q}_P}$ is close to the identity map, hence continuous. Thus, the moduli space  $\ov{\M}_{0,l}(P,\sigma)$ has the structure of an orbifold in the topology given by the gluing maps. 
  
  The smooth orbifold structure is  given by the two lemmas in the next subsection. In these two lemmas we prove more. Namely, we  construct smooth atlases on both  $\ov{\M}_{0,l}(P,\sigma)$ and $\ov{\M}_{0,l}(B,\sigma_B)$ compatibly with  $\ov{\pi}$.

 %Now let $\S_P$ such that $\S=\F_{m}(\S_P)$ and $\S_B=\pi(\S_P)$ are stable.  
%In other words,  we have diffeomorphisms:
%\begin{equation*}\Psi_P:=\phi_P\circ (1\times f^P): U_P\lra \wt{U}_P\hspa\text{and}\hspa \Psi_B:=\phi_B\circ (1\times f^B):U_B\lra \wt{U}_B
%\end{equation*}
%such that $$\pi\circ \Psi_P=\Psi_B\circ \pi.$$ From these datas we construct new gluing maps:
%$$Gl_{\mathcal{D}_P}:\left.\L_{\S}\right|_{\wt{U}_P}\lra \M_{0,l}(P,\sigma)\hspa\text{and}\hspa Gl_{\mathcal{D}_B}:\left.\L_{\S}\right|_{\wt{U}_B}\lra \M_{0,l}(B,\sigma_B),$$ defined by %precomposing with the given charts. 
%Precisely, given hypothesis  \eqref{hyp}, we have gluing maps:
%$$Gl^P:\left.\L_{\S}\right|_{U_P}\lra \M_{0,l}(P,\sigma)\hspa\text{and}\hspa Gl^B:\left.\L_{\S}\right|_{U_B}\lra \M_{0,l}(B,\sigma_B),$$
%$$Gl_{\Psi_P}:\left.\L_{\S}\right|_{U_P}\lra \M_{0,l}(P,\sigma,J_P)\hspa\text{et}\hspa Gl_{\Psi_B}:\left.\L_{\S}\right|_{U_B}\lra \M_{0,l}(B,\sigma_B,J_B),$$ 
%which are  compatible with $\pi$. Then set
%\begin{equation*}Gl_{\mathcal{D}_P}:=Gl^P\circ \Psi_P^{-1}\hspa\text{and}\hspa Gl_{\mathcal{D}_B}:=Gl^B\circ \Psi_B^{-1}.
%\end{equation*}
%Note that, by definition,
%$$\pi\circ Gl_{\mathcal{Q}_P}= Gl_{\mathcal{Q}_B}\circ \pi.$$
  
%  One can ask about the compatibility between  $Gl^P$ and $Gl_{\mathcal{Q}_P}$ (resp. $Gl^B$ and $Gl_{\mathcal{Q}_B}$).  In the next section, we actually show that these maps are compatible with respect to the $C^{\infty}$ structure of the top strata $\M_{0,l}(P,\sigma)$ and $\M_{0,l}(B,\sigma_B)$. 

\subsubsection{Structure of orbi-bundle}  Consider now the  fibration context.  We begin by the following observation. Let $\mathcal{Q}_P:=(U_P,Q^P)$ projecting onto  $\mathcal{Q}_B:=(U_B,Q^B)$ in the sense that  $\ov{\pi}_{\S_P}(U_P)=U_B$ and $Q^P$ is of the matrix form \ref{matrixQ}. Also, suppose that both pairs satisfy the assumption \eqref{hyp}, and that they generate charts datas, 
$(U_P,\phi_P,f^P)$ and  $(U_B,\phi_B,f^B)$, for some proper open subsets,  
$\wt{U}_P\subset\M_{\S_P}(P)$ and $\wt{U}_B\subset\M_{\S_B}(B)$ such that $\ov{\pi}_{\S_P}(\wt{U}_P)=\wt{U}_B$. Then, repeating the arguments in the gluing map section, we obtain that 
$$\ov{\pi}\circ Gl_{\mathcal{Q}_P}= Gl_{\mathcal{Q}_B}\circ \ov{\ov{\pi}_{\S_P}}.$$
This implies that the orbifold structures on  the compactified moduli spaces are defined compatibly
with  $\ov{\pi}$ which is continuous, open, and surjective, in the topology of the gluing maps. 
%Using th e gluing maps we introduce a topology basis on the compactified  moduli spaces as follows: an open neighbourhood of $(u_B,\j)$ in $\ov{\M}_{0,l}(B,\sigma_B)$ will, by definition, be the image of some gluing map $Gl_{\S_B}$ previously constructed. We define similarly an open neighbourhood of $(u,\j)$ in $\ov{\M}_{0,l}(P,\sigma)$ as being the image of some gluing map $Gl_{\S_P}$. A neighbourhood of a given point is given by charts of the type  $(\L_{\S_B,\epsilon_B,U_B},Gl_{\S_B})$ (resp. $(\L_{\S_P,\epsilon_P,U_P},Gl_{\S_P})$). Chen and Li have shown (lemma in \cite{CL}), that for small enough $\rho$, the map $(Gl^B)^{-1}Gl_{\mathcal{Q}_B}$ (resp. $(Gl^P)^{-1}Gl_{\mathcal{Q}_P}$) is 
%whenever $p>2$, there exists constants $C^P(\rho)$ and $C^B(\rho)$ converging to  $0$ when $\rho\lra0$ such that
%$$\|Gl^P(u,\j,\rho)-Gl_{\mathcal{D}_P}(u,\j,\rho)\|\leq C^P(\rho)$$ and $$ \|Gl^B(\pi(u),\j,\rho)-Gl_{\mathcal{D}_B}(\pi(u),\j,\rho)\|\leq C^B(\rho).$$
We now construct smooth atlases on both  $\ov{\M}_{0,l}(P,\sigma)$ and $\ov{\M}_{0,l}(B,\sigma_B)$ compatibly with  $\ov{\pi}$. In order to do so we introduce stratum-coverings.

%\begin{rem} Les atlas consid\'er\'es sont particuliers. (pas de compatibilit\'e entre les atlas)
%\end{rem}

\begin{defn} A strata-covering of $\ov{\M}_{0,l}(P,\sigma)$ consists in pairs $(U_P,\epsilon_{\S_P})$ for each stratum data $\S_P$, such that:
\begin{enumerate}[$\bullet$]
\item $U_P$ is a proper open subset of  $\M_{\S_P}(P)$,
\item there exists a well-defined gluing map $Gl_{\S_P}$ with domain  $\L_{\S_P,\epsilon_{\S_P},U_P}$,
\item letting $W_{\S_P}$ be the image $Gl_{\S_P}(\L_{\S_P,\epsilon_{\S_P},U_P})$,
we have that for any two (effective) stratum datas $\S_P$ and $\S_P'$:
$$W_{\S_P}\cap W_{\S_P'}\neq\emptyset\hspa\text{iff}\hspa \S_P\prec\S_P',\hs\text{or}\hs \S_P'\prec\S_P,$$
\item the family $\{W_{\S_P}\}_{\mathcal{D}_{0,l}^{\sigma,J_P}}$ yields an open covering of $\ov{\M}_{0,l}(P,\sigma)$.
\end{enumerate}
\end{defn}

%Analogously we define a stratum covering for $\ov{\M}_{0,l}(B,\sigma_B,J_B)$.

\begin{lem}\label{stratumcovering} There exists strata-coverings $(U_{P},\epsilon_{\S_P})$ and $(U_{B},\epsilon_{\S_B})$ for   $\ov{\M}_{0,l}(P,\sigma)$ and $\ov{\M}_{0,l}(B,\sigma_B)$, such that  $\ov{\pi}_{\S_P}(U_{P})=U_{B}$.
\end{lem} 

\proof\hs  The proof is an induction on the stratum datas in  $\mathcal{D}^B:=\mathcal{D}_{0,l}^{\sigma_B,J_B}$ and $\mathcal{D}^P:=\mathcal{D}_{0,l}^{\sigma,J_P}$. 
%We start with the lower strata in $\ov{\M}_{0,l}(B,\sigma_B)$ and show inductively that the first three  points of the  
%definition are verifyied for every stratum in  $\ov{\M}_{0,l}(P,\sigma)$ projecting over the former lower strata. We then iterate the procedure 
%using the partial order on $\mathcal{D}^B$. 
%Nous montrons comment l'it\'eration proc\`ede pour\`a pr\'esent pour
Let  $\S_{B,0}$ be the set of lowest strata in $\mathcal{D}^B$. For $\S_B$ in $\S_{B,0}$ set $U_B=\M_{\S_B}(B)$. Since $\S_B$ is minimal  $U_B$ is compact and there exists $\epsilon_{\S_B}$ and a  gluing map  $Gl_{\S_B}$ defined on  the restriction of $\L^*_{\S_B,\epsilon_{\S_B}}$ to $U_{B}$.
Furthermore,  the minimal strata are isolated and for each $\S_B\in\S_{B,0}$ we can choose a small enough  $\epsilon_B$ so that the resulting   gluing neighbourhoods never intersect.
Now let  $\S_{P,0,0,\S_B}$ be the set of lowest strata in $\mathcal{D}^P\cap S^{-1}_{\pi_*}(\S_B)$, where  $\S_B$ is minimal.  For any $\S_P\in \S_{P,0,0,\S_B}$ set $U_P=\M_{\S_P}(P)$. 
 Argumenting as above,  there is  $\epsilon_{\S_P}$ and $Gl_{\S_P}$, with domain  $\L^*_{\S_P,\epsilon_{\S_P},U_P}$,  such that
 $\ov{\pi}\circ Gl_{\S_P}=Gl_{\S_B}\circ \ov{\ov{\pi}_{\S_P}}$.
 Once again, we can choose the  $\epsilon_{\S_P}$  such that $W_{\S_P}\cap W_{\S_P'}=\emptyset$, 
for any two strata in $\S_{P,0,0,\S_B}$. 
 
Define inductively  $\S_{B,k}$ as being the set of minimal strata in $\mathcal{D}^B\bs\S_{B,k-1}$, and  $\S_{P,k,m_k}$ as being the set of minimal strata in $\mathcal{D}^P\cap S^{-1}_{\pi_*}(\S_{B,k})\bs \S_{P,k,m_k-1}$. Suppose that each pairs $(U_B,\epsilon_{\S_B})$ and $(U_P,\epsilon_{\S_P})$, for $\S_B\in \S_{B,n}$ and $\S_P\in \S_{P,n, m_n}$ with $n\leq k-1$ and $m_n\leq m_k-1$, have been chosen so that the induction holds. Set
$$W_{\S_B,\S'_B}:= Gl_{\S_B,\S'_B}(\L_{\S_B,\S'_B,\epsilon_{\S_B},U_B}).$$
Then, for  $\S_B'\in \S_{B,k}$ we can choose a proper open subset $U'_B$ such that 
$\{W_{\S_B,\S'_B}|\S_B\prec\S_B'\}\cup U'_B$, 
 is a covering for $\M_{\S_B'}(B)$. 
 Furthermore, there is  $\epsilon_{\S_B'}$ and a gluing map 
 $$Gl_{\S_B'}:\L^*_{\S_B',\epsilon_{\S_B'},U'_B}\sra \M_{0,l}(B,\sigma_B),$$
 and we can make sure that for all $\S_B'\in\S_{B,k}$ and $\S_B\in\bigsqcup^k_{i=0}\S_{B,i},$ the intersection  $W_{\S_B'}\cap W_{\S_B}$ is empty  unless
 $\S_B\prec\S_B'$ (by choosing smaller $\epsilon_{\S_B'}$ and $\epsilon_{\S_B}$ if necessary).

 Since the  gluings commute with the projection, it suffices to fix  $\S_B\in \S_{B,k}$,  and to apply the arguments given for $k=0$  to the elements of  $\S_{P,k,m_k,\S_B}:=\S_{P,k,m_k}\cap S^{-1}_{\pi_*}(\S_B)$. Set 
$$ W_{\S_P,\S'_P}:= Gl_{\S_P,\S'_P}(\L_{\S_P,\S'_P,\epsilon_{\S_P},U_P}),$$ 
where $\S'_P$ projects onto $\S_B$. Note that $W_{\S_P,\S'_P}$ also projects onto $W_{\pi(\S_B),\pi(\S'_B)}$. 
We can choose  $U_{\S_P'}$ such that $\{W_{\S_P,\S'_P}|\S_P\prec\S_P'\}\cup U'_P$, covers
$\M_{\S_P'}(P)$, and for a well chosen $\epsilon_{\S_P'}$ we have a map  $Gl_{\S_P'}$ which image does not intersect the neighbourhoods obtained so far  unless it comes from a stratum $\S_P$ such that $\S_P\prec\S_P'$. 
\qed\\

%\begin{rem}
%For fixed stable stratum data  $\S_B$ set:
%$$\ov{\M}^v_{\S_B}(P):=\bigsqcup_{\{\S_P\,\,\text{stable}\,\,|\F_{\pi}(\S_P)=\S_B\}}\M_{\S_P}(P).$$
%%By the construction above, this stratified space {\rm(}dont chaque composante est une orbi-vari\'et\'e par {\rm (\ref{strongreg})}{\rm)}, est recouvert par ses strates dans le sens donn\'e par la d\'efinition ci-haut, au-dessus de $U_B$.
%\end{rem}

A strata-covering gives an atlas. If the transition functions were to be smooth we would directly have a smooth orbifold structure on the considered moduli spaces. However, this may be hard to show and even not true in full generality. Instead, we show (cf \cite{CL}) that for each stratum data there are charts $Gl_{\S_B}$  such that the composition $Gl_{\S_B}\circ Gl_{\S_B,\S_B'}^{-1}$ is smooth for every 
 $\S_B'\prec\S_B$, which provides, not canonically, a smooth atlas. 

%Nous d\'esirons \`a pr\'esent montrer que ces recouvrements, qui forment un atlas des espaces de modules compactifi\'es, sont compatibles de fa\c{c}on lisse. Autrement dit que
%$$Gl_{\S_B}\circ Gl_{\S_B,\S'_B}^{-1}\et Gl_{\S_P}\circ Gl_{\S_P,\S'_P}^{-1}$$
%sont $C^{\infty}$ pour toutes cartes apparaissant dans les atlas correspondants. Nous avons:

\begin{lem} There are strata-coverings  $(U_P,\epsilon_{\S_P})$ and $(U_B,\epsilon_{\S_B})$, and gluing maps, $Gl_{\S_P}$ and $Gl_{\S_B}$, compatible with $\ov{\pi}_{\S_P}$, such that for every stratum datas $\S_P$ and $\S_B$,  the maps $Gl_{\S_P}$ and $Gl_{\S_B}$ coincide with any other gluing maps  
$$Gl'_{\S_P}=Gl_{\S'_P}\circ Gl_{\S'_P,\S_P}^{-1}\and  Gl'_{\S_B}=Gl_{\S'_B}\circ Gl_{\S'_B,\S_B}^{-1},$$ where $\S_P'\prec\S_P$ and $\S_B'\prec\S_B$.
\end{lem}

\proof\hs  The proof is again by induction. Let $\S_{B,k}$ and $\S_{P,k,m_k}$ as in lemma \ref{stratumcovering}. We see that the result holds for $\S_{B,0}$ and $\S_{P,0,0}$. 
Suppose it is true for all  $\S_B\in \S_{B,n}$ and $\S_P\in \S_{P,n, m_n}$ 
such that $n\leq k-1$ and $m_n\leq m_k-1$. 
Let  $\S_B\in \S_{B,k}$ and set $W_{\S_B}:=\cup_{\S_B'\prec\S_B} W_{\S'_B,\S_B}$.
Let  $Gl'_{\S_B}(\S'_B)$ be the gluing map induced by $\S'_B\prec\S_B$. 
Recall that this map is defined above $W_{\S'_B,\S_B}$. 
We must show that
\begineq \label{egalglu}
Gl'_{\S_B}(\S'_B)=Gl'_{\S_B}(\S''_B),
\endeq
on $W_{\S'_B,\S_B}\cap W_{\S''_B,\S_B}$. 
But this latter intersection is non-empty if and only if  $\S''_B\prec\S'_B$. But from the induction, $Gl'_{\S'_B}(\S''_B)=Gl_{\S'_B}$ on $W_{\S''_B,\S'_B}\cap U'_B$. Thus 
\begineqn 
Gl'_{\S_B}(\S''_B)&=&Gl_{\S''_B}\circ Gl_{\S''_B,\S_B}^{-1}\nonumber\\
&=&Gl_{\S''_B}\circ Gl_{\S''_B,\S'_B}^{-1}\circ Gl_{\S''_B,\S'_B} \circ Gl_{\S''_B,\S_B}^{-1}\nonumber\\
&=&Gl_{\S'_B}\circ Gl_{\S''_B,\S'_B} \circ Gl_{\S''_B,\S_B}^{-1},\nonumber
\endeqn
giving (\ref{egalglu}). 
 As a result, we obtain a gluing map $Gl'_{\S_B}$ defined on $W_{\S_B}$. 
Now, given a gluing map  $Gl''_{\S_B}$ on $U_B$, 
we derive a third map $Gl_{\S_B}$, which is obtained as an interpolation between $Gl'_{\S_B}$ and $Gl''_{\S_B}$ using a cut-off function. This ends the induction for $\ov{\M}_{0,l}(B,\sigma_B)$. 

We explain in details how to interpolate   $Gl'_{\S_B}$ and $Gl''_{\S_B}$. By definition  $Gl'_{\S_B}(\S_B')$ is of type-2 with domain:
% Let $\S_B'$  such that  $\S_B'\prec\S_B$ and consider $Gl'_{\S_B}(\S_B')$ on $W_{\S'_B,\S_B}$, the corresponding gluing map. By definition
$$W_{\S'_B,\S_B}= Gl_{\S'_B,\S_B}(\L_{\S'_B,\S_B,\epsilon_{\S'_B},U'_B})\equiv W_{\S'_B,\S_B}(\epsilon_{\S_B'}).$$
For  $Gl_{\S'_B,\S_B}$ admissible, the associated chart data is a triple $(V:=\prgl(\L_{\S_B',\S_B,\epsilon_{\S'_B}}),\phi,F)$ where:
$$F:V\sra \M_{\S_B}(B),\hspa (u_{B,\rho},\j_{\rho}) \mapsto \exp_{u_{B,\rho}}Q^B_{u_{B,\rho}}f_{\S_B}(u_B,\j,\rho).$$
%for $(u_{B,\rho},\j_{\rho})= \prgl (u_B,\j,\rho)$ an element of the bundle $\L_{\S'_B,\S_B,\epsilon_{\S'_B}}$ restricted to $U'_B$.
Consider  a function
$$\nu:\L_{\F_{B}(\S'_B),\F_{B}(\S_B),\epsilon_{\S'_B}}\sra \R, \hspa \nu(\j,\rho)\mapsto \begin{cases}
0 &\text{si $|\rho|\leq 0.5\epsilon_{\S_B'}$}\\
1 &\text{si $|\rho|\geq 0.75\epsilon_{\S_B'}$.}
\end{cases}
$$
We now ``glue" the domains of the chart datas for  $Gl''_{\S_B}$ and $Gl'_{\S_B}(\S_B')$, i.e we glue $U_B$ and $V$. The new domain
$$V':=\text{Im}\left(\exp_{u_{B,\rho}}\left(\nu(\j,\rho)Q^B_{u_{B,\rho}}f_{\S_B}(u_B,\j,\rho)\right)\right),$$
 coincides with $V$ for  $|\rho|\leq 0.5\epsilon_{\S'_B}$, and with $U_{\S'_{B}}$ when $|\rho|\geq 0.75\epsilon_{\S'_B}$.
We can make sure that  the pair,  $(V',\{Q^B_{u_B}|u_B\in V'\})$, satisfies hypothesis \eqref{hyp} since  $V'$ is a uniform deformation between $V$ and $U_B$ ($\nu$ does not depend on  $u_B$).
%Indeed $(V,\left.Q^B\right|_{V})$ already verifies those hypothesis, que sur la partie dans $U_B$ ces hypoth\`eses sont directement satisfaites aussi,  $\nu$ does not depend on  $u_B$, ce qui nous donne une d\'eformation uniforme. 
Denote by  $Gl_{\S_B}$ the gluing map arising from the pair $(V',\{Q^B_{u_B}|u_B\in V'\})$. Then,  $Gl_{\S_B}$ coincides with $Gl'_{\S_B}(\S_B')$ on $V\cap V'= W_{\S'_B,\S_B}(0.5\epsilon_{\S_B'})$, and with $Gl''_{\S_B}$ on $U_B\cap V'$. In particular, this map extends to  $U_B \bs W_{\S'_B,\S_B}(\epsilon_{\S_B'})$. Regarding $\ov{\M}_{0,l}(P,\sigma)$ we proceed similarly. 
Fixing  $\S_B\in\S_{B,k}$, and applying the same arguments as above to the elements 
$\S_P\in  \S_{P,k, m_k,\S_B}$, ends the proof.
\qed\\

\subsection{The product formula revisited}
In what follows, we will assume hypothesis \eqref{strongreg} is verified for a given fibered structure  $J_P=(J_B,J,H)$ and classes $\sigma$, $\sigma_B\neq0$ such that  $\pi_*\sigma=\sigma_B$. Hence, $\ov{\pi}$ is fibration when restricted to the top stratum of $\ov{\M}_{0,l}(B,\sigma_B)$. In addition,  the lower  strata in $\ov{\M}_{0,l}(B,\sigma_B)$
%\bs\M^*_{0,l}(B,\sigma_B)\and \bs\M^{**}_{0,l}(P,\sigma) 
and $\ov{\M}_{0,l}(P,\sigma)$
are of codimension at least two. The same applies for the lower strata in  $\ov{\pi}^{-1}(u_B,\mathbf{x})$, for every $(u_B,\mathbf{x})\in\ov{\MMC}_{0,l}(B,\sigma_B)$.

\begin{prop} 
%Suppose there is a fibered structure  $J_P=(J_B,H,J)$ and  classes $\sigma$ and $\sigma_B$, with respect to which hypothesis \ref{strongreg} is satisfied, %that $\coker D^v$ is zero on every stratum of $\ov{\MMC}^v_{0,k}(P,\sigma;J_P)$ for some regular triple $(J_b,J_B,H)\in \mathcal{P}_{reg}$. 
Under the assumptions above the product formula is obtained using integration over the fibers of  $\ov{\pi}$. 
\end{prop}

\proof\hs  Set $$\MMC^P_{0,l}:=\MMC^{**}_{0,l}(P,\sigma)\hsp\text{and}\hsp \MMC^B_{0,l}:=\MMC^{*}_{0,l}(B,\sigma_B),$$ 
and consider the evaluation maps,
$$ev^P:\ov{\MMC}^P_{0,l} \sra P^l\,\,\, ,\hspa ev^B:\ov{\MMC}^B_{0,l}\sra B^l\,\,\, ,\hspa ev_{(u_B,\mathbf{x})}:\ov{\pi}^{-1}(u_B,\mathbf{x})\sra F^l, $$
where $(u_B,\mathbf{x})\in\MMC^B_{0,l}$. These  are  smooth with respect to the gluing topology.

Now, let $c_i^P$, $c_i^B$  and $c_i^F$ ($i=1,...,m$) be classes  as in condition  \eqref{cond11}. 
 We can represent these classes by submanifolds after multiplying them by  well chosen integers, if necessary.  
  %Since the inferior strata of the moduli spaces  are of codimension at least two, the preimages under  $ev^P$, $ev^B$, and $ev_{(u_B,\mathbf{x})}$ of the induced product  submanifolds are generically  given by a finite number  points in $\M^P_{0,l}$, and $\M^B_{0,l}$. 
Now represent the Poincar\'e duals of  $c_i^P$, $c_i^B$  and $c_i^F$, by differential forms,  $\alpha_i^P$, $\alpha_i^B$, and $\alpha_i^F$, compactly supported in a small enough tubular neighbourhoods around the submanifolds. Then the  pull-backs of these forms along $ev^P$, $ev^B$ or $ev_{(u_B,\mathbf{x})}$, are also compactly supported. 
  By definition, $$\alpha_i^P=\pi^*(\alpha_i^B), \hspa i=m+1,...,l.$$
 Moreover, since $c_i^B=pt$:
  $$\alpha_i^B=\text{vol}(B),\hspa i=1,...,m.$$ 
%$PD(c_i^P)=\pi^*(PD(c^B_i))$ for $i=k+1,...,l$. 
If  $\mathcal{N}_i$ denotes a tubular neighbourhood around the fiber above $c_i^B=pt$, for $i=1,...,m$, and if  $\rho_i:\mathcal{N}_i\sra F$ is a deformation retract associated to this normal neighbourhood,  we obtain that: 
\begin{equation}\label{expression}\alpha^P_i=\pi^*\alpha_i^B\wedge \rho_i^*\alpha_i^F=\pi^*\text{vol}(B)\wedge  \rho_i^*\alpha_i^F.\end{equation}
%$$PD(c^P_i)=PD(F_{c_i^B})\wedge \rho_i^*PD(c_i^F)=\pi^*PD(c_i^B)\wedge  \rho_i^*PD(c_i^F),$$
%(Here the compact support of $\rho_i^*PD_F(c_i^F)$ should lie in the compact support of $PD_P(F_{c_i^B})$). 
Here we have to make sure that the support of $\pi^*\text{vol}(B)$ does not strictly contain the support of  $\rho_i^*\alpha_i^F$, but this can be realized by decreasing the support of    $\text{vol}(B)$ if necessary.  By definition, (see \cite{R}) we have that:
\begin{equation*}\label{GWP}\la c_1^P,...,c_l^P\ra_{0,l,\sigma}^P= \int_{\MMC^P_{0,l}}(ev^P)^*\left( \bigwedge_{i=1}^l \alpha_i^P\right)= \int_{\MMC^B_{0,l}} \ov{\pi}_* (ev^P)^* \left(\bigwedge_{i=1}^l \alpha_i^P\right), 
\end{equation*}
where $\ov{\pi}_*$  stands for integration along the fibers. Using \eqref{expression}, the equation above then equals:
$$ \int_{\MMC^B_{0,l}} \ov{\pi}_* (ev^P)^* \left(\bigwedge^k_{i=1}( \pi^*\alpha^B_i\wedge  \rho_i^*\alpha_i^F))\wedge\bigwedge^l_{i=k+1}\pi^*\alpha^B_i \right). $$

Let  $ev^P_i$ be the  projection of  $ev^P$ on the $i^{\text{th}}$ component of $P^l$ and define $ev^B_i$ and $ev_{(u_B,\mathbf{x}),i}$ similarly. For  $i=1,...,l$, we have $\pi\circ ev^P_i= ev^B_i\circ \ov{\pi}$ so that: 
\begin{equation*}\la c_1^P,...,c_l^P\ra_{0,l,\sigma}^P=(-1)^{\alpha} \int_{\MMC^B_{0,l}}\ov{\pi}_* \left(\bigwedge^l_{i=1} \pi^*(ev_i^B)^*\alpha^B_{i}\wedge \bigwedge^k_{i=1}(ev_i^P)^* \rho_i^*\alpha_i^F \right),
\end{equation*}
where $\alpha=\sum^{l}_{i=k+1} \text{deg}\alpha^B_i\sum^k_{i=1}\text{deg}\alpha_i^F$ is odd (by a simple dimension argument). 
Furthermore, since $\pi_* (a \wedge \pi^* b)=(\pi_* a)\wedge b$, for any form $a$ and $b$ we finally have:
%L'unique possibilit\'e pour laquelle l\''equation ci-dessus est non n\'ecessairement nulle est lorsque $$\sum^{l}_{i=k+1} \text{deg}\alpha^B_i=\sum^{l}_{i=1} \text{deg}\alpha^B_i=2n_B+2c_1(\sigma_B)+2l-6$$ qui est pair.
\begin{equation*} \la c_1^P,...,c_l^P\ra_{0,l,\sigma}^P= \int_{\MMC^B_{0,l}} \left((ev^B)^*\bigwedge^l_{i=1}\alpha^B_{i}\right)\wedge \ov{\pi}_* \left(\bigwedge^k_{i=1}(ev_i^P)^* \rho_i^*\alpha_i^F \right),
\end{equation*}
where the term  involving integration over the fibers of $\pi$
% $$\pi_* \left(\bigwedge^k_{i=1}(ev_i^P)^* \rho_i^*\alpha_i^F \right)$$
must be a  function  $\psi$ on $\MMC^B_{0,l}$  given by (cf  \eqref{evdiagram}):
$$\psi(u_B,\mathbf{x})=\int_{\ov{\pi}^{-1}(u_B,\mathbf{x})}ev^*_{(u_B,\mathbf{x})}\left(\bigwedge_{i=1}^k\alpha^F_i\wedge\bigwedge_{i=k+1}^l\text{vol}(F)\right)=n_{\alpha}$$
%If  $C$  is the image of a $J_B$-holomorphic map  $u_B$, then $\psi(u_B,\mathbf{x})$ coincides with
%$$\sum_j GW^{P_C}_{0,l}(\iota_F^{P_C}(c^F_1),...,\iota_F^{P_C}(c^F_l),\sigma_j),$$
where $n_{\alpha}$ is given in \eqref{relationfibre}.  Hence, $\psi$  does not depend on $(u_B,\mathbf{x})$ and we can withdraw it from the integral. This ends the proof.  
%Then the proposition follows since: 
%$$\la c_1^B,...,c_l^B\ra_{0,l,\sigma_B}^B:= \int_{\MMC^B_{0,l}} \left((ev^B)^*\bigwedge^l_{i=1}\alpha^B_{i}\right).$$
\qed

\subsection{An example of non-triviality of the induced fibration} %Consider the rank 2 holomorphic vector bundle $\pi:V:=\mathcal{O}_{\CP^2}(1)\oplus \C\lra \CP^2$, and let $P$ be its projectivization with projection to $\CP^2$ denoted by $\pi$ again. Then $P$ together with $\pi$ is a $(\CP^1,\om_{FS})$ Hamiltonian fibration over $( \CP^2,\om_{FS})$.
Consider the  Hamiltonian fibration $\pi:P\sra \CP^2$ with fiber $(F:=\CP^1,\om_{FS})$,
%$$(\CP^1,\om_{FS})\into \mathbb{P}(\mathcal{O}_{\CP^2}(1)\oplus \C)\stackrel{\pi}{\lra}( \CP^2,\om_{FS}).$$
where  $P$ is the projectivization of the rank 2 holomorphic vector bundle, $$\pi:V:=\mathcal{O}_{\CP^2}(1)\oplus \C\sra \CP^2.$$  Let $J_0$ be the standard complex structure on  $\CP^2$ which is compatible with $\om_{FS}$. Let $J_P$ be the integrable structure on $P$ induced by $J_0$, the structure of complex fibration on $V$ and the holomorphic hermitian connection on $V$ (inducing the coupling form here.).  

Let $h\in H^2(P,\Z)$ denote the pull-back under $\pi$ of the positive generator in  $H^2(\CP^2,\Z)$, Poincar\'e dual to the class $L\in H_2(\CP^2,\Z)$ of 
a line in $\CP^2$. Let also   $\xi\in H^2(P;\Z)$ be the first Chern class of  the dual of the tautological bundle over $P$.  It is standard   that:
$$ H^*(P;\Z)\cong \frac{\Z[h,\xi]}{\{h^3=0,\hs \xi^2+h\xi=0\}}.$$
Then,  $H_2(P;\Z)$ is generated by the duals of the classes  $h^2$ and $h\xi$. Let  $L_0\in H_2(P,\Z)$  denote the  Poincar\'e dual of  $h^2+h\xi$. If $\pi_*$ represents integration along the fibers of $\pi$, then $\pi_*L_0=L$.  
%since $\pi_*(h\cup\xi)=h$ and $\pi_*(h^2)=0$.

 The map $\pi$ induces the projection
$$\ov{\pi}:\ov{\M}(P,L_0,J_P)\sra \ov{\M}(\CP^2,L,J_0)\cong(\CP^{2})^*.$$
The source moduli space is made of two strata: $\S_0$ the top stratum of simple maps, and $\S_1$ the stratum consisting of stable maps having two components, one being a $\pi_*$-stable component  representing the class of the exceptional divisor, $PD(h\xi)$, and the other being a $\pi_*$-unstable component representing the class $PD(h^2)=[F]$. Observe that the second stratum contains only irreducible elements. 

\begin{lem}
The fibered complex structure $J_P$ is regular and parametric for the class $L_0$.
\end{lem}

\proof   First, since $L$ is $J_0$-indecomposable, $J_0$ belongs to $ \J_{irr}(L)$.  Consider $u\in\S_0$. From the exact sequence \eqref{splitting}, we only have to verify that $\coker D^v_{u}=0$,  but 
$$\coker D^v_{u}= H^{0,1}(S^2,u^*TP^v)=H^{0,1}(S^2,\O_{\CP^2}(1))\cong H^0(S^2,\O(-3))=0.$$ 
Now consider $(u_1,u_2, y_{12},y_{21})\in \S_1$, where $u_1$ denotes the $\pi_*$-stable component, and  $u_2$  denotes the $\pi_*$-unstable component. Again,  we have to check that $\coker D^v_{u_1}$ and $\coker D^v_{u_2}$ vanish. Moreover, we have to show that for every line  $\ell\in (\CP^2)^*$ the edge evaluation map
$$ev:\ov{\pi}^{-1}(\ell)\cap\S_1 \sra \pi^{-1}(b)\times  \pi^{-1}(b),\hspa (u_1,u_2, y_{12},y_{21})\mapsto (u_1(y_{12}), u_1(y_{21})),$$
is transverse to the diagonal $\Delta_{ \pi^{-1}(b)}$, where $b:=\pi(\ov{b})$ and $\ov{b}$ denotes  the unique intersection point  between  $u_1$ and $u_2$.  Again, for every holomorphic map representing $[F]$,  $\coker D^v_{u_2}$ can be identified to $H^{0,1}(S^2,T\CP^1)$ which vanishes,  and similarly,
$$\coker D^v_{u_1}= H^{0,1}(S^2,u_1^*TP^v)=H^{0,1}(S^2,\O(-1))=0.$$ Finally,   transversality of  $ev$ follows since for every  $X_0\in T_{\ov{b}}\pi^{-1}(b)$ there exists a  holomorphic vector field on $\CP^1\cong\pi^{-1}(b)$ with value  $X_0$ at $\ov{b}$. 
\qed\\

Let $\wt{\CP}^2$ denote the blow-up of $\CP^2$ at one point. Then  the restriction $\left.P\right|_{\ell}$ of $P$ to a line $\ell$ in  $\CP^2$ is identified to $\wt{\CP}^2$. The fibers of $\ov{\pi}$ are then given by  $\ov{\M}(\wt{\CP}^2,L_0,J)$ where $J$ is the complex structure associated to the Hirzebruch surface $\mathbb{P}(\O_{\CP^1}(1)\oplus \C)$, and $L_0$ is the  class represented by the zero section. Furthermore, 
  %the set holomorphic sections of the sheave $\O_{\CP^1}(1)$ which is given by
 $$\M^*(\wt{\CP}^2,L_0,J)\cong H^0(\CP^1,\O_{\CP^1}(1))\cong\C\LA u,v\RA,$$
where $[u:v]$ stands for  the homogeneous coordinates on $\CP^1$. The Gromov closure then corresponds to adding the line at infinity so that $\ov{\M}(\wt{\CP}^2,L_0,J)=\CP^2$.  

Next, we show that $\ov{\pi}$ is  non-trivial, more precisely that  it is a $\CP^2$-fibration  obtained as the  projectivization of a non-trivial rank two holomorphic bundle over  $(\CP^2)^*$. Consider the incidence variety:
$$W:=\left\{(p,\ell)\in \C P^2 \times (\CP^2)^*\,\,|\,\, p\in \ell\right\},$$
and let $\pi_1$, $\pi_2$ denote the projections on the first and second factors. Note that 
%$p_2:W\lra  (\CP^2)^*$ is the tautological bundle over $(\CP^2)^*$, i.e 
$(W,\pi_2)$ is the projectivization of $\O_{(\CP^2)^*}(-1)\oplus \C$. Consider the direct image sheave over $(\CP^2)^*$,
$$\mathcal{R}:=\pi_{2_*}\pi_1^*\O_{\CP^2}(1).$$  
which germ at $\ell\in (\CP^2)^*$ is given by
$$H^0(\pi^{-1}_2(\ell),\left.\O_{\CP^1}(1)\right|_{\pi^{-1}_2(\ell)})\cong H^0(\CP^1,\O_{\CP^1}(1)),$$
Hence $\mathcal{R}=\S_0$. 
Let $\mathcal{D}$ be a line in $(\CP^2)^*$. Then $\pi^{-1}_2(\D)$ can be identified to  $\wt{\CP}^2$, where the blown-up point is given by the intersection of all the lines generated by  $\mathcal{D}$. 
In fact, 
  $$\pi_2^{-1}(\mathcal{D})\cong \frac{D^+\times \CP^1\sqcup D^-\times \CP^1}{(\lambda,[u:v])\sim(\lambda^{-1},[\lambda u:v]),\hsp\text{$\lambda\neq 0$}},$$
 where $D^+$ and $D^-$ respectively denote the complements of $[0:1]$ and $[1:0]$ in $\CP^1$. Hence, the restriction of $\mathcal{R}$ to $\mathcal{D}$ is the direct sum of  $\O_{\CP^1}(-1)$ with  $ \C$, 
%$$\left.\mathcal{R}\right|_{\mathcal{D}}=\frac{D^+\times \C^2\bigsqcup D^+\times  \C^2}{(\lambda,( u,  v))\sim(\lambda^{-1},(\lambda u, v))\hsp\text{for $\lambda\neq 0$}},$$ 
% $$\left.\mathcal{R}\right|_{\mathcal{D}}=\O_{\CP^1}(-1)\oplus \C.$$
so that $\pi$ is non-trivial. 

As an example, we compute 
$$\la pt,pt\ra_{0,2,L_0}^{P}.$$ 
Since   the homology of the fiber injects in the  homology of the total space, the product formula simplifies to give:
%fondamentale $[\CP^1]$ est donn\'ee par $PD(h^2)\neq0$. Alors la formule produit se simplifie et nous avons:
\begin{equation*}\label{prodexample}\la pt,pt\ra_{0,2,L_0}^{P}=\la pt,pt\ra^{\CP^2}_{0,2,L}\cdot \la pt,pt \ra_{0,2,B}^{P_C}
\end{equation*}
where $C$ is the image of a holomorphic map in $\CP^2$ representing a line passing through two points, and $B$ is the Poincar\'e dual of the restriction of   $h^2+h\xi$ to $P_C$. This latter class corresponds to the sum of the fiber and the exceptional divisor in the blow up. It follows that both members of the above product give 1, hence
%There is only one such line. 
%Furthermore, the Poincare dual of the restriction to $\PC=\wt{\CP}^2$ of the class $h^2+h\xi$, actually corresponds to the sum of the exceptional divisor with the fiber in the blow up, so that the right member in the product \eqref{prodexample} gives 1, and we conclude that
$$\la pt,pt\ra_{0,2,L_0}^P=1.$$

\bibliographystyle{plain}
\bibliography{bibarticle}

\nocite*

\end{document}